\documentclass[a4paper,10pt]{amsart}
\usepackage[T1]{fontenc}
\usepackage{amssymb}

\newcommand{\vsp}{\vspace*{0.5cm}}
\newcommand{\D}{\mathcal{D}'}

\newcommand{\E}{\mathcal{E}'}
\newcommand{\e}{\mathcal{E}}
\newcommand{\U}{\mathcal{U}}
\newcommand{\A}{\mathcal{A}}

\newcommand{\R}{\mathbf{R}}
\newcommand{\B}{\mathcal{B}}

\newcommand{\s}{\mathcal{S}}

\newcommand{\F}{\mathcal{F}}

\begin{document}
\begin{flushleft}

\title{On partially hypoelliptic operators. \\ Part I: Differential operators}
\author{T. Dahn \\ \emph{Lund University}}

\begin{abstract}
This article gives a fundamental discussion on variable coefficients, self-adjoint, formally partially hypoelliptic differential operators. A generalization of the results to pseudo differential operators is given in a following article in ArXiv. Close to (\cite{Ni_72}), we give a construction and estimates of a fundamental solution to the operator in a suitable topology. We further give estimates of the corresponding spectral kernel.
\end{abstract}
\maketitle

\section{Preliminaries}
\label{sec:Prel}
For elliptic operators $P(D)$ with constant coefficients, we have
the following  result (\cite{Mal}):
\newtheorem{Weyl}{Proposition }[section]
\begin{Weyl}[Weyl's lemma]
In order for a distribution $\varphi$ defined in an open set to be an
infinitely differentiable function, it is necessary and sufficient that
$P(D)\varphi$ is an infinitely differentiable function in this open set.
\end{Weyl}

The hypoelliptic operators can be characterized as the class of
operators, for which this proposition holds (see Def. \ref{HE}). We
will refer to this proposition, as Weyl's criterion for
hypoellipticity. For these operators, the fundamental solutions are infinitely differentiable outside the
origin.

\vsp

The singular support of a distribution $u \in \D($ $\Omega)$,
$\mbox{sing supp u}$, is
defined as the set of points in $\Omega$ which do not have a
neighborhood in which $u$ is an infinitely differentiable function.
Note that hypoellipticity in general is very sensitive to change of topology,
however we will consider,

\newtheorem{HE}[Weyl]{Definition }
\begin{HE}[Hypoelliptic constant coefficients operator] \label{HE}
  An operator $P(D)$ is said to be hypoelliptic if and only if
  for every open set $\Omega \subset \R^{\textit{n}}$ and $\varphi \in
  \D($ $\Omega)$ we have
  $$\mbox{sing supp } P(D)\varphi = \mbox{sing supp } \varphi$$
\end{HE}
We use the term homogeneously hypoelliptic for an operator such that $P(D)u=0$ in $\D($ $\Omega)$
implies $u \in C^{\infty}(\Omega)$. There are several analogues to hypoellipticity that will be dealt with and
several condition on an operator which are equivalent with
hypoellipticity as in Weyl's criterion (see \cite{Ho_LPDO} sec. 11.1), and we will just mention one of these:
\newtheorem{slowosc}[Weyl]{Proposition }
\begin{slowosc}[\cite{Mal},Th.II.1.2] \label{so}
  For an operator $P(D)$ to be hypoelliptic, it is necessary and
  sufficient that the polynomial $P(\xi)$ does not vanish outside a
  compact set and that for every multi-index $\alpha$ such that $\alpha \neq 0$
  we have
  $$ \lim_{\xi \rightarrow \infty} \frac{P^{(\alpha)}(\xi)}{P(\xi)}=0.$$
\end{slowosc}
A consequence of this proposition is that the polynomials corresponding
to hypoelliptic operators, have the property of slow oscillation, that
is
$$ \frac{P(\xi + \eta)}{P(\xi)}= \sum_{\alpha,\beta} c_{\alpha,\beta}
\frac{P^{(\alpha)}(\xi)}{P(\xi)}P^{(\beta)}(\eta) \rightarrow 1 \quad \xi
\rightarrow \infty$$
We will also frequently use for the following proposition:
\newtheorem{Mal}[Weyl]{Proposition }
\begin{Mal}[\cite{Mal},Prop.II.1.5] \label{Mal}
  If $P(D)$ is a hypoelliptic operator, then there are two positive numbers,
  $c$ and $C$, such that
  $$ \mid P^{(\alpha)}(\xi) \mid^2 \leq C (1 + \mid \xi \mid^2)^{-c} (1
  + \mid P(\xi) \mid^2 )$$
  for all $\xi \in \R^{\textit{n}}$ and for every
  multi-index $\alpha$, $\mid \alpha \mid \neq 0$,
\end{Mal}

Given a hypoelliptic operator $P(D)$, we will say that the operator $Q(D)$
is weaker than $P(D)$, $Q \prec P$, if $Q(\xi)/P(\xi)$ is bounded by a
constant in $\R^{\textit{n}}-$ infinity. We will say that $Q(D)$ is strictly weaker than
$P(D)$, $Q \prec \! \prec P$, if the quotient $\frac{Q}{P}$ tends to $0$, when $\xi \rightarrow \infty$
in $\R^{\textit{n}}$.
 According to
proposition \ref{so} a condition equivalent with hypoellipticity is
that $P^{(\alpha)}(\xi) \prec \! \prec P(\xi)$, $\alpha \neq 0$. A condition
equivalent with ellipticity for a differential operator $P$ of order $m$
is that it is stronger than any differential operator of order not
exceeding that of $P$ (\cite{Ho_LPDO}).

\vsp

For variable coefficients operators, we use the same criterion for
hypoellipticity as in Def. \ref{HE},
\newtheorem{HEvar}[Weyl]{Definition }
\begin{HEvar}[Hypoelliptic variable coefficients operator]
  Given an open set $\Omega \subset \R^{\textit{n}}$, a differential operator
  $P(x,D)$ with infinitely differentiable coefficients in $\Omega$, is
  said to be hypoelliptic in $\Omega$, if it has the
  property that
     $\mbox{ sing supp }P(x,D)u =\mbox{ sing supp }u$,
  for every $u \in \D($ $\Omega)$.
   \end{HEvar}

We will briefly indicate the techniques used in the proof of the partial
analogue. The arguments will be extensively exemplified in the
study.
\newtheorem{reginx}[Weyl]{Definition}
\begin{reginx}[ Partially regular distribution \cite{Miz_repr}]
  A distribution $T$ is said to be regular in x, in an open set
  $\Omega \subset \subset \R^{\textit{n}} \times \R^{\textit{m}}$ if for every
  product of open sets $\Omega_x \times \Omega_y \subset \Omega$, the
  distribution $< T(x,y), \alpha(y)>_y$ is in $\e(\Omega_{\textit{x}})$, for every
  $\alpha$ in $\mathcal{D} (\Omega_{\textit{y}})$.
\end{reginx}
The norm $\parallel \varphi \parallel_{s,t}$ can be defined, for
$\varphi \in \E$ and $s,t$ real numbers, as $$ \parallel \varphi \parallel_{s,t}^2= \int
\mid \widehat{\varphi}(\xi, \eta) \mid^2 (1 + \mid \xi \mid^2)^s (1 + \mid \eta
\mid^2)^t d \xi d \eta $$
A distribution $f$ is said to be of order $(s,t)$, if for every $\alpha
\in \mathcal{D}(\R^{\textit{n + m}})$, $\parallel \alpha f \parallel_{s,t} < \infty$.
\newtheorem{order}[Weyl]{Proposition }
\begin{order}[\cite{Gard},section 2] \label{ord}
  A distribution $f \in \D(\Omega)$, $\Omega$ an open set in
  $\R^{\nu}$, is regular in x if and only if, for every open set $W \subset \subset \Omega$ and for every real
  $s$, there is a number $t$, dependent on $s$ and $W$, such that
  $f$ is of order (s,t)
\end{order}
\newtheorem{PFHE}[Weyl]{Definition }
\begin{PFHE}[Partially hypoelliptic operator \cite{Miz_repr}] \label{PFHE}
  A differential operator $L(x,y,D_x,D_y)$, is said to be
  hypoelliptic in $x$ if every distribution solution $u$ to
  $Lu=f$ is regular in $x$, when $f$ is regular in $x$.
\end{PFHE}
Using the proposition \ref{ord},we see that given a distribution
$f$ of order $(s,-N)$ ($f$ is assumed partially regular) and $Lu=f$
with $u$ of order $(s',-N)$, $s'$
arbitrary, to conclude that the operator $L$
is partially hypoelliptic, we only have to prove that $u$ is
of order $(s'',-N')$ for some suitable $N'$ and $s'' > s$.

\vsp

We will limit our study to finite order distributions and for this
reason we do not initially separate between the concept of
hypoelliptic operators and homogeneously hypoelliptic operators.

\vsp
Finally, for the polynomial corresponding to an operator with constant
coefficients, $P(D_x,D_y)$, we can give the following equivalent conditions for
partial hypoellipticity;
\newtheorem{PHE}[Weyl]{Proposition }
\begin{PHE}[\cite{Ho_LPDO}, section 11.2] \label{PHE}
  The following conditions are necessary and sufficient for a
  polynomial to be hypoelliptic in x
  \begin{enumerate}
  \item[i)]
     $P^{( \alpha)}(\xi,\eta)/P(\xi,\eta) \rightarrow 0$ if $\alpha \neq 0$ and
      $\xi \rightarrow \infty$ while $\eta$ remains bounded.
    \item[ii)]
      $P$ can be written as a finite sum
          $$P(\xi,\eta)=\sum_{\alpha'=0}P_{\alpha}(\xi)\eta^{\alpha}$$
          where $P_0(\xi)$ is hypoelliptic (as a polynomial in
          $\xi$) and $P_{\alpha}(\xi)/P_0(\xi) \rightarrow 0$
          when $\xi \rightarrow \infty$ if $\alpha \neq 0$.
   \end{enumerate}
\end{PHE}

\section{ Mizohata's representation}
\label{sec:Miz_repr}
For polynomials over $\R^{\textit{n}}$, $P,Q$ with constant coefficients, assuming both are
HE (hypoelliptic), we can define the following equivalence-relation:
$P \sim Q$ if there exist $C,C'$ positive constants such that for all $ \xi \in
\R^{\textit{n}}$    $$ C \leq \frac{1 + \mid P(\xi) \mid^2}{1 + \mid Q(\xi)
  \mid^2} \leq C' $$

\vsp

Assuming that $P$ is partially hypoelliptic, we will later show, that it is sufficient to
consider the real and imaginary parts of the operator separately. We will thus, if nothing else is indicated,
assume that the operators we are considering are self-adjoint. Note that for constant coefficients operators
we always have $P \sim P^*$. An operator $L(x,y,D_x,D_y$) with coefficients in $C^{\infty}(\Omega),
\Omega$ an open set in $\R^{\textit{n}+\textit{m}}$ is called partially
formally hypoelliptic
(PFHE)
of type $M$, where $M=M(D_{x})$ is hypoelliptic, if for $(a,b)\in W$, $W$ a neighborhood of $(x^{0},y^{0})
\in \Omega$, the contact operator (the frozen operator)
\begin{enumerate}
   \item[($*$)] $\qquad \qquad L(a,b,D_x,D_y$) $\sim_x$ $M(D_x$)
\end{enumerate}
The notation $\sim_x$ will be explained below. We note that since
$M$ is assumed HE, this criterion implies that $L(a,b,D_x,D_y$)  is HE in $x$.
The class of operators equivalent with a given operator constitute a finite
dimensional vector space, (after adding 0 to the class) . This means that $M$ has the representation
$M(D_x)=\sum^r_{j=1}M_j(D_x)$, where $M_j(D_x) \sim M(D_x)$ for every $j$ and $r$ a finite integer.
Further, for a
constant coefficients operator $P(D_x)=\sum_{k=1}^r N_k(D_x)$ such
that $N_k(D_x) \sim M(D_x)$ for every $k$, we must have $M_j \sim N_k$, for every $j,k$.
\newtheorem{repr_in_equiv}{ Lemma }[subsection]
\begin{repr_in_equiv} \label{repr_in_equiv}
Assume that $M(D_{x})$ is a constant coefficients hypoelliptic operator and that
${\mathbf{P}}=( P_1, \ldots ,P_r )$
is a vector of constant
coefficients operators equivalent in strength with $M$. Assume
$P=\sum_{j=1}^r c_jP_j$, $c_j$ in $\mathbf{C}$ for all $j$,
is such that for all $j$, $P_j$ is weaker than $P$. Then $P$ and $M$ are equally strong.
\end{repr_in_equiv}
Proof:\\
Let $P^{-2}=( \sum_j \mid P_j \mid^2 )$. Clearly $P \prec
P^{-}$ and $P^{-} \prec M$. Since $$\mid
\frac{M(\xi)}{P_j(\xi)}\frac{P_j(\xi)}{P^{-}(\xi)} \mid < C$$
for large $\mid \xi \mid$, we have $P^{-} \sim M$.
Finally $$\mid
\frac{M(\xi)}{P_j(\xi)}\frac{P_j(\xi)}{P(\xi)}
\mid < C$$
for $\mid
\xi \mid$ large. $\Box$

\vsp

$\bf{Remark:}$ Assume $P=\sum_{j=0}^r P_j$, a decomposition in hypoelliptic operators, such that $P_j$ is
weaker than $P$ for every $j$. If there is an operator in this development, say $P_0$, such that $P_0 \sim
P^{-}$, then $P_0 \prec P \prec P^{-}$, that is $P \sim P^{-}$ and $P$ is hypoelliptic.

\vsp

The criterion ($*$) should now be understood using the representation
$$L(a,b,D_x,D_y)=\sum_jc_j(a,b)N_j(D_x)Q_j(D_y)$$
where $N_j,Q_j$ are constant coefficient operators, such that $N_j$
is weaker than $L'=\sum_j c_j(a,b) N_j$, for any $j$, as
\begin{enumerate}
    \item[($*'$)] $\qquad \qquad N_j(D_x)\sim$ $M(D_x)\quad \text{ for
        all } j$
\end{enumerate}
  L is said to be PFHE in $\Omega$ (particularly FHE in x), if the condition ($*$) is satisfied for every
($x_0,y_0)\in \Omega$.
Note that if we do not assume the polynomials $P,Q$ hypoelliptic, we
can use the following criterion for equivalence found in (\cite{Ho_LPDO}, sec. 10.3.4)
$$ C < \widetilde{P}(\xi)/\widetilde{Q}(\xi) < C' \qquad \xi \in \R^{\textit{n}}
$$  where $\widetilde{P}(\xi)^2=\sum_{\mid \alpha \mid \geq 0} \mid
P^{(\alpha)}(\xi) \mid^2$.

\vspace{.5cm}

\newtheorem{const_str}[repr_in_equiv]{ Lemma }
\begin{const_str}
Assume $L(\xi)=P_0(\xi') +\sum_j P_j(\xi')Q_j(\xi'')$ is the polynomial for a constant
coefficients partially hypoelliptic operator, (as in Prop. \ref{PHE}
ii)) and
$M(\xi')=\sum_jM_j(\xi')$ with a development in equivalent operators as in Lemma \ref{repr_in_equiv}, the polynomial for
a constant coefficients operator equivalent with $P_0$ in strength. Then $L \sim \sum_jM_jQ_j$.
\end{const_str}
Proof:\\ It will turn out that the real zero's of the polynomial
do not contain necessary information concerning regularity
behavior to the operator and we may assume (possibly by adding a
parameter and/or by squaring the expressions), that the polynomial
has no real zero's, for $\xi$ large.
\begin{itemize}
  \item[i)] Let $P^{+}=\sum_j (P_0+P_j)Q_j$. Then $P^{+} \sim \sum_j
    M_jQ_j$, since we trivially have
    $$
    \frac{\mid (P_0+P_j)(\xi')Q_j(\xi'') \mid}{\mid M_j(\xi')Q_j(\xi'')
    \mid}\frac{\mid M_j(\xi')Q_j(\xi'')\mid}{\mid \sum_jM_j(\xi')Q_j(\xi'')\mid} < C, \qquad
    \mid \xi \mid \text{ large }, $$
    and analogously $\sum_j M_jQ_j \prec P^{+}$.
  \item[ii)] We have $L \sim P^{+}$. It immediately follows that $L \prec P^{+}$. The
    condition $M \sim P_0$, means that for $\xi''$ fixed and bounded
    $$ \mid \frac{M(\xi')}{L(\xi',\xi'')} \mid < C \qquad \text{ for } \mid \xi'
    \mid \qquad \text{ large }$$ so $M_jQ_j \prec L$. Finally
    $$
    \frac{\mid (P_0+P_j)(\xi')Q_j(\xi'')\mid}{\mid M_j(\xi')Q_j(\xi'') \mid}\frac{\mid
    M_j(\xi')Q_j(\xi'') \mid}{\mid L(\xi) \mid} < C \qquad \mid \xi
    \mid \text{ large }$$
    We conclude $P^{+} \sim L$$\Box$
\end{itemize}

Note also that $L(x,y,D_x,D_y)$ can be written as an operator with
constant strength  $$
L(x_0,y_0,D_x,D_y) + \sum_{j=1}^r d_j(x,y)R_j(D_x,D_y)$$ where the
coefficients $d_j(x,y)$ are uniquely determined, vanish at
$(x_0,y_0)$, and have the same regularity properties as the
coefficients to $L(x,y,D_x,D_y)$. Further, the $R_j$'s are
constant coefficients operators weaker than $M(D_{x})N(D_{y})$,
where $M$ is the type operator and $N$ a constant coefficients
operator, adjusting the lower order term behavior.

\vspace{.5cm}
The representation of principal interest in this study is however the
following, due to Mizohata (section 2,\cite{Miz_repr}):
\begin{equation}
  L(x,y,D_x,D_y)=P(x,y,D_x)+\sum_{j=1}^rP_j(x,y,D_x)Q_j(D_y)
  \label{L1} \end{equation}
where
\begin{enumerate}
\item $P$ is FHE in $x$, in the sense that $P(x,y,D_x)=\sum_ja_j(x,y)M_j(D_x)$
  in an open set $W \subset \Omega$, where
$a_j\in C^{\infty}$ , $M_j \sim M$ and HE for all $ j$ ($ j=1,\ldots,r$).
\item $P_j \prec \! \prec {}_x M$, in the sense that $P_j(x,y,D_x)=\sum_kb_k(x,y)N_k(D_x)$ in $W$, where
$b_k\in C^{\infty}$, $N_k  \prec \! \prec M$ for all $k$ ($k=1,\ldots,r$).
\item $Q_j$ is a constant coefficients operator for all $ j$.
\end{enumerate}
We assume in what follows that $\text{deg}_{(x)}(L) >$ 0.
We shall in this study assume that $L$ is defined according to
Mizohata's representation on a fixed compact set $K \subset
\Omega$ and that it is defined as the type-operator $M$, outside
$K$.
\newtheorem{Prop1}[repr_in_equiv]{Proposition }
\begin{Prop1}
 If an operator  can be represented according to Mizohata
(\ref{L1}), then it is FHE in $x$ of type $M$. Conversely, any operator FHE
in $x$ of type $M$ can be written on the form above.
\end{Prop1}
Proof:  The fact that the operator with Mizohata's representation
is HE in $x$, is proven in \cite{Miz_repr}. It can also be shown,
that given operators $P$ and $Q$ with constant coefficients, such
that $Q  \prec \! \prec P$, we get $P \sim P+aQ$, for all $a \in
\mathbf{C}$. This means that $L
\sim_x M$ in $K$ for the
operator $L(a,b,D_x,D_y$) with constant coefficients. We conclude that $L$ is FHE in $x$.

\vsp

It is always possible to write the polynomial corresponding to the operator $L$ as
\begin{equation}
L(x,y,\xi,\eta)=P(x,y,\xi)+\sum_{\mid\alpha\mid >
  0}b_{\alpha}(x,y)P^{(\alpha)}(\xi)\eta^{\alpha} \label{Miz} \end{equation}
For a fixed $\eta$, we know that $L$ must be FHE in $x$ and of type
$M$, particularly if $\eta=0$,
$P(x,y,\xi$) must be FHE in $x$. According to the definition, this implies that $P(a,b,\xi$)
is HE in $x$, for a fixed $(a,b)$ in $K$, consequently $P^{(\alpha)}
\prec \! \prec M$ and finally
$b_{\alpha}(a,b)P^{(\alpha)}(\xi)\eta^{\alpha}  \prec \! \prec {}_x M(\xi)$.$\Box$

\section{ Generalized Sobolev spaces}
\label{sec:Gen_Sob}
We start with a generalization of the Sobolev spaces, written as $H(Q_1,\ldots,Q_r;\Omega)$,
where $u \in L_2(\Omega)$ and in the distributional sense, $Q_iu \in
L_2(\Omega)$,for all $i$.
Let $\mathbf{Q}$ stand for the constant coefficient operators ($Q_1,\ldots,Q_r$).
This produces separable Hilbert spaces $H(\mathbf{Q},\Omega)$, for the norm
$$\parallel u \parallel ^2_{H(\mathbf{Q},\Omega)}=\parallel u \parallel^2_{L_2(\Omega)}+
\sum_{i=1}^r \parallel Q_iu \parallel^2_{L_2(\Omega)} \qquad u \in L_2(\Omega).$$
Using Parseval's relation and applying weight functions to these
spaces, we get Hilbert spaces of compactly supported distributions, $H^{s,t}_K(\mathbf{Q})$, $H^{t}_K(\mathbf{Q})$, for real
numbers $s,t$.
The same construction holds for variable coefficients operators and we
get Hilbert spaces
$H^{s}_K(\mathbf{P})$, short for $H^{s}_K(P_1,\ldots,P_r)$, with the norm
\begin{equation} \parallel u \parallel^2_{H^{s}_K(\mathbf{P})}=\parallel u
\parallel^2_{s}+\sum_{j=1}^r\parallel P_j(x,D_x)u \parallel^2_{s}
\qquad u\in H^{s}_K \label{Lions} \end{equation}
Proof:(\cite{Lions}) \\
It is sufficient to study the case $s = t = 0$. The right side
in (\ref{Lions}) is a pre-Hilbert scalar product. If the left side
is 0, then $\parallel u \parallel^2_{L^{2}} = 0$, so $u=0$. Obviously, $L^{2}_K(\mathbf{P})$ $\subset L^{2}_K$ algebraically.
Also if $\varphi \rightarrow 0$ in $L^{2}_K(\mathbf{P})$, then $\varphi \rightarrow 0$ in $L^{2}_K$,
so the inclusion holds topologically. There remains to prove completeness
for the norm (\ref{Lions}).
Assume $\varphi_k$ a Cauchy sequence in $L^{2}_K(\mathbf{P})$, then $\varphi_k, P_i \varphi_k$ are
Cauchy sequences in $L^{2}_K$, which is a complete space, so $\varphi_k \rightarrow \varphi$
and $P_i \varphi_k \rightarrow \psi_i$ in $L^{2}_K$. Further, $P_i \varphi_k \rightarrow P_i \varphi$ in $\D$
and we conclude that $\psi_i = P_i \varphi$. $\Box$

\vsp

Note that it follows from the proof of Lemma \ref{repr_in_equiv}, that if we
assume the operators $P_j$ hypoelliptic for all $j$, then the
condition that $P=\sum_j P_j$ ( a development in weaker operators ) is
hypoelliptic, is equivalent with the condition that $P \sim P^{-}=(\sum_j \mid P_j
\mid^2)^{1/2}$. Particularly, if $P$ is hypoelliptic then the corresponding norms are equivalent.

\vsp

We now consider the operator as acting on the Hilbert spaces $H^{s,t}_K$, where
$s,t$ are real numbers, $K$ is a compact set and $$H^{s,t}_K=\{f\in{\E}(K) \quad (1+\mid \xi \mid^2)^{s/2}(1+\mid\eta\mid^2)^{t/2}\hat{f}\in L^{2}_{\xi,\eta}\}$$
$H^{s,t}_{\Omega}$ is defined as the inductive limit as $K$ varies in $\Omega$.
We also consider the Fr\'echet spaces
$H^{s,t}_{loc}(\Omega)=\{f\in \s '(\Omega);\ \ \varphi f \in
H^{s,t}_{\Omega} \ ,\text{ for all } \varphi
\in \  C^{\infty}_0(\Omega)\} $. According to (\cite{Miz_repr}), any distribution $T$ in $H^{s,t}$ can be mollified with
a tensor product of test functions,
that is for $\alpha,\beta \in \mathcal{D}$ non-negative with $\int \alpha(x)dx = \int \beta(y)dy = 1$,
$\alpha_{\epsilon}(x) = \epsilon^{-n} \alpha(x / \epsilon)$ and $\beta_{\epsilon}(y) = \epsilon^{-m} \beta(y/\epsilon)$,
we have $\alpha_{\epsilon} \otimes \beta_{\epsilon}*T \rightarrow T$
in $H^{s,t}$, as $\epsilon \rightarrow 0$.
We can also use a partial mollifier, that is
$T*'\alpha_{\epsilon}(x,y) \rightarrow T$ in $H^{s,t}$. An
argument similar to the proof above now gives a generalized Sobolev
norm, corresponding to the variable coefficients operator
P in (\ref{L1}), formally hypoelliptic in $x$.
$$ \parallel u \parallel^2_{H^{s,t}_K(P)}=\parallel u
\parallel^2_{s,t}+\sum_{j=1}^r\parallel P_j(x,y,D_x)u \parallel^2_{s,t}
\qquad u\in H^{s,t}_K $$
The norm corresponding to the Mizohata's representation (\ref{L1}) becomes
{\small $$\parallel u \parallel^2_{H^{s,-N}_K(L)}=\parallel u \parallel^2_{s,-N}  +\parallel P(x,y,D_x)u
\parallel^2_{s,-N}+
\sum_j\parallel P_j(x,y,D_x)Q_j(D_y)u \parallel^2_{s,-N}$$}
for $ u\in H^{s,-N}_K$ \par

\subsection{ The spaces $H^{s,-N}_K$}
\label{sec:Sob_sp}
 We will use the notation (\cite{Mal}) $H^{p_s,t}_K$, for
$\{f \in \E(\textsl{K}); \textit{p}_{\textit{s}}(\xi)$ $(1+\mid \eta \mid^2)^{t/2}\hat{f} \in
L^{2}_{\xi,\eta} \}$, where $p_s$ is a weight function over $\R^{\textit{n}}$.

\vsp

\newtheorem{Lem1}{Lemma  }[subsection]
\begin{Lem1}
We have the following equivalence for the weight function
$p_s(\xi)=(1+\mid\xi\mid^2)^{s/2}(1+\mid M(\xi)\mid^2)$ :
For the operator $P(x,y,D_x)$ FHE in $x$ of type $M(D_x)$, there is a
compact neighborhood $K$ of $(x_0,y_0)$, such that
$$f \in H^{s,-N}_K(P) \  \Leftrightarrow \ f\in H^{p_s,-N}_K$$
$s$ is here assumed real and $N$ integer.
\end{Lem1}
Proof:(\cite{Mal} Lemma III.1.3)\\
Let $W$ be a compact neighborhood of $(x^0,y^0)$ included in $K$. The
conditions on $P$ allow us to write
$P(x,y,D_x)=\sum_ia_i(x,y)P_i(D_x)$ with
$a_i\in C^{\infty}$ and $P_i \sim M,\quad  P_i$ with constant coefficients. For the implication
$f\in H^{p_s,-N}_K \Rightarrow
P_i(D_x)f\in H^{s,-N}_W$, we use the trivial inequality \\
(**)\hspace{3.6cm} $\parallel Mf \parallel_{s,-N} \leq C_K\parallel f \parallel_{p_s,-N}$.
\\ From the condition $P_i \sim M$ and consequently
$P_i \prec  M$, we conclude that (\cite{Miz_repr} Cor. Lemma 2.2),  $\parallel P_i(D_x)f\parallel_{s,-N}
\leq C_K \parallel M(D_x)f \parallel_{s,-N}$. Thus $P_i(D_x)f\in H^{s,-N}_W$.
Finally, $H^{s,-N}_W$ is a space of local type (\cite{Miz_repr} Cor.2  Prop. 2.1), which means that
$P_i(D_x)f\in H^{s,-N}_W \Rightarrow a_i(x,y)P_i(D_x)f\in H^{s,-N}_W$,
for all $a_i\in C^{\infty}$. This establishes the implication
$f\in H^{p_s,-N}_K \Rightarrow$ $P(x,y,D_x)f\in H^{s,-N}_K$. The opposite
implication is proven in \cite{Miz_repr} Prop. 2.4 for the local spaces, but with
the additional assumption that $f$ has compact support, the result
follows on $K$ using an appropriate test function.$\Box$

\vspace{.5cm}

We also note, that the same implications can be established for the
operator $L$, on the local spaces (\cite{Mal} Theorem III.1.4), although this
requires an adjustment of the order in the "bad" variable. Thus
$$f \in H^{s,-N}_K \quad \text{and} \quad L(x,y,D_x,D_y)f \in H^{s,-N'}_K \Leftrightarrow \ f \in H^{p_s,-N}_K$$

\subsection{ Sobolev's embedding theorem}
\label{sec:Sob_emb}
In section \ref{sec:Sob_sp}, we proved that $\parallel \cdot
\parallel_{H^{s,-N}_K(P)}$ is equivalent to
$\parallel \cdot \parallel_{H^{p_s,-N}_K}$, where in the last norm
$P$ denotes the single operator
with variable coefficients formally hypoelliptic in $x$, according to section \ref{sec:Miz_repr}. If $Mf \in H^{s,-N}_K$, then by
the regularity property of Sobolev spaces, we have $f \in H^{s+\sigma,-N}_K$, for some
positive number $\sigma$. Thus $f \in
H^{s,-N}_K(M)$. Immediately
$H^{s,-N}_K(M) = H^{p_s,-N}_K$. We conclude that we have the following
norm equivalences, $\parallel \cdot \parallel_{H^{s,-N}_K(M)}$
$\sim \parallel \cdot \parallel_{H^{p_s,-N}_K}$
$\sim \parallel M \cdot \parallel_{H^{s,-N}_K} $.
We note in particular:
\begin{equation}
\parallel Mf \parallel_{s,-N} \leq C_K \parallel f
\parallel_{H^{s,-N}_K(P)} \qquad  f \in H^{s,-N}_K \label{O1}
\end{equation}
It can easily be proved that for an operator $L$, partially formally hypoelliptic in $x$ and of type $M$,
the iterated operator $L^r$, $r$ any integer, is partially formally hypoelliptic in $x$ and of
type $M^r$.
It thus admits representation $$L^r(x,y,D_x,D_y)=P_{(r)}(x,y,D_x)+\sum_j
P_{(r),j}(x,y,D_x)Q_{(r),j}(D_y)$$ and a direct calculation shows that
$P_{(r)}=P^r$, the iterated operator FHE in $x$ being of type $M^r$.
\\[.5cm]

The inequality (\ref{O1}) above now becomes $$\parallel M^rf \parallel_{s,-N}
\leq C_{K,r}\parallel f \parallel_{H^{s,-N}_K(P^r)} \qquad f \in H^{s,-N}_K$$  Since
$(1+ \mid M(\xi) \mid^2)^r \leq C(1+ \mid M(\xi) \mid^{2r})$ and
using the well known
inequality for hypoelliptic operators (Proposition \ref{Mal} for r=1)
\begin{equation} (1+ \mid \xi \mid^2)^{kr} \leq C(1+ \mid M(\xi) \mid^2)^r \label{M1} \end{equation}
for some positive $k,C_K$ and for all $\xi \in \R^{\textit{n}}$, we get
$$\parallel f \parallel_{kr+s,-N} \leq C_K\parallel f \parallel_
{H^{s,-N}_K(P^r)} \qquad f \in H^{s,-N}_K$$
We can assume $k \leq 1$.

\vsp

{$\bf{Remark:}$ We could also show that the condition on $M$: $\mid M(\xi) \mid \geq
C\mid \xi \mid^{\varrho}$ gives $(1+ \mid \xi \mid^2)^{\varrho} \leq C(1+ \mid M(\xi)
\mid^2)$, so if $r$ is integer $(1+ \mid \xi \mid^2)^{r \varrho} \leq$ $ C(1+ \mid M(\xi)
\mid^{2r})$. In what follows we assume this condition satisfied.}

\vsp

In the following modified Sobolev's embedding theorem, we use the
notation $\U_{-\textsl{N}}$, for a tempered fundamental
solution of the operator $(1- \Delta)^{N/2}$. According to  Schwartz, $$
  \U_{\textit{l}}= \frac{\pi^{\textit{l/2}}}{\Gamma(\textit{l/2})}\mbox{Pf.}\Big[
  (\frac{\textit{r}}{\textit{2} \pi})^{\frac{\textit{l-m}}{\textit{2}}}\textrm{K}_{\frac{\textit{m-l}}
  {\textit{2}}}(\textit{r}) \Big]
  \qquad (\textrm{1}-\Delta)^{\textit{k}} \U_{\textit{l}} = \U_{\textit{l-2k}}
$$ except for $l$ even integer $\leq 0$, when
$\U_{\textit{2k}}=(\textrm{1}-\Delta)^{\textit{k}}\delta_{\textit{0}}$.
$K$ is an analytic function outside the origin, it is non
negative for $l \geq 0$ and it is exponentially decreasing towards
0 in infinity. This means that it is in $\mathcal{O}_C$ and
convolution with any distribution in $\s '$ is well defined.
(for notation, definitions and results, see \cite{Trev} chapter 30)

\vsp

\newtheorem{Prop2}[Lem1]{Proposition }
\begin{Prop2} \label{Prop2}
For $P$ an operator formally hypoelliptic in $x$ with coefficients in $C^{\infty}(\Omega)$,
for $f\in H^{kr+s,t}_K$, $K$ a compact set in $\Omega$, $s,t$ real numbers and for $kr+s
\geq \frac{n}{2}+l$, there is an integer $N$ and a constant $C_K$, such that
$$ \sup_{K, \mid \alpha \mid \leq l} \mid (D^{\alpha}_xf*''{\U}_{-\textsl{2N}})(x,y) \mid \leq
C_K\Big( \parallel P(x,y,D_x)^rf \parallel_{s,-N}+ \parallel f \parallel_{s,-N} \Big)$$
Here ${\U}_{-N}$ is defined as $(1-\Delta)^{-N/2}\delta_0$, the Laplacian is taken in $\R^{\textit{m}}$
\end{Prop2}
Proof: In general we have $v\in L^1_{\xi} \Rightarrow
{\F}^{-1}_{\eta}v\in C^0_{x}$ (here $\F_{\eta}$ is the partial
Fourier transform acting on $x \rightarrow \xi$), so if we can show that
$\F_{\eta}D_{x}^{\alpha}u(x,\eta)(\xi) \in L^1_{\xi,\eta}$, $\quad \text{for}
\quad u=(1+ \mid \eta \mid^2)^{-N}\F_{\xi}f(\xi,\eta)$, assuming
that $f\in H^{kr+s,-N}_K$, we must have that $D_{x}^{\alpha}u \in C^0_{x}$.
Thus $D_x^{\alpha}f*''{\U}_{-2N} \in C^0_{x,y}$, for
$\mid \alpha \mid \leq l$ and $\F{\U}_{-2N}=(1+ \mid \eta \mid^2)^{-N}$.
Now $\F_{\eta}D_{x}^{\alpha}u(x,\eta)(\xi)=\xi^{\alpha}\F_{\eta}u(\xi,\eta)$.
This gives using H\"olders inequality:
$$\int\!\!\!\int (1 + \mid \eta \mid^2)^{-N/2} \mid \F_{\eta}D_{x}^{\alpha}u(x,\eta)(\xi) \mid d\xi d\eta \leq$$
$$\Big( \int\!\!\!\int ( \mid \xi \mid^{\alpha})^2(1+\mid \xi
\mid^2)^{-(kr+s)} (1 + \mid \eta \mid^2 )^{-N} d\xi d\eta
\Big)^{1/2}$$
$$ \times \Big( \int\!\!\!\int (1+ \mid \xi \mid^2)^{kr+s} \mid \F_{\eta}u \mid^2 d\xi d\eta \Big)^{1/2}$$
For the first integral, $\mid \xi^{\alpha} \mid^2 \leq (1+ \mid \xi \mid^2)^
{\mid \alpha \mid} \leq (1+ \mid \xi \mid^2)^l$ and $(1+\mid \xi
\mid^2)^{kr+s-l}(1 + \mid \eta \mid^2)^{-N} \in
L^1_{\xi,\eta}$, if $2(kr+s-l) > n$ and $N > m$, that is $u\in C^l_{x}$.
For the second integral, $\F_{\eta}u=(1+\mid \eta \mid^2)^{-N/2}\F f$,
 so the last integral is  finite through the conditions on f.$\Box$
\\[.5cm]

We note that the same claim holds for the local spaces. If $x_0$
is fixed in ${\Omega}$
and $\varphi \in C^{\infty}_0(\Omega)$ with $\varphi=1$ in a
neighborhood $K$ of $x_0$ and $u\in H^{kr+s,-N}_{loc}$, we must have $u=\varphi u
\in H^{kr+s,-N}_K$. The inclusion follows from the proof above and the fact that the neighborhood is arbitrary.

\section{ The operator $\mbox{Re } L$}
\label{sec:Re_L}
We will later make use of the notion of a $\U_{\textit{t},-\textsl{N}}$-formally self-adjoint operator.
By this we mean that the operator is self-adjoint in the weighted scalar product
$$ (Lu,v)_{\textit{t,-N}} = (u,Lv)_{\textit{t,-N}} \qquad
\U_{\textit{t,-N}}(\textsl{D}_{\textit{x}},\textsl{D}_{\textit{y}})=(\textrm{1}-\Delta_{\textit{x}})^{\textit{t/2}}
(\textrm{1}- \Delta_{\textit{y}})^{-\textit{N/2}}\delta_{\textit{0}}$$
We introduce the notion of an operator partially formally self-adjoint, meaning
that the formally hypoelliptic part of the operator is $\U_{\textit{t,-N}}$-formally self-adjoint,
that is if $u,v \in H^{t,-N}_K(Q)$, $t$ a real number and $N$ integer
{\small $$ \sum_{j=1}^r (P_j(x,y,D_x)Q_j(D_y)u,Q_j(D_y)v)_{t,-N} =
\sum_{j=1}^r(Q_j(D_y)u,P_j(x,y,D_x)Q_j(D_y)v)_{t,-N}$$ }\par
assuming both sides are finite. In this study we will usually
assume the operator L formally self-adjoint as well as partially formally self-adjoint.
The proof of the following Lemma is close to \cite{Ni_72} Lemma 3.
\newtheorem{Lem2}{Lemma }[subsection]
\begin{Lem2}
If L is formally hypoelliptic in x of type M as in (\ref{L1}) and partially formally
self-adjoint, then $\mbox{Re }L $ is formally hypoelliptic in $x$ and of type $M$.
\end{Lem2}

\vsp

Assume $L=\sum_{j=0}^{r} P_j(x,y,D_{x}) Q_j(D_{y})=\sum_k b_k(x,y)N_k(D_x)Q_k(D_y)$,$Q_0=I$, our partially hypoelliptic,
variable coefficients operator, such that $Q_j,N_k$ are real, constant coefficients operators. Further,
that $L$ is self-adjoint in $H^{s,t}(\R^{\nu})$, that is
$$ ( v_{s,t}(D)L(x,y,D',D'')u,v)_{L^{2}}=(v_{s,t}(D)u,L(x,y,D',D'')v)_{L^{2}} \qquad u,v \in H^{s,t}$$
It is for our study sufficient to consider weights $v_{0,-t}(D'')=(1-\Delta_{y})^{t}$. Assuming $PQ$
is any of the $P_jQ_j$, $j \geq 1$, in the development of $L$, we have
$$
0=(v_{0,t}(D'')P(x,y,D')Q(D''))^*-\overline{(v_{0,t}(D'')P(x,y,D')Q(D''))}$$
$$ +\overline{(v_{0,t}(D'')P(x,y,D')Q(D'')}-(v_{0,t}(D'')P(x,y,D')Q(D''))=$$
$$T-2i \mbox{ Im }(v_{0,t}(D'')P(x,y,D')Q(D''))$$
Assume further $P$ self-adjoint in $H^{s,t}$, that is $L$ partially self-adjoint. Then,
$$ 0=Q(D'')\big[ P(x,y,D'),v_{0,t}(D'') \big] + \big[ Q(D''),v_{0,t}(D'') \big]P(x,y,D')+$$ $$v_{0,t}(D'') \big[
Q(D''),P(x,y,D') \big] $$
where $\big[ A,B \big]=AB-BA$. Obviously, $0=Q(D'')\big[
P(x,y,D'),v_{0,t}(D'') \big]+v_{0,t}(D'') \big[Q(D''),P(x,y,D') \big]$ and
\begin{equation} \label{SAH}
Q(D'')v_{0,-t}(D'')P(x,y,D')v_{0,t}(D'')=P(x,y,D')Q(D'')
\end{equation}
Let
$$T(x,y,D)=Q(D'')\big[
P(x,y,D'),v_{0,t}(D'') \big]+v_{0,t}(D'') \big[Q(D''),P(x,y,D') \big]+$$
$$2iv_{0,t}(D'')(\mbox{ Im }P(x,y,D'))Q(D'')$$
and $R(D'')=(1-\Delta'')^{t}Q(D'')$, then $R(D'')P(x,y,D')v_{0,t}(D'')-P(x,y,D')Q(D'')=0$ can be written
\begin{equation} \label{Leib} \sum_{\alpha'' \neq 0} D^{\alpha''}_{y}P(x,y,D')R^{(\alpha'')}(D'')v_{0,t}(D'')=0 \end{equation}
The sum on the left side in (\ref{Leib}),
can be divided into two partial sums, $\sum' + \sum''=\sum_{m < \mid \alpha'' \mid} + \sum_{m \geq
\mid \alpha'' \mid}$, where $m = \mbox{ deg }Q(\xi'')$. We now have $\sum'(x,y,D) \prec \prec
P(x,y,D')Q(D'')$ and since $\sum'+\sum''=0$, also $\sum''(x,y,D) \prec \prec P(x,y,D')Q(D'')$.
Put
$$ v_{0,-t}(D'')T(x,y,D)={\sum}^{'}(x,y,D)+{\sum}^{''}(x,y,D)+2 i (\mbox{ Im }P(x,y,D'))Q(D'').$$ Obviously, $v_{0,t}\sum',v_{0,t}\sum''
\prec \prec PQ$ and if we choose $m < 2t$, $v_{0,t}(D'')(\mbox{ Im }P(x,y,D'))Q(D'') \prec \prec
P(x,y,D')Q(D'')$, that is $T(x,y,D) \prec \prec P(x,y,D')Q(D'')$ and we conclude

\newtheorem{SAH}[Lem2]{ Proposition}
\begin{SAH}
Assume $L(x,y,D)$ a variable coefficients, partially formally hypoelliptic operator, self-adjoint and
partially self-adjoint  over $H^{s,t}(\R^{\nu})$, then $\mbox{ Re }L(x,y,D) \sim L(x,y,D)$
\end{SAH}
Note that with the additional assumption that $P(x,y,D')Q(D'')=$ $Q(D'')P(x,y,D')$, we have $\mbox{
Im }(P(x,y,D')Q(D''))=0$.

\vsp

We are also interested in the behavior of $\mbox{Re } L$ at the infinity. Consider first
the type polynomial $M$. Using the inequality for hypoelliptic operators (\ref{M1}), we see that $\mid M(\xi) \mid
\rightarrow \infty, \ \mid \xi \mid \rightarrow \infty$. If we assume $M$ with
real coefficients and also $M(\xi) \geq 1$, we have $M(\xi) \rightarrow \infty,
\mid \xi \mid \rightarrow \infty$ . We will adopt these assumptions throughout this study.

\newtheorem{Lem3}[Lem2]{Lemma }
\begin{Lem3}
With $L$ as in the previous lemma and defined on an open connected set $\Omega$, we
have that $\mbox{Re } L(x,y,\xi,\eta)$ does not change sign at the $\R^{\textit{n}}-$infinity.
\end{Lem3}
Proof: First assume that $L$ has constant coefficients as an operator hypoelliptic in $x$.
According to \cite{Ho_LPDO} (sec.11.1), if $d(\xi,\eta)$ denotes the distance from
$(\xi,\eta)\in \R^{\textit{n}+\textit{m}}$
to the surface $ \{ (\zeta,\upsilon)\in \mathbf{C}^{\textit{n}+\textit{m}};$ $L(\zeta,\upsilon)=0 \} $,
it follows that $d(\xi,\eta) \rightarrow \infty \ \text{as} \  \xi \rightarrow \infty$,
while $\eta$ remains bounded. Thus $\mbox{Re } L$ cannot change
sign at the
$\R^{\textit{n}}$-infinity.
\\[.5cm]

For the variable-coefficients case, we make use of the assumption that $\Omega$ is
connected. We know that for every fixed $(x,y) \in \Omega$, $\mid
\mbox{Re } L(x,y,\xi,\eta)
\mid \rightarrow \infty$, as  $\xi \rightarrow \infty$, while $\eta$ bounded. Let
$E\pm =$ $\{ (x,y)\in \Omega; \mbox{Re } L(x,y,\xi,\eta) \rightarrow
\pm \infty $, $\quad \xi \rightarrow \infty $,  $ \eta \  \text{bounded} \}$. We have
already proved that $E_+ \bigcap E_- = \emptyset$. Since $\mbox{Re } L$ is continuously dependent
on $(x,y)$, the sets $E \pm$ must be open, but this would give a separation of
$\Omega$, which contradicts the assumption that $\Omega$ is connected! We can
without loss of generality assume that $\mbox{Re } L(x,y,\xi,\eta) \rightarrow \infty$,
at the $\R^{\textit{n}}$-infinity. $\Box$

\newtheorem{imroot}[Lem2]{ Lemma }
\begin{imroot}
If $P$ is a reduced (cf. section \ref{sec:Lin}) and real operator on $L^{2}$, then $P^{1/N}$ is defined as an operator on $L^{2}$
for some positive integer $N$.
\end{imroot}

Proof:
We can assume $P(\xi) \rightarrow \infty$ as $\mid \xi \mid \rightarrow \infty$. It is sufficient
to study the polynomial for $\mid \xi \mid$ large, that is if $\tilde{P}$ is the polynomial adjusted to
a constant on $\mid \xi \mid \leq R$, we have $D(P_0)=D(\tilde{P_0})$. The
minimal operator, $\tilde{P}_0$, can thus be considered as a positive, selfadjoint operator
on $L^{2}$ and the "Gelfand-Naimark-theorem" gives unique existence of a selfadjoint operator
$\tilde{P}_0^{1/N}$, for some $N$.$\Box$

\vsp

Assume $P=P_1+iP_2$ homogeneously (geometrically) hypoelliptic on $\D$. This means particularly that
$\frac{P_1(\xi + \vartheta)}{P_1(\xi)} \rightarrow 1$ and $\frac{P_2(\xi + \vartheta)}{P_2(\xi)}
\rightarrow 1$ as $\mid \xi \mid \rightarrow \infty$ and $\mid \vartheta \mid \leq A$ for some
constant $A$ and $\xi,\vartheta$ real. If $P_1/P \rightarrow \alpha=1/(1+iC)$, with $C$ a real
constant, then $P_2/P_1 \rightarrow -C$ as $\mid \xi \mid \rightarrow \infty$. If $\mid \xi
\mid^{\sigma} \leq P_1(\xi)$, $P_2(\xi) \leq \mid \xi \mid^{\delta}$, then $\mid P_2^{1/N}/P_1 \mid
\leq \mid \xi \mid^{\delta/N-\sigma}$ as $\mid \xi \mid \rightarrow \infty$. If $\sigma >
\delta/N$, we have that the corresponding constant $C$, can be selected as $0$ and $P_1 + iP_2^{1/N}$ is hypoelliptic.
Conversely, any polynomial operator, hypoelliptic in $\D$, can be constructed in this way.

\vsp

Assume $P$ with $\sigma > 0$ and real. If $E$ is a fundamental solution to $P$, we have
$\parallel Ef \parallel_{L^{2}} \leq C \parallel f \parallel_{L^{2}}$. Assume $f$ with first
derivatives in $L^{2}$, since $P$ is real we have that $\overline{\delta}E=E\overline{\delta}$,
why $E: H(\Omega) \rightarrow H(\Omega)$ for an open set $\Omega \subset \R^{\textit{n}}$. As
$P(D)u=f$ can be written $Ef=u$ in $L^{2}$, we conclude that $P$ is analytically hypoelliptic in
$L^{2}$.
Assume $P(D_{x})$ an operator with constant coefficients and $I_E$ a parametrix. Then, $$\overline{\delta}PI_E-PI_E\overline{\delta}=
\big( \overline{\delta}PI_E-P\overline{\delta}I_E \big)+\big(P\overline{\delta}I_E-PI_E\overline{\delta} \big)$$
The first expression is trivially $0$. If $u \in H \cap L^{2}(\Omega)$, $\overline{\delta}u=0$, so $- \overline{\delta}PI_E(u)=P\big(\overline{\delta}I_E-I_E\overline{\delta} \big)$. Since $PI_E=I-\gamma$, $\gamma$ regularizing, we have
$$\overline{\delta}\gamma u=P \big(I_E\overline{\delta}-\overline{\delta}I_E \big)u \in C^{\infty}.$$ Thus, if
$P$ is homogeneously hypoelliptic in $L^{2}$, $\Psi=I_E\overline{\delta} - \overline{\delta}I_E$ is regularizing. Finally,
$\parallel \overline{\delta}I_E u \parallel \leq C \parallel \overline{\delta} u \parallel$ for a constant $C$,
given that the set $\Omega$ is bounded. Thus $I_E: H(\Omega) \rightarrow H(\Omega)$ for bounded $\Omega$.

\newtheorem{A-HE}[Lem2]{ Proposition }
\begin{A-HE} \label{A-HE}
  If $P$ is a constant coefficients operator, hypoelliptic in $L^{2}$, then $P$ is analytically (homogeneously)
  hypoelliptic in $L^{2}$. Conversely, if $P$ is analytically (homogeneously) hypoelliptic in $L^{2}$, $P$ is
  hypoelliptic in $L^{2}$.
\end{A-HE}

\section{ Some remarks on fundamental solutions}
\label{sec:Fund_sol}
Assuming that $\Omega$ is an open set in $\R^{\textit{n}+\textit{m}}$, we let $\Omega_{\varphi}=$
$ \{ (x,y) \ ; \{ (x,y) \} - \mbox{supp }(\delta_0 \otimes \varphi)$
$ \subset \Omega \} $, where $\delta_0 \otimes \delta_0= \delta_0$ is the
Dirac measure
in $ \R^{\textit{n}+\textit{m}}$ and $\varphi \in C^{\infty}_0 (\R^{\textit{m}})$. We
then know
that the partial convolution $u *'' \varphi = u*(\delta_0 \otimes \varphi)$
is well defined in $\Omega_{\varphi}$, for $u \in {\D}(\Omega)$. In the same way
we can define $u*'\varphi$ in $$ \{ (x,y); \{ (x,y) \}-  \mbox{supp }
\varphi \times {0} \subset \Omega \},$$ for $\varphi \in C^{\infty}_0(\R^{\textit{n}})$.
The partial convolution $u*''v$, can also be defined for $u \in {\D}$ and
$v \in {\E}$ in $\Omega_v$. We note that:
\begin{equation}
u \in H^{s,-N}_{loc}(\Omega) \Rightarrow u*''\varphi \in
H^{s,0}_{loc}(\Omega_{\varphi}) \label{P1} \end{equation}
(cf. \cite{Boi}). For any partial differential operator $Q$ with constant coefficients, we
have a fundamental solution
$E_0$ in $ {\D}^F(\R^{\textit{n}}) $, that is distributions of finite order (and a $E_0 \in
B_{\infty,\tilde{Q}}^{loc}$, \cite{Ho_LPDO}, sec. 10.2 ). Assume for instance $ Q(D_y)E_0 = \delta_0 $. We then
have for $u \in {\E}(\Omega)$,  $$ Q(u*'' E_0) =u*( \delta_0 \otimes QE_0 )=u.$$ For a hypoelliptic operator $M(D_x)$
with constant coefficients, we have a fundamental
solution as above, now also with $\mbox{sing supp } F_0= \{ 0 \}$ and as
before $M(u *' F_0)=u \  \text{in}
\  \Omega$. By multiplying $E_0$ with a $\chi \in C^{\infty}_0(\R^{\textit{m}})$,
$ \chi = 1$ in a neighborhood of a compact neighborhood of $(0,0)$, $K_y$, we have
$\chi E_0 \in
{\E}$. Define a linear operator $E: {\E} \rightarrow {\E}$ by Ef = $f*''\chi E_0$
for $f \in {\E}(\Omega)$. Define $F: {\E} \rightarrow {\E}$ as $Fu=u*'\vartheta F_0$, $F_0$ corresponding to $M$ and
$\vartheta=1$ in a neighborhood of $K_x$ as above. Obviously, we have
$$E: H^{s,t'}_W \rightarrow H^{s,t}_W(Q), \quad
F: H^{s,t'}_W \rightarrow H^{p_s,t'}_W,$$ for some $t' > t$
and $EFf=$ $f*(F_0 \otimes E_0)$ on $W$. Here $W$ is a compact set in $\Omega$.

\vsp

For the operator $P(x^0,y^0,D_x)$, we can now construct a fundamental solution
with singularity in $(x^0,y^0)$ as follows. Let $G_0$ be the fundamental solution, modified to ${\E}$,
corresponding to $$ P(x^0,y^0,D_x)=\sum_jP_j(x^0,y^0,D_x)$$ on a compact neighborhood of
($x_0,y_0$), $W$ and $Gf=f*'G_0$. If
$f=\delta_{(x_0,y_0)}$ and $I$ is the trivial mapping, we get the fundamental solution to $P$ as
$GI\delta_{(x_0,y_0)}=G_0 \otimes \delta_{y_0}$. It can also be shown that a linear mapping $F$
can be constructed for a formally hypoelliptic operator, which gives a fundamental solution to
the variable coefficients operator $P(x,y,D_x)$, but we will use another
method to produce this solution in chapter \ref{sec:Levi}. We have the following lemma

\newtheorem{Fund}{Lemma  }[subsection]
\begin{Fund}
  There is a fundamental solution with singularity in $(x_0,y_0)$, to a
  constant coefficients differential operator  on $\R^{\textit{n}+\textit{m}}$,
  independent of the variables in $\R^{\textit{m}}$-space,
  $P(x_0,y_0,D_x)$, on the form
  $G_0 \otimes \delta_{y_0}$, where $G_0$ is a fundamental solution
  with singularity in $x_0$, to the operator,
  considered as an operator on $\R^{\textit{n}}$.
\end{Fund}

\section{ Realizations}
\label{sec:Realiz}
We assume $\Omega$ open and connected and of dimension $n>1$. We
can define a linear mapping $\tilde{M_0}$ on $C^{\infty}_0(\Omega)$,
as $\tilde{M_0}\varphi=
M(D_x) \varphi$. This is a symmetric operator, densely defined in $L^{2}(\Omega)$.
The minimal operator, $M_0$, or the strong extension, is defined as
the closure of
$\tilde{M_0}$ in $L^{2}(\Omega)$. This is a closed, linear, symmetric, densely
defined mapping on $L^{2}(\Omega)$. The adjoint, $M^{*}_{0}$, is a closed, linear
densely defined operator and $u \in \mathcal{D}_{M^{*}_{0}}$ if and only if in
distribution sense
$M^{*}_{0}u=M(D_x)u \in L^{2}(\Omega)$. $M^{*}_{0}$ is called the maximal operator or
the weak extension. $M_0$ is self-adjoint if and only if $M^{*}_{0}$ is also symmetric.
\\[.5cm]

For a self-adjoint extension, $\A$, of $M_0$ we must have
$M_0 \subseteq \A \subseteq$ $M_0^{*}$ and $\A$ is called a self-adjoint realization of $M(D_x)$ in $L^{2}(\Omega)$.
The minimal- and maximal operators can also be defined for the operator
$P(x,y,D_x)$ just as above. We will later on rationalize the
calculations by comparing norms for realizations corresponding to
different operators. For this we need:
\newtheorem{Lem4}{Lemma }[subsection]
\begin{Lem4}
For the operator $P$, FHE in $x$ of type $M$, we have that the
domains of the respective closures in $H^{0,-N}_K$ of the operators coincide. That is
$$\mathcal{D}(M \tilde{\otimes} I)=\mathcal{D}(P \tilde{\otimes} I)$$ In particular (\cite{Schech},Ch. 5) $\mathcal{D}({(M \otimes I)}_0)$ $=\mathcal{D}({(P \otimes I)}_0)$.
For Mizohata's representation $L$, we have that the domains of $\tilde{L}$ and
$P \tilde{\otimes} I$ are partially coinciding and $\mathcal{D}(\widetilde{L}) \subset $ $\mathcal{D}( \textsl{P}
\tilde{\otimes} \textsl{I})$.
\end{Lem4}

Proof:  Assume for $\varphi_n \in C^{\infty}_c(K)$,where $C^{\infty}_c$ denotes $\{ f \in C^{\infty} ; \mbox{ supp }f \subset K \}$
, a dense subset of $H^{s,t}_K$, $K$ is a
compact set in $\Omega$, where $P$ is defined,
$\varphi_n \rightarrow 0$ and $P(x,y,D_x) \varphi_n \rightarrow f$ in $H^{0,-N}_{K}$.
If we can show that
\begin{equation}
 \parallel M(D_x) \varphi \parallel_{0,-N} \leq C_1
\parallel \varphi \parallel_{H^{0,-N}_{K}(P)} \label{M_P} \end{equation}
we see that $M(D_x)\varphi_n$ is convergent in $H^{0,-N}_{K}$ and since
$M \tilde{\otimes} I$ is a closed operator, $M(D_x) \varphi \rightarrow 0$.
If in addition
\begin{equation} \parallel P(x,y,D_x) \varphi \parallel_{0,-N} \leq C_2
\parallel \varphi \parallel_{H^{0,-N}_K(M)} \label{P_M} \end{equation}
and since we know that the right hand side tends to $0$, we must have that
$ P(x,y,D_x) \varphi_n \rightarrow 0$ in $H^{0,-N}_K$, that is $f=0$. The
inequalities (\ref{M_P}),(\ref{P_M}) follow immediately from what was said in
section \ref{sec:Sob_sp}. The second statement also follows from section \ref{sec:Sob_sp}, assuming
the rightmost norms finite
$$\parallel P(x,y,D_x) \varphi \parallel_{0,-N} \leq C'_0 \parallel
\varphi \parallel_{H^{0,-N}_K(M)} \leq C'_1 \parallel \varphi
\parallel_{H^{0,-N'}_K(L)}$$
 and
$$\parallel L(x,y,D_x,D_y) \varphi \parallel_{0,-N'} \leq C'_2
\parallel \varphi \parallel_{H^{0,-N}_K(M)} \leq C'_3 \parallel\varphi
\parallel_{H^{0,-N}_K(P)}.\Box$$

\section{ Some remarks on Schwartz kernels }
\label{sec:S-kernels}
Given $f \in C^{\infty}$ $(X \times Y)$, $X,Y$ open sets in
$\R^{\nu}$, we can define an integral operator on $\mathcal{D}$ $(Y)$
according to
$$ I_f (\psi) (x)= \int f(x,y) \psi(y) d y $$
which is continuous $\mathcal{D}$ $(Y) \rightarrow C^{\infty}$ $(X)$. This is a
regularizing operator, it can be extended to a continuous linear
operator $\E$ $(Y) \rightarrow C^{\infty}$ $(X)$. More generally, we
have Schwartz kernel theorem, that is given $K \in \D$ $(X \times Y)$, we
can define a continuous, linear operator $\mathcal{D}$ $(Y) \rightarrow
\D$ $(X)$ according to
\begin{equation}
  < I_K (\psi), \varphi > = < K, \varphi \otimes \psi > \qquad \psi \in
  \mathcal{D}(\textsl{Y}), \varphi \in \mathcal{D}(\textsl{X}) \label{Schwartz}
\end{equation}
and conversely, given the continuous, linear operator $I_K$, we have
unique existence of a kernel $K \in \D$ $(X \times Y)$ such that
(\ref{Schwartz}) holds.

For a constant coefficients differential operator $Q$, let
$Q'=Q(D_{x}),Q''=Q(D_{y}),\overline{Q''}=Q(-D_{y})$. For test functions $f,g \in
C^{\infty}_0(\R^{\nu} \times \R^{\nu})$,
  $$I_{g}(I_{Q'f}(\phi))(x)=I_{g}(Q'I_{f}(\phi))(x)=I_{\overline{Q''}g}(I_{f}(\phi))(x)$$
Assume $f \in C^{\infty}_0(\R^{\nu} \times
\R^{\nu})$ such that for constant coefficients operators,
$\overline{Q'}f=Q''f$. Assume $P,Q,R$ such operators, then
\begin{itemize}
  \item[v)] $I_{\big[ \big[ P''f,Q''f \big], R''f \big]}=I_{\big[
      \big[ f,P''f \big],Q''R''f \big]}=I_{\big[ \big[ f,f
      \big],P''Q''R''f \big]}$
\end{itemize}

\vsp

Finally, we will have use for the concept of analytic functionals as in \cite{Mart}.
Assume $V$ a complex analytic variety. For analytic functionals
defined on $V$, we have
$$I_F: H(V_y) \rightarrow H'(V_x) \qquad F \in H'(V_x \times V_y) $$
If $V$ is countable in the infinity, it can be shown that $H'(V)$ is a nuclear $({\mathcal{DS}})$-space (the dual
topology to the $(\mathcal{FS})$-topology, that is Frechet with Schwartz property). Iteration of the integral operators $I_F$
is, for our applications, defined by convolution. We note that (\cite{Mart} Ch. 1, Proposition 1.1)
\emph{Any analytic functional defined on $V$, can be represented by a measure with compact support}

In the same way as in $\D$, we can using nuclearity in $H'$,
determine the kernel to the evaluation functional, as the parameter
$x$ varies. For $\mbox{ Id }_x \in H'$, we have $< \mbox{ Id }_x,\psi >=
\int\phi(z) \psi(z) dz =< \delta_x(\phi),\psi >$, with a $\phi \in
C^{\infty}_0(\R^{\nu})$ and $\psi \in H(\mathbf{C}^{\nu})$. Let
$\Delta \in H'(\mathbf{C}^{\nu} \times \mathbf{C}^{\nu})$, be the kernel
representing $\mbox{ Id }_x$. We note, that in
a neighborhood of $x$, $V$, we have for this kernel,
$<I^N_{\Delta}(\varphi),\psi >_V=<I_{\big[ \Delta,\Delta
\big]_N} (\varphi),\psi >_V=0$, with $\varphi,\psi \in
H(\mathbf{C}^{\nu})$, for some $N$ (cf. \cite{jag_II}).

\section{ The spectral kernel }
 According to the spectral theorem, any self-adjoint realization as in
 the previous section, can be
 represented as a spectral operator of scalar type:
$$ \A = \int \lambda dE_{\lambda} \qquad \lambda  \text{\ real} $$
For the iterated operators, we also have realizations that are of
scalar type:
$$ \A_{L^r} =  \int \lambda dE_r(\lambda)
\qquad  \text{where} \quad E_r(\lambda^r)= E(\lambda) \ \text{ , r is
  assumed odd}  $$
Assuming the constant coefficients operator $L$ defined on $\R^{\nu}$, since $L$ is assumed formally self-adjoint,
 the realization
is related to the symbol of the operator as $\A_L = {\mathcal F}L(\xi,\eta)\F^{-1}$.
Let $\{ E_{\lambda} \}$ be the spectral family corresponding to the realization
$\A_L$ with spectrum $\sigma ( \A_L )= \{ L(\xi,\eta) \quad \xi,\eta \quad  \text{ real} \}$. We can
assume the spectrum left semi-bounded. We then know that $E_{\lambda}$ is related to
$\omega=\chi_{(-\infty,\lambda)} \circ L$, through $ \mathcal{F} ( E_{\lambda}u ) = \omega \mathcal{F} u$
with $\mathcal{F}^{-1} \omega \in \D$ (\cite{Ni}).

\vsp

For the type polynomial $M(\xi)$, we know that $M(\xi) \rightarrow \infty$ as
$\mid \xi \mid \rightarrow \infty$, this means that the set $G_{\lambda}=\{ \xi; M(\xi) < \lambda \}$
is bounded for every $\lambda$. Then, for the spectral kernel $ e_{\lambda}(x_0,x) =$
$\F^{-1}\chi_{G_{\lambda}}(x_0-x) \in$ $C^{\infty}( \R^{\textit{n}} \times  \R^{\textit{n}})$ (\cite{Ni_92}).
 We have
\begin{equation} e_{\lambda}(x_0,x) = (2 \pi)^{-n} \int_{M(\xi) < \lambda}
e^{i(x_0-x) \cdot \xi} d \xi \label{e_reg} \end{equation}

\vspace{.5cm}

Let $(x,y)=(x',x'') \in \R^{\nu}$,  $\nu=\text{n}+\text{m}$.
For the contact polynomial $L(\xi)$, that is the polynomial
corresponding to the frozen operator $L(x_0,D_{x'},D_{x''})$, this argument only holds for the hypoelliptic part
of the polynomial. Consider instead
$$L(x_0,D_x)=\sum_jc_j(x_0)P_j(D_{x'})Q_j(D_{x''}),$$ where $P_j \sim M$ and
$P_j,Q_j$ are constant coefficients, finite order operators. We can
write $R(P,Q)=L$, where $R$ is a constant coefficients polynomial.
Let $T_1=P \otimes Id_1, T_2=Id_2 \otimes Q$, where $Id_j,j=1,2$ are
identity operators in respective complex Banach space. It can be verified for the joint spectrum,
$\sigma(T_1,T_2)=\sigma(T_1)\sigma(T_2)=\sigma(P)\sigma(Q)$. Further
$T_1T_2=T_2T_1$. Finally, $\sigma[R(T_1,T_2)]=R[\sigma(T_1,T_2)]$,
(\cite{SchechII}).

\vsp

Assume $E_{\lambda,\mu}$ the spectral resolution corresponding to the
operator $L$ and consider the restriction to $\mathcal{D}(\R^{\nu}
\times \R^{\nu})$. Let $E_{\lambda},E_{\mu}$ be the resolutions
corresponding to $T_1$ and $T_2$, restricted to $\mathcal{D}$. According to what has been said
above, for
$v \in \mathcal{D}$, \\ $\F$ $(E_{\lambda,\mu}v)=$ $w\F$ $v$ with
$w=w_{\lambda}w_{\mu}$ and $E_{\lambda}v=$ $v*'\F^{-\text{1}}$ $w_{\lambda}$
and $E_{\mu}v=v*''\F^{-\text{1}}$ $w_{\mu}$.
The tensor product of projections
$E_L(J \times K)=E_{T_1}(J)E_{T_2}(K)$, is well defined for
Baire sets $J,K$ in $\R$ ($Q$ is assumed to have real
coefficients). Assume $E_{\lambda}$ has spectral kernel
{\small
$$t^1_{\lambda}(x,y)=e_{x',\lambda} \otimes \delta_{x''}= \frac{1}{(2
  \pi)^{\nu}} \int_{P(\xi') < \lambda} e^{i (x-y) \cdot \xi} d \xi$$
and $E_{\mu}$ with spectral kernel $t^2_{\mu}=\delta_{x'} \otimes e_{x'',\mu}$.
Note that $e_{x',\lambda}(x',y')=$ \\ $\F^{-\text{1}}$ $w_{\lambda}(x'-y')$ and
analogously for $e_{x'',\mu}$.
We then have, for $f \in \mathcal{D}(\R^{\nu}_{\text{y}})$
$$E_{\lambda,\mu}(f)(x)=I_{e_{x,\lambda,\mu}}(f)(x)=I_{t^1_{\lambda}}(I_{t^2_{\mu}}(f))(x)=I_{\big[
  t^1_{\lambda},t^2_{\mu} \big]}(f)(x)=I_{e_{x',\lambda} \otimes e_{x'',\mu}}(f)(x)$$}
We then know for the spectral function corresponding to the operator
$P$ (defined as the type operator outside a compact set), $e_{x',\lambda} \in C^{\infty}(\R^{\textit{n}} \times \R^{\textit{n}})$ and
for the spectral kernel corresponding to the operator $Q$ (defined as
the identity operator outside the same compact set),
$e_{x'',\mu} \in \D(\R^{\textit{m}} \times
\R^{\textit{m}})$. Using that
$E_{\lambda,\mu}$ is an orthogonal projection,
$$ E^2_{\lambda,\mu}=I_{\big[ e_{x,\lambda,\mu},e_{x,\lambda,\mu}
  \big]}=I^{'2}_{e_{x',\lambda}}I^{''2}_{e_{x'',\mu}}=I'_{e_{x',\lambda}}I''_{e_{x'',\mu}}=E_{\lambda,\mu}$$
Particularly, if $\widetilde{e}_{x'',\mu}= e_{x'',\mu}*_x\varphi*_y\psi$, for test
functions with support in a neighborhood of the origin in $\R^{\textit{m}}$,
$$I_{\big[ e_{x'',\mu}*_x\varphi, e_{x'',\mu}*_y\psi
  \big]}=I_{\widetilde{e}_{x'',\mu}}$$
and since $\widetilde{e}_{x'',\mu} \in C^{\infty}(\R^{\textit{m}} \times
\R^{\textit{m}})$, we can define
$\widetilde{E}_{\lambda,\mu}=I_{e_{x',\lambda} \otimes
  \widetilde{e}_{x'',\mu}}(f)(x)$ for $f \in L^{2}(\R^{\nu})$.
This means that we have partial regularity, for the spectral kernel.

\section{ Topology. Notation and fundamentals}
  \label{sec:Levi_lemmas}
  In what follows, the presentation is more or less a modification of the article
  (Nilsson \cite{Ni_72}) and we refer to this work for proofs and arguments.
  Assume the operator $L(x,D_x)$ defined as partially formally hypoelliptic in a neighborhood of $x_0$, included in $K$,
  $K$ a compact subset in $\Omega \subset \R^{\nu}$
  and as the type operator with constant coefficients outside $K$.
  We also assume, as has been mentioned earlier, that the type polynomial $M$ is real
  and such that $M(\xi') \geq 1$, for all $ \xi' \in \R^{\textit{n}}$ and for positive constants
  $C$ and $\varrho$,  $M(\xi') \geq C \mid \xi' \mid^{\varrho}$, for all $ \xi' \in \R^{\textit{n}}$.
  According to these conditions, we get for this modified operator $
  \mbox{Re }L(x,\xi) \rightarrow \infty$
  as $\xi' \rightarrow \infty$ while $\xi''$ bounded.

\vsp

  Next we will construct and estimate a fundamental solution
  to the variable coefficients operator $L(z,D_z)-\lambda$ and for this purpose we need
  an estimate of the distribution $$\alpha_{\lambda}(x,z) = [P^x(D_{z'})-L(z,D_z)]K^+_{\lambda}(x,z)$$ Here $K^+_{\lambda}(x,z)$
  is a fundamental solution with singularity in $x$, corresponding to the modified operator
  $P(x,D_{z'})-\lambda$,$\  \lambda$ is assumed large and negative. So $K^+_{\lambda}(x,\cdot) \in
  H^{t,-N}_{loc}(\R^{\nu})$, can be written as
  $\chi_{x_0}(x)h_{\lambda} \otimes \delta_{x''}(x,z) +$ $(1 - \chi_{x_0}(x))t_{\lambda}
  \otimes \delta_{\textsl{0}}(x,z)$,
  where $\chi_{x_0} \in C^{\infty}(K)$ assumes the value 1 in a
  neighborhood of $x_0$,  $t_{\lambda}$ is the fundamental solution to
  the type operator
  and $h_{\lambda}$ the fundamental solution to
  the operator $P(x,D_{z'}) - \lambda$ (section
  \ref{sec:Fund_sol}). With this
  modification $\alpha_{\lambda}(x,\cdot) \in \E(\R^{\nu})$.
  \\[.5cm]

  Consider the norms $$ M_{\alpha}(u) = \int \mid \xi \mid^{\alpha} \mid
  \mathcal{F} u(\xi) \mid d \xi $$ over $\mathcal{S}'(\R^{\nu})$, for any multi-index $\alpha$,
  where $\mid \xi \mid^{\alpha} = \sum_{\beta \leq \alpha} \mid \xi^{\beta} \mid$.
  These norms $M_{\alpha}$ define Banach spaces $\B_{\alpha}$ with the test functions
  $C^{\infty}_0$ as a dense subset. Also consider the operator norms
  $$ N^{\alpha,\beta}(L) = \sup_{0 \neq u \in \B_{\alpha}}
    \frac{M_{\beta}(Lu)}
    {M_{\alpha}(u)} $$ over linear mappings L from $\B_{\alpha}$
    to $\B_{\beta}$
    \\[.5cm]

\section{ Complex translations in the polynomials }
\label{sec:transl}
      Now to the distribution $\alpha_{\lambda}(x,z)$
      $=\beta_{\lambda}(x,z) +r_{\lambda}(x,z)$, where
      \begin{equation}
      \beta_{\lambda}(x,z)=(P^x(D_{z'}) -
      P(z,D_{z'}))K^+_{\lambda}(x,z) \label{beta-lambda}
      \end{equation}
      Since
      $(1-\Delta'')^{-t}\delta_{x''}$ has Fourier transform in
      $L^1(\R^{\textit{m}})$, for $t > m$, we see that
      $\beta_{\lambda}*_x''\U_{-2t}(x,z)$
      $\in (\B_{\textit{0}})_{\textit{x}}$, where $\U_{-t}$ is according to
      proposition \ref{Prop2}. For the last term, we know that
      $r_{\lambda}(x,z)=\sum_jP_j(z,D_{z'})Q_j(D_{z''})K^+_{\lambda}(x,z)$
       weighted with $(1-\Delta'')^{-t'}$, $t' > max_j \{ deg(Q_j) \}+m$, is in $(\B_{\textit{0}})_{\textit{x}}$.
        Further for $t' > max_j \{deg(Q_j) \}$, $\F_z \big(r_{\lambda}*_x''\U_{\textit{-2t}'} \big)$
         is essentially bounded. The following proposition and proof is close to \cite{Ni_72} Lemma 8.
          Let $\alpha_{\lambda,N,t}(x,y)=\alpha_{\lambda}*''_x \U_{\textit{-2N}} *''_{\textit{y}}
           \U_{\textit{-2t}}(\textit{x,y})$ and analogously for $\beta_{\lambda,N,t}$.
      \newtheorem{Prop4}{Proposition }[subsection]
      \begin{Prop4} \label{Prop4}
         There are positive numbers c and $\kappa_0$ such that
           \begin{equation} N^{{\alpha},{\alpha}}(\text{exp}({\kappa} \mid {\lambda}
           \mid^b(z_j'-x_j'))I_{{\alpha}_{\lambda,N,t}}(x,z) )=O(1) \mid {\lambda} \mid^{-c} \qquad
           \lambda \rightarrow - \infty  \end{equation} for every
           j, $1 \leq j \leq n$, for every real
           ${\kappa}$ such that $\mid {\kappa} \mid \leq {\kappa}_0$,
           for every multi-index ${\alpha}$ and for some positive numbers
           N,t. The norms $N^{\alpha,\alpha}$ are
           taken with respect to z.
    \end{Prop4}

    {$\bf{Remark:}$ Here $b$ is associated to the type polynomial according to
    \begin{equation} \mid D^{\alpha}M(\xi') \mid \leq C(1 + \mid M(\xi') \mid)
    ^{1-b \mid {\alpha} \mid} \label{b} \end{equation}}
\newtheorem{Prop1,Lem1}[Prop4]{Lemma  }
\begin{Prop1,Lem1} \label{L2}
    For $\beta_{\lambda}(x,z)$ in (\ref{beta-lambda})
    $$ N^{\alpha,\alpha}(I_{\beta_{\lambda,N,t}}(x,z))=O(1) \mid \lambda \mid^{-c} \quad
    \lambda \rightarrow - \infty$$ with $N,t$ as in the
    proposition.
\end{Prop1,Lem1}

\vsp

Proof: We write $\beta_{\lambda}$ on the form $$
\beta_{\lambda}(x,z)=\sum_j \Big[b_j(x) - b_j(z)
\Big]N_j(D_{z'})K^+_{\lambda}(x,z)$$ where $b_j$ are the
$C^{\infty}$-coefficients to the operator $P(z,D_{z'})$ (the
coefficients to $\beta_{\lambda}(x,z)$ vanish in $x=z$), $N_j$ a
constant coefficient operator equivalent in $x'$ with the type
operator and $K^+_{\lambda}(x,z)=h_{\lambda} \otimes
\delta_{x''}$, the fundamental solution to the operator
$P(x,D_{z'})-\lambda$. Expansion in a Taylor series gives
$\beta_{\lambda}(x,z)=\sum_{\lambda,\mu}F_{\lambda,\mu}(x,z)$,
where each term is on a form $$ F_{\lambda}(x,z)=\sum_{0<\mid
\beta \mid < k} i^{\mid \beta
  \mid}(\beta !)^{-1} \Big(D^{\beta}b(z)
\Big)(x-z)^{\beta}N(D_{z'})K^+_{\lambda}(x,z) + R_{\lambda}(x,z)$$
and
{\small $$ R_{\lambda}(x,z)=\frac{1}{(k-1)!} \Big( \int^1_0 (1-k)^{k-1}
\frac{d^k}{dt^k} b(z + t(x-z)) dt \Big) (x-z)^{\beta}
N(D_{z'})K^+_{\lambda}(x,z)$$} \par
We now study the two mappings
$$ L_1: u \rightarrow (D^{\beta}b)u \qquad
 L_2: v \rightarrow const \int
(x-z)^{\beta}N(D_{z'})K^+_{\lambda}(x,z)v(z)dz$$
and we aim to prove a weighted analogue to
\begin{equation}
    N^{\alpha,\alpha}(L_2 \circ L_1)=O(1)\mid \lambda \mid^{-c} \quad
    \lambda \rightarrow - \infty \label{est1}
\end{equation}
\begin{equation}
  N^{\alpha,\alpha}(I_{R_{\lambda}})=O(1)\mid \lambda \mid^{-c}
    \quad \lambda \rightarrow - \infty \label{est2}
\end{equation}
Proof of (\ref{est1}):\\
The first mapping is immediate, since we already know the coefficients
are in $\B_{\alpha}$ (or constant), which is a Banach algebra. So
$N^{\alpha,\alpha}(L_1) < \infty$. For $L_2$, we set
$K^+_{\lambda}(x,z)=H_{\lambda}(x,z)+ \widetilde{H}_{\lambda}(x,z)$,
where $$ H_{\lambda}(x,z)=\frac{1}{(2 \pi)^{\nu}} \int \frac{e^{i(z-x) \cdot
  \xi} d \xi}{M(\xi') - \lambda}$$
We begin by proving the first claim for $H_{\lambda}(x,z)$. Let $$
p(x)=\int (x-z)^{\beta}N(D_{z'})H_{\lambda}(x,z)\varphi(z)dz$$
where $\varphi \in C^{\infty}_0(\R^{\nu})$. Assume the
partial convolution $\varphi*'f(x)$ defined by $$\int
\varphi(z',x'')f(x'-z')dz'$$ then $$p(x)= {\varphi}^{\vee}*' \Big[
(-x')^{\beta'}N(D_{z'})\F^{\textrm{-1}'}_{\textit{x}} \big[
\frac{\textrm{1}}{\textsl{M}(\xi')-\lambda} \big] \Big]*''\Big[
(-\textit{x}'')^{\beta''}\delta_{\textit{0}} \Big]$$ Considering
the partial Fourier transforms $$\F' \textit{p}
(\xi',\textit{x}'')=\big(\F'\varphi^{\vee} \big)
\textsl{D}^{\beta'}_{\xi'} \Big[
\frac{\textsl{N}(\xi')}{\textsl{M}(\xi')-\lambda} \Big]*'' \Big[
(-\textit{x}'')^{\beta''}\delta_{\textit{0}} \Big]$$ $$\F''
\textit{p} (\textit{x}',\xi'')=(\F'' {\varphi}^{\vee})*'
\Big[(-\textit{x}')^{\beta'}\textsl{N}(\textsl{D}_{\textit{z}'})\F^{\textrm{-1}'}_{\textit{x}}
\big[ \frac{\textrm{1}}{\textsl{M}(\xi')-\lambda} \big] \Big]
\Big[ \textsl{D}^{\beta''}_{\xi''} \textrm{1} \Big]$$ we have $$
\F \textit{p}(\xi) =\F\varphi(-\xi)
\Big[\textsl{D}^{\beta'}_{\xi'} \big[
\frac{\textsl{N}({\xi'})}{\textsl{M}({\xi'})-{\lambda}} \big]
\Big]$$ and the outer bracket is estimated by Nilsson to be
$O(1)\mid \lambda \mid^{-c}$ as $\lambda \rightarrow - \infty$.
For the composed mapping $L_2 \circ L_1(\varphi)$, in order to
make the norm $M_{\alpha}$ finite, we apply a weight operator in
the ''bad'' variable. This weight is here (through duality), acting
as an operator with constant coefficients, so $p(x)$ is replaced
by $$ \big[\textsl{D}^{\beta}b \varphi \big]^{\vee} *'
\Big[(-x')^{\beta'}N(\textsl{D}_{z'})
\F^{\textrm{-1}'}_{\textit{x}} \big[
\frac{\textrm{1}}{\textsl{M}(\xi')-\lambda} \big] \Big]*''
\Big[(\textrm{1}-
\Delta_{\textit{x}''})^{-\textsl{N}}\delta_{\textit{0}} \Big]$$ by
choosing $N > \mid \alpha'' \mid$ and using that $I_{\big[
\beta_{\lambda}*''_{\textit{x}} \U_{-\textit{2}\textsl{N}}
\big]}( \varphi )(x)= I_{\big[ \beta_{\lambda}
\big]}(\varphi)*''_{\textit{x}}
\U_{-\textit{2}\textsl{N}}(\textit{x})$, we have proved that
$N^{\alpha,\alpha}(\big[(\textsl{D}^{\beta}b)(x-z)^{\beta}N(\textsl{D}_{z'})H_{\lambda}\big]*''\U_{-\textit{2N}}(\textit{x,z}))
= \textsl{O}(\textrm{1})\mid \lambda \mid^{-\textit{c}}$ as
$\lambda \rightarrow - \infty$. Consider now
$\widetilde{H}_{\lambda}(x,z)$, note that it is sufficient to
consider the case when $x \in K$,  $$\widetilde{H}_{\lambda}(x,z)= \frac{1}{(2
\pi)^{\nu}} \int \underbrace{\Big[
  \frac{1}{P^x(\xi') - \lambda} - \frac{1}{M(\xi')-\lambda} \Big]}_{G_{\lambda}(x,\xi')}
  e^{i(z-x) \cdot \xi} d \xi$$
With much the same arguments as in the previous case we get for
$$ q(x)= \varphi^{\vee}*' \Big[
(-x')^{\beta'}N(\textsl{D}_{z'})\F^{\textrm{-1}'}_{\textit{x}}\textsl{G}_{\lambda}(\textit{x},\xi')
\Big]*'' \Big[ (-\textit{x}'')^{\beta''}\delta_{\textit{0}} \Big]$$
that
$$ \F \textit{q} (\xi)= \F \varphi(-\xi) \Big[ \textsl{D}^{\beta'} \big(
\textsl{N}(\xi') \textsl{G}_{\lambda}(\textit{x},\xi') \big) \Big]$$
where the bracket according to Nilsson is $O(1)\mid \lambda \mid^{-c}$
as $\lambda \rightarrow - \infty$ so by applying the same weight as
above we get
$$ N^{\alpha,\alpha}(L_2 \circ L_1*''\U_{-\textit{2} \textsl{N}})=\textsl{O}(\textrm{1})\mid \lambda \mid^{-\textit{c}} \qquad
\lambda \rightarrow - \infty$$
Proof of (\ref{est2}):\\
Now for the second claim, we can write
$R_{\lambda}(x,z)=\sum_{\lambda,\mid \beta \mid=k}S_{\lambda,\beta} $ (a finite sum),
where each term is on the form
$S_{\lambda,\beta}=F_{\beta}(x,z)\widetilde{G}_{\lambda,\beta}(x,z)$. The
conditions on the coefficients give that $F_{\beta} \in
C^{\infty}(\R^{\nu} \times \R^{\nu})$ and with bounded
derivatives. Thus $R_{\lambda}(x,z)=0$ for $\mid x \mid,
\mid z \mid$ large. Using the tensor form
$K^+_{\lambda}(x,z)=h_{\lambda}(x',z') \otimes \delta_{x''}(z'')$ we
choose in this case to give a separate estimation of
$\widetilde{G}_{\lambda}$ in the ''good'' variable. According to
Nilsson
$$ \mid \textsl{D}^{\gamma'}_{x'}((x'-z')^{\beta'}N(\textsl{D}_{z'})h_{\lambda}(x',z'))
\mid \leq C \mid \lambda \mid^{-c}(1 + \mid x'-z' \mid)^{-2(n+1)}$$
for $\mid \gamma' \mid \leq n+1+\mid \alpha \mid$. Using the inequality
$$  (\Omega) \qquad \mid x-z \mid^{\alpha} \leq C \mid x \mid^{\alpha} \mid z
\mid^{\alpha}$$ we get $$ \mid \textsl{D}^{\gamma'}_{x'}\F_{\textit{x}}
(\textsl{R}_{\lambda}(\textit{x}',\textit{z}') \otimes \delta_{\textit{x}''}) \mid \leq \text{C} \mid \lambda
\mid^{-\textit{c}}(\textrm{1} + \mid \textit{x}' \mid)^{-\textit{(n+1)}}$$
Let $r(x)=\int R_{\lambda}(x,z)\varphi(z)dz$
then $$\mid \textsl{D}^{\gamma'}_{x'}r(x) \mid \leq C \mid \lambda \mid^{-c} (1+
\mid x' \mid )^{-(n+1)} M_{\alpha}(\varphi)$$ where we have used that
$\B_{\alpha} \subset \B_{\textit{0}}$. Finally,
$$ M_{\alpha'}(r)=\int \sum_{\gamma' \leq \alpha'} \mid
\widehat{\textsl{D}^{\gamma'}r} \mid dx \leq C \mid \lambda
\mid^{-c}M_{\alpha}(\varphi)$$ for $\lambda$ large and negative. If we
apply the same weight operator as before
we get $ N^{\alpha,\alpha}(r*\U_{-2N}) =\textsl{O}(1) \mid \lambda
\mid^{-c}$ as $\lambda \rightarrow - \infty$. This proves the
second claim and we have proved the lemma.$\Box$
\newtheorem{Prop1,Lem2}[Prop4]{Lemma  }
\begin{Prop1,Lem2} \label{L3}
  For $\beta_{\lambda}(x,z)$ in (\ref{beta-lambda}), we have
  $$ N^{\alpha,\alpha}(\exp(\kappa \mid \lambda \mid^b(z_j' -
    x_j'))I_{\beta_{\lambda,N,t}}(x,z))=\textsl{O}(1) \mid \lambda \mid^{-c} \quad
    \lambda \rightarrow - \infty$$ with $j,\kappa,\alpha$ and $t$ as in the
    proposition.
\end{Prop1,Lem2}
Proof: follows immediately from the results in (\cite{Ni_72}).

\vsp

Proof of the proposition:
We now consider the distribution (with compact support)
$$r_{\lambda}(x,y)= \sum_j b_j(z) R_j(D_{z'}) Q_j(D_{z''})K^+_{\lambda}(x,z)$$
where $b_j$ are the $C^{\infty}$-coefficients corresponding to $L$,
$R_j$ are constant coefficients operators
strictly weaker than $M$ in $x'$, $Q_j$ are constant coefficients
operators and $K^+_{\lambda}$ is the ''parametrix'' as before. We start by
assuming that $\kappa=0$. We make the
same approach as in the proof of Lemma \ref{L2}, noting that the Taylor series now contains a
term $b_j(x)R_j(D_{z'})h_{\lambda}(x',z') \otimes Q_j(D_{z''})\delta_{x''}(z'')$.
But we still have that $b(x) \in \B_{\alpha}$, so
$N^{\alpha,\alpha}(L_1) < \infty$ (if $b(x)$ is constant there is
nothing to prove).

\vsp

For the term corresponding to $H_{\lambda}(x,z)$, we get with the
same calculations as in the proof of (\ref{est1}) $$ \F
\textit{p}(\xi)=\F\varphi(-\xi) \Big[\textsl{D}^{\beta'}
\big[ \frac{\textsl{R}(\xi')}{\textsl{M}(\xi')-\lambda} \big]
\Big] \Big[ \textsl{D}^{\beta''} \textsl{Q}(\xi'') \Big]$$ for
some test function $\varphi$. The last bracket can be handled by
applying a weight operator of (''Sobolev'')-order $-2t$, where $t
> deg_{\R^{\textit{m}}} Q + \mid \alpha'' \mid$. For the
first bracket, Nilssons result still holds (according to the
conditions on $R$, also for the term where $\beta=0$), that is $$
\mbox{ess. sup } \mid \textsl{D}^{\beta'} \Big[
\frac{R(\xi')}{M(\xi')-\lambda} \Big] \Big[
\frac{\textsl{D}^{\beta''}Q(\xi'')}{(1+ \mid \xi'' \mid^2)^t}
\Big] \mid = \textsl{O}(1)\mid \lambda \mid^{-c} \quad \lambda
\rightarrow - \infty$$ For the second term corresponding to
$\widetilde{H}_{\lambda}(x,z)$, we get the same results. So we conclude with this
modification $$ N^{\alpha,\alpha}(L_2 \circ
L_1*''_{x}\U_{-\textit{2N}})=\textsl{O}(\textrm{1})\mid
\lambda \mid^{-\textit{c}} \qquad \lambda \rightarrow - \infty$$
Finally, for the analogue to (\ref{est2}), we estimate the
''good'' variable separately, $$ \mid \textsl{D}^{\gamma'}_{x'}
\mathcal{F'}_{\textit{x}}\textsl{R}_{\lambda}(\textit{x}',\textit{z}')
\mid \leq \textsl{C} \mid \lambda \mid^{-\textit{c}} (\textrm{1} +
\mid \xi' \mid^{\textrm{2}} )^{-\textit{(n+1)}}$$ We replace
$r(x)$ by $$ \tilde{r}(x)= \int R_{\lambda}(x',z') \otimes \big(
(1 - \Delta_{z''})^{-t} Q(D_{z''}) \delta_{x''}(z'') \big)
\varphi(z)dz $$
 with $t > deg_{\R^{\textit{m}}}Q + \mid \alpha'' \mid$ and by use of the inversion formula, we reach
 an estimate
 $$ \mid \textsl{D}^{\gamma'}_{x'}\tilde{r}(x) \mid \leq C' \mid \lambda
 \mid^{-c}(1 + \mid \xi' \mid^2 )^{-(n+1)} M_{\alpha}(\varphi)$$
as before. We get the conclusion with an additional weight
operator ($N > \mid \alpha'' \mid$) acting on the $x''$-variable, that is
$$ N^{\alpha, \alpha}(\tilde{r}*''\U_{-\textit{2N}})= \textsl{O}(\textrm{1})\mid \lambda
\mid^{-\textit{c}} \qquad \lambda \rightarrow - \infty$$
Finally, the case $k \neq 0$ can be treated  exactly as in the proof
of Lemma \ref{L3}, if we
substitute $R$ for $N$ and $M$ for $P^x$ in the right side of
$(\Theta)$, allowing $\beta'=0$, since Lemma 1 in \cite{Ni_72} can be applied
without modifications  on the strictly weaker operators.$\Box$

\vsp

{ $\bf{Remark:}$\\ Concerning  the translations, we are content with the
fact, that in the direction of the ''bad'' variable, translation
is in general more difficult to handle, however we know that (Mizohata \cite{Miz_repr}
Cor. 2) the mapping $$ \e \times
\textsl{H}^{\textsl{0},-\textsl{N}}  \ni
(\textit{c(z)},\textsl{T}) \rightarrow \textit{c(z)}\textsl{T} \in
\textsl{H}^{\textsl{0},-\textsl{N}}$$ is continuous with the
additional assumption that $c(x)=0$. Further $\parallel c(z)T
\parallel_{0,-N} \leq \epsilon \parallel T \parallel_{0,-N}$,
where $\epsilon$ can be chosen arbitrarily small, as the support
for $T$ tends to $\{ x \}$. Thus $$ \parallel c_j(z'')Q(D_{z''})
\varphi \parallel_{0,-N} \leq \epsilon
\parallel Q(D_{z''})\varphi \parallel_{0,-N} \leq \epsilon C \parallel
\varphi \parallel_{0,-N'}$$
We note that $\epsilon$ can be chosen as
$$ \epsilon= \Big( \mbox{sup}_{\mbox{supp } \varphi} \mid c(z) \mid + d(
\mbox{supp }
  \varphi,{x}) \Big)$$ where $d$ denotes the greatest distance between
  $x$ and $\mbox{supp } \varphi$.
We see that $\mid c_j(z'') \mid =\mid e^{i \gamma (z_j''-x_j'')}-1 \mid \leq 2$ for
$\gamma$ real and that $ \mid c_j(z'') \mid \rightarrow 0, \quad z_j''
\rightarrow  x_j''$, $1 \leq j \leq m$.
Since the functions $(1 + \mid \xi'' \mid^2 )^{-N}$ are very sensitive to
complex translations, we must here assume that $N=0$.
This implies (study the mollifier $\varphi_k=\exp(ix
\cdot \xi) \psi( x/k) / k^{n/2}$, as $k \rightarrow \infty$ for $\psi \in C^{\infty}_0$ and
$\parallel \psi \parallel_{0,0}=1$)
for all $ \xi'' \in \R^{\textit{m}}$, for $\beta'' \neq 0$ and
$d_Q=\mbox{ deg }_{\R^{\textit{m}}} Q$
$$ \mid D^{\beta''}Q(\xi'' - \gamma e_j'') - D^{\beta''}Q(\xi'') \mid
\leq \epsilon C( 1 + \mid \xi'' \mid^2 )^{d_Q \mid \beta'' \mid}$$
Through the analytical properties of
$Q(\xi'')$, we have that the
inequalities can be extended to $\gamma = -i \kappa \mid \lambda
\mid^b$,
for $\kappa$ real with $\mid \kappa \mid \leq \kappa_0$ (although we have to assume $ \lambda  \neq \pm \infty$).}
Note that a better result can be found in \cite{Ho_LPDO} (Lemma 3.1.5, (1963) ).

\vsp

However, since $\alpha_{\lambda}$ only involves the Dirac measure on
the ''bad'' side, it is trivial that the proposition implies that, for every test function $\varphi$
    $ \in C^{\infty}_0(\R^{\textit{m}})$
    $$N^{\alpha,\alpha}(\exp({\kappa}\mid \lambda
    \mid^b(z_j-x_j))I_{\alpha_{\lambda}*_{\textit{x}}''\varphi)}=\textsl{O}(\textrm{1})\mid \lambda
    \mid^{-\textsl{c}} \qquad \lambda \rightarrow -\infty$$ with $1 \leq j \leq
    \nu$, for all $\alpha$ and for $\kappa$ real with $\mid k \mid \leq k_0$.

\vsp

    We also need an estimate of $K^+_{\lambda}(x,z)$. We introduce
    the integrals
    $$ T^{2\alpha}_t(\lambda) = \int \frac{\xi^{2\alpha} d \xi}
    {(M(\xi') - \lambda)(1 + \mid \xi'' \mid^2)^t} $$
    According to \cite{Ni_72} Lemma 9, we have:\\

    \vsp

    \newtheorem{Prop5}[Prop4]{Proposition }
    \begin{Prop5} \label{Prop5}
         There is a positive constant $\kappa_0$ such that, when $1 \leq j \leq n$
         and $\alpha$ is any multi-index, we have for all real numbers $\kappa$ with
         $\mid \kappa \mid \leq \kappa_0$ and for some positive t
         \begin{equation} M_{\alpha}(\text{exp}(\kappa \mid \lambda \mid^b (y_j' - x_j'))
         K^+_{\lambda}*_y''\U_{-\textit{2t}}(\textit{x,y}))=\textsl{O}(\textrm{1})\textsl{T}^{\textsl{2} \alpha}_{\textit{t}} (\lambda) \end{equation}
         $\lambda \rightarrow - \infty$, where the norm $M_{\alpha}$
         is taken with
         respect to the variable x and where the estimate is uniform with respect to
         $ y \in \R^{\nu}$. Further b is the positive number corresponding to M
         as (\ref{b}).

         \vsp

         When $\beta \leq \alpha$, we also have, uniformly in $y
         \in \R^{\nu}$, for some positive $N$,
         \begin{equation} M_{\alpha}(D^{\beta}_y K^+_{\lambda}*_{\textit{x}}''\U_{-\textit{2N}}(\cdot,\textit{y}))=
         \textsl{O}(\textrm{1})\textsl{T}^{\textsl{2} \alpha}_{\textsl{N}}(\lambda) \end{equation}
         $\lambda \rightarrow - \infty$.
         \end{Prop5}

\section{ The complex set of lineality}
\label{sec:Lin}
A constant coefficients operator $P(D)$, is called reduced, if for
$\eta \in \mathbf{C}^{\textit{n}}$,
$$ P(\xi + t \eta) - P(\xi)=0  \qquad \mbox{ for all } \xi \mbox{ in } \mathbf{C}^{\textit{n}} \mbox{ and all } \textit{t} \in \R \quad \Rightarrow \eta = \textsl{0} $$
 We have use for the class
$\mathcal{H}_{\sigma}$ of polynomials $P$, such that $\mid \xi \mid^{\sigma}
\leq C \mid P(\xi) \mid$ for large $\mid \xi \mid$, for some constant $C$.
Consider first, as before, complex translations in one variable, that is $\eta=i \eta_0 e_j''$, $1 \leq j \leq m$ , where $e_j''$
corresponds to a standard base in $\mathbf{C}$ $^m$ and $\eta_0$ is some real number. The set of $\eta''$, such that $L(\xi + t \eta'') - L(\xi)=0$, for our
partially hypoelliptic operator $L=P_0+R$, is a zero-dimensional analytic (algebraic) set, that can locally, be given as
the zero-set corresponding to a polynomial. That is (disregarding points $\xi$ where $R(\xi)=0$ and assuming $R$ dependent on all variables)
$$ \Delta(L) = \{ \eta'';\quad R(\xi',\xi'' + t \eta'') - R(\xi',\xi'')=0 \ \forall t, \ \forall \xi \} = Z_{\varphi}$$
where $\varphi(\eta'')$ is a polynomial in a complex variable. The set $\Delta(L)$ will consist of a finite number of isolated points, that
 may cluster at the boundary of the local domain given for $\eta''$. We could say that the operator $L$ is reduced,
 with respect to $\varphi$. Further, for any truly reduced operator $Q$, we have that
$$ \Delta(L)=\{ \eta''; \quad Q(\xi',\xi''+ \varphi(t \eta'')) - Q(\xi',\xi'')= 0 \}$$
Let's assume $Q \in \mathcal{H}_{\sigma}(\R^{\nu})$ with $\sigma > 0$ and hypoelliptic in $L^{2}$. We can also
assume, that a sufficiently large number of iterations of the operator $Q$, gives an operator hypoelliptic in
$\mathcal{D'}$. Thus, if $1 \leq j \leq m$ is arbitrary, we have that $L^N$ is hypoelliptic, with respect to $\varphi$.

\vsp

If we write the condition on reducedness,
$$ \big[ e^{i<x,t \eta''>} - 1 \big]R(D)=0 \quad \Rightarrow \eta'' = 0$$
and interpret the bracket as a hypoelliptic operator, dependent on a parameter $t,$ on $L^{2}$, using the arguments in section \ref{sec:Hyp-infty}, even though $R$,
is not reduced, we will assuming $\big[ e^{i<x,t \eta''>}-1 \big] R(D)=0$, get
$$ \big[ e^{i<x,t \eta''>}-1 \big]^N R^N(D)=0 $$
for some positive integer $N$. This means particularly that
$$R^N(\xi',\xi''+t \eta_1'')\pm \ldots \pm R^N(\xi',\xi''+t \eta_N'')=R^N(\xi',\xi'')$$
The Fredholm theory gives, that the null-space to the translation
operator, is stable from some iteration index $N$, so if $R^N$ is
reduced, we have $\big[e^{i<x,t \eta''>}-1 \big]R^N(D)=0 \Rightarrow
\eta''=0$.
\vsp

Assume now,
$$ ( 1 + c\mid \zeta \mid )^k \mid \varphi_{\lambda} ( \zeta ) \mid^2 = C h(\zeta) $$
for a holomorphic function $h \in \mathcal{O}($ $U')$,$U \subset U'$ an open set. We can assume $\mid h \mid \leq 1$ and the domain $U$ bounded.
Let $g( \zeta ) = C h( \zeta ) - ( 1 + \mid \zeta \mid )^k \mid \varphi_{\lambda} ( \zeta ) \mid^2$, be the holomorphic
function defining $U$. As $\mbox{ord}_{\zeta} h \leq \mbox{ord}_{\zeta} g$, for every $\zeta \in U$, we have that
$g $ is less than a constant times $h$ in modulus. Thus for m = 1, by Rouch$\acute{e}$'s theorem, $g$ is a polynomial. In higher
dimension, we can at least say that the restriction of $g$ to any complex line, is a
polynomial. Let's define a real set $$\Delta_{\mathbf{C}}=\{ \eta \in
\R^{\textit{m}}; \quad \textsl{R}(\xi + \textsl{it} \eta) - \textsl{R}(\xi)=\textsl{0} \quad  \xi \in
\R^{\nu}, \quad  \textsl{t} \in \R \}$$

\newtheorem{Treves}{Lemma }[section]
\begin{Treves}
If $R_{\lambda}$ is reduced for complex lineality $\Delta_{\mathbf{C}}$ (with respect to one dimensional translations),
 then $(1+ \mid \xi \mid)^c \leq C \mid R_{\lambda}(\xi) \mid$, for positive constants $c,C$
and $\xi \in \R^{\nu}$.
\end{Treves}

Proof: ( sketch of a proof )
$R_{\lambda}$ can, as an entire function, be developed as \begin{equation} R(\xi + i t \eta)
- R(\xi)=\sum_{\alpha} \big[ R_{\alpha}(i t \eta) - R_{\alpha}(0)
\big] \xi^{\alpha}= \sum_{\alpha} F_{\alpha}(i t \eta)
\xi^{\alpha} \end{equation}
We can without loss of generality, assume the following argument
takes place outside the complex set of zero's to the polynomial
$R$.
If $\eta \in \Delta_{\mathbf{C}}$, $i \eta$ is a growth vector for
$R$, $\mbox{ ord }_0 F_{\alpha}=+ \infty$ and we can show $\mid
F_{\alpha}(i t \eta) \mid / \mid t \eta \mid^r \leq C_{r,t}$, for
every positive real $r$, for every multiorder $\alpha$ and for $t$ close to $0$. Assume
$\Delta_{\mathbf{C}}=\{ 0 \}$, it then follows (\cite{Chirka},Ch.1, section 5) that
for $\eta \neq 0$
\begin{equation} \label{Chirka} \mid
F_{\alpha}(i t \eta) \mid / \mid t \eta \mid^r > C_{t,\alpha} \quad \text{
for some real } r
\end{equation}
for $t$ close to $0$. Further $ \mid R_{\alpha}( i t \eta )
\mid / \mid t \eta \mid^r > C_{\alpha}$, for $t$ close to $0$.
This means, that $$\mid \xi \mid^v C_{t,\alpha,\eta} \leq
\sum_{\beta \leq \alpha,\mid \alpha \mid = v} \mid R_{\alpha}(i
t \eta) \xi^{\beta} \mid$$
Thus, for $\eta
\notin \Delta_{\mathbf{C}}=\{ 0 \}$,
 particularly, $\mid \xi \mid^v \leq C
\mid R_{\lambda}(\xi) \mid$, for some $\lambda$. Note that $v$ can be chosen as the least $v$ for which $F_{\alpha}\neq 0$, $\mid \alpha \mid =v.\Box$

\vsp

For $h(\zeta)=0$, since $g^t$ is reduced, for some positive
integer $t$, there is a positive $\sigma$, such that $\mid \xi
\mid^{\sigma} \leq \mid g(\xi) \mid$, for large $ \xi \in $
$\R^{\textit{m}}$. Thus, $\mid \xi \mid^{t \sigma -k''} \leq C
\mid \varphi_{\lambda}(\xi) \mid^{2t} $, for some positive $k''$.
For $\mid h(\zeta) \mid > 0$, the inequality for $h$, is
immediate. Assuming the defining polynomial is self-adjoint, we
conclude that $\varphi^t_{\lambda} \in \mathcal{H}$ $_{\frac{1}{2}(t
\sigma - k'')}(\R^{\textit{m}})$. This means, for a
sufficient number of iterations, $\sigma > k''/t$ and the polynomial is reduced. The constant $k'' < 1$, so there is a
positive integer $N$ (possibly smaller than $t$), such that
$\sigma > 1/N$ and $\varphi^N_{\lambda}$ is hypoelliptic in $L^{2}$.
Note that this hypoellipticity is not dependent on $\lambda$.

\vsp

 Obviously, if $\mid \xi \mid^{\sigma} \leq
C\mid \varphi_{\lambda}(\xi) \mid$, for some positive $\sigma$, then the same must hold for $R_{\lambda}$. Since, if $\eta
\in \Delta_{\mathbf{C}}$, the $R_{\alpha}$'s are constants and if $\eta \notin \Delta_{\mathbf{C}}$, the result follows
analogously to the second part of the proof of the previous lemma. Note that for a hyperbolic polynomial, this
implication is not true. Since the set of lineality for such an operator is determined by a homogeneous polynomial
$Q$, if $\varphi_{\lambda}$ is reduced, we get something like $\mid \xi \mid^{- \sigma} \leq C \mid Q(\xi) \mid$,
for a positive number $\sigma$.

\vsp

Assume now, $\mbox{ Re }P$ partially hypoelliptic. This is a self-adjoint operator, which means
that there is a $N_0$ such that $(\mbox{ Re }P)^N$ is hypoelliptic for $N \geq N_0$. The same holds
for the imaginary part. Assume $\mbox{ Im }P \prec \mbox{ Re }P$, then $i(\mbox{ Im }P)^{N'} +
(\mbox{ Re }P)^{N'} \prec P^{N'}$ and we can conclude that $P^{N'}$ is hypoelliptic. Thus the requirement
on self-adjointness in the above argument is not essential in the proof of

\newtheorem{Ruc}[Treves]{ Proposition }
\begin{Ruc} \label{Ruc}
If $P(D)$ is a constant coefficients, (self-adjoint), partially hypoelliptic operator, there is an iteration index
$N_0$, such that $P(D)^N$ is hypoelliptic in $\D$, for all $N \geq N_0$.
\end{Ruc}

\vsp

Let $e^{i<x,\varrho>}$ denote the translation operator in $\parallel \cdot \parallel$, for $\varrho \in \R^{\textit{n}}$, we then
have a remainder operator $R(D)$ with constant coefficients dependent on the coefficients for $P(D)$ and on $\varrho$, such that
$$e^{i<x,\varrho>}R(D)=e^{i<x,\varrho>}P(D)-P(D)e^{i<x,\varrho>}$$
Further,
$$\big[ e^{i<x,\varrho>}-1 \big]P(D)=e^{i<x,\varrho>}R(D)+P(D) \big[
e^{i<x,\varrho>} -1 \big]$$
and we can prove, that there is a $\lambda_0$ large, such that
{\small \begin{equation} \label{HE-in-infty} \parallel e^{i<x_j,\mid \lambda \mid>} R(D)u \parallel \leq
C_{\lambda}
\big( \parallel P(D)u \parallel + \parallel u \parallel \big) \qquad \mid \lambda \mid > \lambda_0 \quad u \in
H^{0,0}_K \end{equation}}
where the constant can be chosen, so that $C_{\lambda} \rightarrow 0$ as $\mid \lambda \mid \rightarrow \infty$.
 For a reduced operator, the distance function grows like
 $\mid \xi \mid$ and we have that $R \prec \prec P$, where the constant may depend on $\lambda$.
For a non-reduced operator, we note that we may
assume $\mid \xi \mid^q \leq \mid \lambda \mid$, as
$\mid \lambda \mid \rightarrow \infty$, since with the opposite condition, the result is trivial.
More precisely, assume $V_{\lambda}=\{ \xi; \quad \mid \lambda \mid \leq  \mid \xi \mid^q \}$. These
sets will become insignificant as $\mid \lambda \mid \rightarrow \infty$, however we prove that $R \prec
\prec P$ on this set, as $\mid \lambda \mid \rightarrow \infty$. The inequality $\mid \xi \mid^{2q}
\leq \mid \lambda \mid d(\xi, \Delta_{\mathbf{C}})$ gives that $\mid \lambda \mid \leq C
d(\xi,\Delta_{\mathbf{C}})$ as $\mid \lambda \mid \rightarrow \infty$. Further $\mid \lambda
\mid^{\sigma} \mid \xi \mid^{\sigma} \leq C d(\xi,\Delta_{\mathbf{C}})$, so
$$\frac{\mid R(\xi) \mid}{1 + \mid P(\xi)\mid} \leq C \frac{1}{\mid \lambda \mid^{\sigma}} \rightarrow 0 \ \mbox{ as } \mid
\lambda \mid \rightarrow \infty$$
This gives, with the first condition, for large $\mid \xi \mid$, $\frac{1}{c}\mid \xi \mid^{\sigma-q} \mid R(\xi +  \mid \lambda \mid e_j)
\mid\leq \mid \lambda \mid + \mid R(\xi) \mid \mid \xi \mid^{q-\sigma}$, with $q$ according to (\cite{Palamadov},Ch.II,§3,Prop.2) and $\sigma$
corresponds to the distance function to the lineality (one dimensional), for the operator $R$. The result
follows, for $\sigma > q$.

\section{The Levi parametrix method}
\label{sec:Levi}
The method used to construct a fundamental solution $g_{\lambda}(x,y)$ to the operator
$L(y,D_y) - \lambda$, when $\lambda$ large and negative, is a modified
version of Levi's parametrix method \cite{John}. Assume $K^+_{\lambda}(x,z)$ a fundamental
solution in $\mathcal{S}'$ to the operator with constant coefficients $P^x_{\lambda}(D_{z'})=P(x,D_{z'}) - \lambda$,
$P^x$ hypoelliptic in $z'$, that is the operator $P(z,D_{z'})$ in section \ref{sec:Miz_repr}
(1) frozen in $x$. This fundamental solution was discussed in section \ref{sec:Fund_sol},
\begin{displaymath}
   K^+_{\lambda}(x,z)= \left \{
\begin{array}{lr}
 h_{\lambda}(x',z')\otimes\delta_{x''}(z'') \qquad x \in K \\
 t_{\lambda}(x',z')\otimes\delta_{x''}(z'') \qquad x \notin K
\end{array} \right.
\end{displaymath}
where $t_{\lambda}$ is the fundamental solution to the type operator regarded
as an operator on $\R^{\textit{n}}$. Note that $t_{\lambda}(x',\cdot),h_{\lambda}(x',\cdot)$ in
$C^{\infty}(\R^{\textit{n}} \backslash$ ${x'})$.
Let
\begin{equation} g_{\lambda}(x,z) = K_{\lambda}^+(x,z) + \Big [u_{\lambda},K_{\lambda}^+ \Big ](x,z) \label{fsol}
\end{equation} where the bracket stands for $\Big [f,g \Big ](x,z)=
\int f(x,y)g(y,z) dy$. Assuming $u_{\lambda}$ can be constructed with certain
given regularity properties, $g_{\lambda}(x,z)$ will be the
fundamental solution to the operator $L_{\lambda}(z,D_z)=L(z,D_z) - \lambda$, that is
\begin{equation} \overline{v(x)} = \int \overline{L_{\lambda}(z,D_z)}
 \overline{\text{[}v(z)\text{]}}g_{\lambda}(x,z)dz \label{gfsl} \end{equation}
Using that $K^+_{\lambda}(x,y)$ is a fundamental solution to
$P^x_{\lambda}(D_{z'})$ and assuming that $v \in \s$ (although it would
be sufficient to assume $v \in \e$)
\begin{equation} \int \overline{L_{\lambda}(z,D_z)} \overline{\text{[} v(z) \text{]}}K^+_{\lambda}(x,z)dz  = \int
(\overline{L_{\lambda}(z,D_{z}}) - \overline{P^x_{\lambda}(D_{z'}})) \overline{\text{[} v(z) \text{]}}
  K^+_{\lambda}(x,z) dz + \label{green} \end{equation}
  $$ \int \overline{P^x_{\lambda}(D_{z'})} \overline{\text{[} v(z) \text{]}}
  K^+_{\lambda}(x,z) dz = - \int \overline{v(z)} \alpha_{\lambda}(x,z) dz + \overline{v(x)}$$
where $$\alpha_{\lambda}(x,z)=(P^x_{\lambda}(D_{z'})-L_{\lambda}(z,D_z))K^+_{\lambda}(x,z)=
\delta_x(z) - L_{\lambda}(z,D_z)K^+_{\lambda}(x,z)$$ Note that
$K^+_{\lambda}$ acts as a ''parametrix'' also to the perturbed operator. We have
for $x' \neq z'$, $-\alpha_{\lambda}(x,\cdot) \in
C^{\infty}$ in $z'$ and for $\mid \beta' \mid < M$, ($M< \varrho - n,
\varrho$ according to section \ref{sec:Sob_emb}) $D^{\beta'}_{x'}\alpha_{\lambda}(\cdot,z)
\in C^0(\R^{\textit{n}})$. Using the expression (\ref{fsol}) in (\ref{gfsl}) and the resulting
equality in
(\ref{green}) and after reversing the order of integration in one of
the integrals, we get
\begin{equation}
 0 = \int \overline{v(z)} \Big( u_{\lambda}(x,z) - \alpha_{\lambda}(x,z) - \Big[
u_{\lambda},\alpha_{\lambda} \Big](x,z) \Big)dz \label{fred-alt}
\end{equation}
The bracket $\Big[ f,g \Big](x,y)$, for $f,g \in
\B_{\textit{0}}$, can be regarded
as a kernel to an operator K $\in \mathcal{L}(\B_{\textit{0}})$, (here $\B_{\textit{0}}$ is the
space of tempered distributions with locally summable Fourier
transforms), such that $K(w)= \int \Big[f,g
\Big](x,z)w(z,y)dz$. Also, the kernel itself can be regarded as an
operator in $\mathcal{L}(\B_{\textit{0}})$.
 In order to construct $g_{\lambda}$, we need to prove that, after
modification, $\alpha_{\lambda}$ is bounded as an operator on $\B_{\alpha}$.

\subsection{ The remainder ${\alpha}_{\lambda}$ weighted is in ${\B}_{\alpha}$}
A closer study of $\alpha_{\lambda}(x,z)$ gives first that $K^+_{\lambda}(x,\cdot)$ has partial Fourier
transforms  in $L^1_{loc}(\R^{\textit{m}})$ and $L^1(\R^{\textit{n}})$
respectively. Further the coefficients of
$$L_{\lambda}(z,D_z)=P_{\lambda}(z,D_{z'}) + R(z,D_z)$$ according to (1), have
derivatives with compact support, which means that if they are not
constant, they are in $\B_{\alpha}$, so
$\F_x(P_{\lambda}(z,D_{z'})K^+_{\lambda}(x,z))(\xi) \in
L^1_{loc}$ and partially in $L^1(\R^{\textit{n}})$. For the part of $\alpha_{\lambda}(x,z)$
involving derivatives in the ''bad'' variables, we get that
\begin{eqnarray}
  R(z,D_z)K^+_{\lambda}(x,z)=
\sum_jP_j(z,D_{z'})Q_j(D_{z''})K^+_{\lambda}(x,z)= \nonumber \\
\sum_{j,l}c_{j,l}(z) \Big( N_{j,l}(D_{z'})h_{\lambda}(x',z') \otimes
Q_j(D_{z''})\delta_{x''}(z'') \Big) \end{eqnarray}
so the conditions, $P_j
 \prec \! \prec {}_{z'} M$ and the coefficients, if not constant in $\B_{\alpha}$, give that the Fourier
transform acting on $z$ and $RK^+_{\lambda}$, is in $L^1_{loc}(\R^{\nu})$ and
partially in $L^1(\R^{\textit{n}})$.

\vsp

For the parameter $x$, we have that $$t_{\lambda}\otimes
\delta_{x''}(x,z)= \Big[ \frac{1}{(2 \pi)^n}\int \frac{\exp \Big[
i(z'-x') \cdot \xi' \Big] d \xi'}{M(\xi') - \lambda} \Big] \otimes \delta_{x''}(z'')$$ so that $\mathcal{F'}_{z}t_{\lambda}(x',z') \in
L^{1}({\R}^{\textit{n}})$. The same results follow for $h_{\lambda}(x',z')$, from the
conditions on the operator $P^x$. So
$K^+_{\lambda}*''_x\U_{-\textit{2t}},$ $(LK^+_{\lambda})*''_{\textit{x}}\U_{-{\textit{2t}}} \in$ $
\B_{\alpha}$, for some $t$.

\vsp

$\bf{Remark:}$ It would be possible to construct, for the contact operators, a fundamental solution in the space of
tempered distributions with Fourier transform in
$L^1_{loc}(\R^{\nu})$. We will however prefer a construction, in the space of analytic functionals.

\vsp

Let
$F_{\lambda}(x,z)$ be the expression within the parenthesis in
(\ref{fred-alt}),
so that $\int \overline{v(z)}F_{\lambda}(x,z) dz = 0 $, for all test functions
$v \in \mathcal{O}$. For a fixed $x$, we then have $F_{\lambda}=0$.
This means that the problem of finding the fundamental solution
$g_{\lambda} \in H'$ (both variables), is
reduced to finding $u_{\lambda}$ such that in $H'(\mathbf{C}^{\nu}_{\text{z}})$
\begin{equation} \alpha_{\lambda}(x,z) = u_{\lambda}(x,z) -
  \Big[u_{\lambda},\alpha_{\lambda} \Big](x,z) \label{ueq} \end{equation}
and such that in $H'(\mathbf{C}^{\nu}_{\text{x}})$
$$ \alpha_{\lambda}(x,z) = u_{\lambda}(x,z) -
  \Big[\alpha_{\lambda},u_{\lambda} \Big](x,z) $$

\subsection{ The remainder $\widetilde{\alpha}_{\lambda}$ is in $\B_{\alpha}$ }
Using a partial mollifier, we see that $\Big[ f,g \Big]$ is a continuous
operator in both variables, on the space of tempered distributions with Fourier
transforms in $L^1(\R^{\textit{n}})$ and
$L^1_{loc}(\R^{\textit{m}})$ respectively, assuming $f,g$ have the
same properties.
The distribution
$\alpha_{{\lambda},N,t}=\alpha_{\lambda}*_x''\U_{\textit{-2}
  \textsl{N}}*_{\textit{y}}''\U_{-\textit{2t}}(\textit{x,y})$,
for $t > \mbox{max deg }Q_j + m + \mid \alpha'' \mid$ and $N > m$,
is from the construction a linear, closed operator on the space
$\B_{\alpha}$. This implies continuity by the closed graph
theorem.
We will now use the double partial regularization meaning
$$ \alpha_{\lambda}*_x''\varphi*_y''\psi(x,y) =
\int\!\!\!\int
\alpha_{\lambda}(x',x''-z'',y',y''-w'')\varphi(z'')\psi(w'')dz''dw''
$$
The continuity implies that
$\widetilde{\alpha}_{\lambda}=\alpha_{\lambda}*_x''\varphi*_y''\psi$
is in $\mathcal{L}(\B_{\alpha})$ (and using the previous paragraph,
it is also in $\B_{\alpha}$), since for $\varphi \in C^{\infty}_0(\R^{\textit{m}})$,
we have $( 1 + \mid \xi'' \mid^2)^N\widehat{\varphi} \in L^{\infty}$
(we denote the corresponding norm $\parallel \cdot \parallel_{\infty,N}$) for every integer $N$.
Thus
$$ M_{\alpha}( \alpha_{\lambda}*_x''\varphi*_y''\U_{-\textit{2t}})
\leq \parallel \delta_{\textit{0}} \otimes \varphi
\parallel_{\infty,\textsl{N}} \textsl{M}_{\alpha}((\textrm{1}
- \Delta_{\textit{x}''})^{-\textsl{N}}\alpha_{\lambda}*_{\textit{y}}''\U_{-\textit{2t}})$$
In the same way
$$M_{\alpha}( \alpha_{\lambda}*_{\textit{x}}''\U_{-\textit{2}\textsl{N}}*_{\textit{y}}''\psi)
\leq \parallel \delta_{\textit{0}} \otimes \psi \parallel_{\infty,\textit{t}} \textsl{M}_{\alpha}( (\textrm{1}
- \Delta_{\textit{y}''})^{-\textit{t}} \alpha_{\lambda}*_{\textit{x}}''\U_{-\textit{2}\textsl{N}})$$
Particularly $\widetilde{\alpha}_{\lambda} \in \B_{\alpha}$ in both variables.

\subsection{ Convergence for the series $u_{\lambda}$ in $\E(\R^{\nu} \times \R^{\nu})$}
\label{sec:construction}
If in the coefficients for $L$, the variable $y''$ is regarded as a
parameter, we can write ${\alpha}_{\lambda}(x,y)$ on the form
\begin{eqnarray} \sum_j
\Big(P^x_j(D_{y'})-P_{j,(y'')}(y',D_{y'}) \Big) h_{\lambda}(x',y') \otimes
Q_j(D_{y''})\delta_{x''}(y'') = \nonumber \\  \sum_jp_{\lambda,j,(y'')}(x',y') \otimes
q_j(x'',y'') & & \nonumber \end{eqnarray}
Proposition
\ref{Prop4} implies an estimate of the regularization, for
$\varphi,\psi \in C^{\infty}_0(\R^{\textit{m}})$
$\alpha_{\lambda,(\sim)}(x,y)=\sum_j \big( p_{\lambda,j,(y'')} \otimes
q_j*''_x \varphi*''_y \psi \big) (x,y)$, $N^{\alpha,\alpha}(I_{\alpha_{\lambda,(\sim)}})=O(1)\mid
\lambda \mid^{-c}$ as $\lambda \rightarrow - \infty$
(and it is not difficult to establish the same estimate for the norms
taken with respect to $x$).
For
$\lambda$ sufficiently large and negative, this operator norm is $<
1$. This makes it possible
to exclude the parameter $y''$ from the calculations,
since according to the conditions on the coefficients,
$\mid c_j(z',y'') - c_j(z',z'') \mid \leq C_1$ for
$y'',z'' \in K$ and $\leq C_2$ for $z'' \notin K$, where we have
constant coefficients. That is, for $\lambda$ sufficiently large,
\begin{displaymath}
    M_{\alpha}\Big( \int p_{\lambda,j,(y'')}(x',z')\otimes
  \widetilde{q}_j(x'',z'') \varphi(z) d z \Big) \leq \end{displaymath}
\begin{displaymath}
 M_{\alpha}\Big( (\int_K + \int_{\R^{\nu} \backslash \textit{K}})( p_{\lambda,j,(y'')} - p_{\lambda,j,(z'')})(x',z') \otimes
\widetilde{q}_j(x'',z'')\varphi (z) d z \Big) + \end{displaymath}
\begin{displaymath}
M_{\alpha}\Big( \int
p_{\lambda,j,(z'')}(x',z') \otimes \widetilde{q}_j(x'',z'')\varphi(z)d z
\Big) <  M_{\alpha}(\varphi)
\end{displaymath}
We will use the same symbol $\alpha_{\lambda,(\sim)}$ for the case
where the parameter $y''$ is ''frozen'' or excluded from the
calculations. Note that, this means that the requirement that the
operator $L_{\lambda}$ is formally partially self adjoint, is not
necessary in the estimates of the fundamental solution, that we shall
give.
Let
$u_{\lambda,(\sim)}$ denote the solution to the equation (\ref{ueq})
corresponding to $\alpha_{\lambda,(\sim)}$. This is an integral equation
of Fredholm type and a theorem from the theory of Fredholm operators on Banach
spaces (\cite{Ho_FnAnalys} Lemma 2.5.4) gives existence of
the solution $u_{\lambda,(\sim)}$, not dependent on the multi-index, in
$\mathcal{L}(\B_{\alpha})$ that is, for $\lambda$ sufficiently large
and negative
\begin{equation} u_{\lambda,(\sim)}= \alpha_{\lambda,(\sim)}+ \Big[
\alpha_{\lambda,(\sim)},\alpha_{\lambda,(\sim)} \Big] + \Big[
\alpha_{\lambda,(\sim)},\alpha_{\lambda,(\sim)} \Big]_2 + \ldots
\label{gseries} \end{equation}
where $\Big [f,f
\Big]_2$ stands for $\Big [\Big [f,f \Big ] ,f \Big ]$ with convergence in
$\mathcal{L}(\B_{\alpha})$. Note that for
$I_{u_{\lambda,(\sim)}}:\B_{\alpha}(\R^{\nu}_{\text{y}})
\rightarrow \B_{\alpha}(\R^{\nu}_{\text{x}})$, when $x$ is
assumed fixed, we get that $I_{u_{\lambda,(\sim)}} \in
  \B_{\alpha}^{'}(\R^{\nu}_{\text{y}}) \subset
  \mathcal{D'}(\R^{\nu}_{\text{y}})$ and as usual, we identify
  $I_{u_{\lambda,(\sim)}}$ with the kernel $u_{\lambda,(\sim)}$.
If we let $\widetilde{q}_j$ be the regularized elements, that is $q_j*''_x\varphi*''_y\psi$, the series ($\ref{gseries}$) becomes
\begin{eqnarray} \sum_j p_{\lambda,j,(y'')} \otimes \widetilde{q}_j + \sum_{j,k} \Big[ p_{\lambda,j,(y'')},p_{\lambda,k,(y'')} \Big]' \Big[
\widetilde{q}_j ,\widetilde{q}_k \Big]'' + \nonumber \\ \sum_{j,k,l} \Big[ \Big[ p_{\lambda,j,(y'')},p_{\lambda,k,(y'')} \Big], p_{\lambda,l,(y'')}
\Big]'\Big[ \Big[ \widetilde{q}_j ,\widetilde{q}_k  \Big]'',
\widetilde{q}_l \Big]'' +..\nonumber \end{eqnarray}
where the brackets are taken over $\R^{\textit{n}}$ and
$\R^{\textit{m}}$ respectively. Let $I$ be the index set $\{ i_0,i_1, \cdots,i_N, \cdots
\}$, where every index assumes values in $\{1, \cdots,r \}$. For the
  brackets over $\R^{\textit{n}}$, let for instance $ \Big[ p_{\lambda,i_0,(y'')},p_{\lambda,i_2,(y'')}
  \Big]_2=$ $ \Big[ \Big[ p_{\lambda,i_0,(y'')},p_{\lambda,i_1,(y'')}
  \Big],p_{\lambda,i_2,(y'')} \Big]$. The remainder, in the series (\ref{gseries}) corresponding to
$u_{\lambda,(\sim)}$ is then on the form
\begin{equation}
  \sum^{\infty}_{\mid I \mid= N+ 1} \sum_I  C_{i_0} \ldots C_{i_{\mid I \mid}}\Big[
p_{\lambda,i_0,(y'')},p_{\lambda,i_{\mid I \mid + 1},(y'')} \Big]'_{\mid I \mid +
    1}(x',y') \varphi \otimes \psi_{i_{\mid I \mid + 1}}(x'',y'') \label{regularized-u-lambda}
  \end{equation}
   where $C_i= \int \psi_i(z'')
  \varphi(z'') d z''$ and where
  $\psi_i(z'')=Q_i(D_{z''})\psi(z'')$.  Note
  that for every index, $\psi_{i}$ involves an application of an operator of
  finite order. The product of constants can be estimated through the estimated
  behavior of the ''good'' side, for $\lambda$ sufficiently large and negative.
  That is, we have $\mid R_{\lambda,(\sim)} \mid \leq const. e^Ce^{-\kappa
    \mid \lambda \mid^b}$, where $C_{i_j} < C$ for every $i_j$.
  Note that if the test functions are assumed to have support in a
  neighborhood of $0$, then the convolution is defined, for $y''$ in a neighborhood of $x''$. This is sufficient for the study of
  the fundamental solution.

\vsp

  We claim that we have convergence
  $u_{\lambda,(\sim)}(x,\cdot) \rightarrow u_{\lambda}(x,\cdot)$ in
  $\E$, as
  $\varphi \otimes \psi \rightarrow \delta_{x''}$. Using the equation
  following (\ref{ueq}) and the parallel argument, the claimed
  convergence can be established in both variables separately.
  Let $u_{\lambda,N,(\sim)}$ be the $N+1$ first terms in
(\ref{gseries}), then
  immediately
{\small \begin{equation}  u_{\lambda,N,(\sim)} - \alpha_{\lambda,(\sim)}= \Big[
\alpha_{\lambda,(\sim)},\alpha_{\lambda,(\sim)} \Big]_1 + \ldots + \Big[
\alpha_{\lambda,(\sim)},\alpha_{\lambda,(\sim)} \Big]_N = \Big[
u_{\lambda,N-1,(\sim)}, \alpha_{\lambda,(\sim)} \Big] \label{u-lambda-series}
\end{equation}} \par

We have established convergence as $N \rightarrow \infty$, for the left
side in (\ref{u-lambda-series}). According to (\cite{Ni_72}, Lemma 10),
the ''good'' part in $u_{\lambda,N,(\sim)}$, that is the '-bracket in
(\ref{regularized-u-lambda}) (with the sum taken from 1 to $N$) is in
$C^0(\R^{\textit{n}} \times \R^{\textit{n}})$. Noting that these partial sums have
compact support, we get that  $u_{\lambda,N}(x,\cdot) \in
\E(\R^{\nu})$. In section \ref{sec:Paley-Wiener}, we argue that this means that $u_{\lambda}(x,\cdot)$ is a linear
form and a distribution in $\E(\R^{\nu})$ (Banach-Steinhaus theorem).

\subsection{ Regularity properties for the series $u_{\lambda}$ }
\label{sec:Reg-u-lambda}
In this section we argue that even though $u_{\lambda}$ will
contain differential operators of infinite order, it is for a
study of the regularity properties, sufficient to study a finite
number of terms, in the development of $u_{\lambda}$, that is the
action of finite order differential operators.
We start by noting the following trivial proposition.

\emph{
Assume $K \in \E(\R \times \R)$ the
kernel to a continuous integral operator, $I_K:\E \rightarrow \E$ (or $B_{\alpha} \rightarrow B_{\alpha}$),
such that $\mid \xi \mid^{\sigma}
\leq \mid \F_{\textit{y}}$ $K(x,y) \mid ,\quad \mid \xi \mid$
large, for $\sigma >0$.
Then $M_{\alpha}(u) \leq M_{\textsl{0}}(I_K(u))$, $\mid
\alpha \mid \leq \sigma/2$, for $u$ a distribution, such that the
right side is finite.}

The proof is immediate. Thus, such an operator is hypoelliptic on
$\E$ (or $B_{\alpha}$). We could choose $K(x,y)=P(D)\delta_{x}(y)$, for a reduced
differential operator with constant coefficients. Further, we
have the following proposition.

\emph{
Assume $K$ as in the previous proposition and $P(D)$ a constant coefficients
differential operator with $\mid \xi \mid^{\nu} \leq C \mid P(\xi)
\mid ,\quad \mid \xi \mid \text{ large }$ and some positive constant $C$. Then, $\mid \xi
\mid^{\nu+\sigma} \leq C \mid \F_{\textit{y}}$ $P(D_x)K(x,y) \mid ,\quad
\mid \xi \mid \text{ large }$ and $\nu > 0$.}

This means, that also $I_{P(D)K}$ is hypoelliptic on $\E$ (or $B_{\alpha}$).
From the section \ref{sec:Lin}, follows that the iterated
polynomial $Q_{i_j}^N$, is reduced with $\sigma \geq 1$, for sufficiently large
$N$. This behavior will be stable for further iteration. Thus,
the remainder operator, as in (\ref{u-lambda-in-B-alpha}), is
hypoelliptic in $\E$ (or $B_{\alpha}$), particularly
$$ I_{R_{\lambda,N}}(u) \in C^{\infty} \Leftrightarrow u \in
C^{\infty}$$
This does however not mean, that it can be represented by a  very regular
kernel. That is, if $I_K$ is hypoelliptic, then $\mbox{ sing supp }(I_K(u))=\mbox{ sing
supp }(u)$ and $ \mbox{ sing supp }(u) = \mbox{ sing supp }(I_K(u)-u)$. If $K$ is very regular, $I_K(u)-u=\gamma \in
C^{\infty}$, where we can assume $\gamma$ has nontrivial support.
The latter proposition is stronger: Assume for instance $K$ a
fundamental solution to a hypoelliptic operator $P$. If $K$ is assumed hypoelliptic, then $
\mbox{ sing supp }(u - Pu) = \mbox{ sing supp }(u)$, for every $u \in \D$ and
if the propositions were equivalent, we would have
$P-I$ is regularizing.
For $\varphi \in \E$ (or $B_{\alpha}$), $\mbox{ sing supp
}(I_{u_{\lambda}}-I_{u_{\lambda,N}}(\varphi))=$ $\mbox{ sing supp
}(I_{\alpha_{\lambda}}(\varphi)) \subset \mbox{ sing supp
}(I_{K^+_{\lambda}}(\varphi))$ and $\mbox{ sing supp
}(I_{g_{\lambda}-K^+_{\lambda}-u_{\lambda,N}}(\varphi))=\mbox{
sing supp }(I_{K^+_{\lambda}}(\varphi))$. Thus, to study
propagation of singularities, in the formal solution
$g_{\lambda}$, we only have to study a finite number of terms, the
remainder will have a hypoelliptic action. If $K$ is the kernel
corresponding to a hypoelliptic integral operator $\E
\rightarrow \E$ and if $U \in \E($ $\Omega \times
\Omega)$, for a bounded open set $\Omega \subset \R^{\nu}$, then $I_U: C^{\infty} \rightarrow \E$.
Further, the composition is defined $I_K(I_U(\varphi))(x)=I_{\big[ K,U
\big]}(\varphi)(x)$, for $\varphi \in C^{\infty}$. Let $U=\delta_{x}$,
then using nuclearity for $\E($ $\Omega)$, we have $\big[
K,\delta_{x} \big] \in \E($ $\Omega \times \Omega)$ and $x \in
\Omega$. Thus, to prove that $R_{\lambda,N}$ is a kernel
in $\E$, it is sufficient to prove that the corresponding
integral operator maps $\E \rightarrow \E$. We will
however use a slightly different approach.

\subsection{ Some remarks on the generalized Paley-Wiener theorem }
\label{sec:Paley-Wiener}
  For $K(x,y) \in {\e}^{'(\textsl{0})}$ $(X \times Y)$ and $X,Y$ open sets in $\R^{\nu}$, $\mu \in
{\E}^{(\textsl{0})}$ $(Y)$ we have
$$\int K(x,y) \mu(d y) \in {\E}^{( \textsl{0})}(\textsl{X})$$
For $K_1,K_2 \in {\E}^{( \textsl{0})}$ $(X \times Y)$, further
{\small $$ I_{K_1}(I_{K_2}(\mu))=\int K_1(x,y) \big( \int K_2(y,z)  \mu(d z)
\big)(d y)=\int \int K_1(x,y)K_2(y,z)dy  \mu(d z)$$}
this iteration can be repeated infinitely many times, so for $N \geq
0$ $$ I^N_K : {\E}^{(\textsl{0})}(\textsl{Y}) \rightarrow {\E}^{(\textsl{0})}(\textsl{X})$$
Particularly, if $K=F \otimes G$, with $F,G \in {\E}^{(\textsl{0})}$ in $X$
and $Y$ respectively, we have $I_K(\mu)(x)=CF(x)$, with $C=G(\mu)$.

\vsp

Using (\ref{regularized-u-lambda}) and section (\ref{sec:S-kernels}),
{\small$$I_{u_{\lambda,N,(\sim)}}(x,y)=\sum^N_{\mid I \mid=1,I}\big[,\big]'_{\mid I
    \mid + 1}(x',y') \varphi \times \ldots \times \varphi(Q_{i_0} \psi,\ldots,
  Q_{i_{\mid I \mid}} \psi) I_{Q_{i_{\mid I \mid +1}} \widetilde{\delta_{x''}}}(x'',y'')$$}
Now, $\varphi$ can be regarded as a measure in ${\E}^{(\textsl{0})}$ and
$\mid < \varphi , Q_{i_j} \psi > \mid \leq C \sup_K \mid Q_{i_j} \psi \mid$
for all $j$ and a compact set $K$. We can find a constant $C$, such
that
$\mid \mathcal{F''}$ $I_{u_{\lambda,N,(\sim)}} \mid \leq e^C$ and this estimate still
holds if we let $\varphi \rightarrow \delta_0$. Let's write
$\mathcal{F''}$ $I_{u_{\lambda,N,\psi}}$ for this limit. Using the results from
section \ref{sec:transl}, $\F$($I_{u_{\lambda,N,\psi}})$=$ O(1)e^{C}$, for $\lambda$
large and negative. Using the
terminology for the generalized Paley-Wiener theorem, we can say that \\
$\F$( $I_{u_{\lambda,N,\psi}}$) is an entire analytic function of
exponential type. Thus \\ $\mid \F$
($I_{u_{\lambda,N,\psi}})(\zeta) \mid \leq e^C$, $\zeta \in
\mathbf{C}^{\nu}$ and that $I_{u_{\lambda,N,\psi}} \in
{\E}^{(\textsl{0})}(\R^{\nu})$. We note that the same
estimate holds, as $N \rightarrow \infty$. We assumed above that $x$
was fixed, but a completely analogous argument holds, if we assume
instead that $y$ is fixed. We then have that
$\F$ ($I_{u_{\lambda,N,\varphi}}$) is of exponential type and
$I_{u_{\lambda,N,\varphi}} \in {\E}^{(\textsl{0})}(\R^{\nu})$.

\vsp

As $\psi \rightarrow \delta_0$ (or analogously if $\varphi
\rightarrow \delta_0$), we have convergence in $\mathcal{O}_M$ for
the Fourier transform and
we write $\F$($I_{u_{\lambda,N}}$) for the limit. We can prove an
  estimate for the Fourier-Laplace transform
  $$ \mid \F(\textsl{I}_{\textit{u}_{\lambda,\textsl{N}}})(\zeta) \mid \leq \textsl{C} ( \textsl{1} + \mid
    \zeta \mid )^{\textsl{M}}
    \textit{e}^{\textsl{H}_{\textsl{K}}(\mbox{ Im } \zeta)} \qquad
    \zeta \in \mathbf{C}^{\nu}$$
for all $N$, where $M \geq 0$ is dependent on $N$, $K$ is a compact, convex subset of
$\R^{\nu}$, $H_K$ the support function for $K$. We say
that $\F$( $I_{u_{\lambda,N}} )$ $\in (\mbox{ Exp } \mathcal{O}_M)$. From
Paley-Wiener-Schwartz theorem ((\cite{Ho_LPDO}, part I)
, Theorem 7.3.1), it follows that
$I_{u_{\lambda,N}} \in {\E}(\R^{\nu})$. We prove in section \ref{sec:constr_g_lambda},
that $I_{(u_{\lambda}-u_{\lambda,N})^{\sim}} \rightarrow 0$, as $N \rightarrow \infty$. So, using
  Banach-Steinhaus theorem, it follows
  that the limit as $N \rightarrow \infty$, $I_{u_{\lambda}} \in
  {\E}(\R^{\nu})$. Using appropriate test functions, we have
  \begin{equation} u_{\lambda}=\sum^N_{\mid I \mid = 1} \sum_{I} \Big[ p_{\lambda,i_0,(y'')},p_{\lambda,i_{\mid I \mid +1},(y'')}
  \Big]'_{\mid I \mid + 1} \otimes Q_{i_0} \cdots Q_{i_{\mid I \mid +
      1}} \delta_{x''} + R_{\lambda}  \label{u-lambda-in-B-alpha} \end{equation}
Since the convergence is established in the Fr\'echet space ${\E}$,
in both variables separately, it follows that we have convergence in
$\E(\R^{\nu} \times \R^{\nu})$.

\subsection{ Construction of fundamental solution to $L_{\lambda}$ in $H'(\mathbf{C}^{\nu} \times \mathbf{C}^{\nu})$ }
\label{sec:constr_g_lambda}

We have already seen (section \ref{sec:construction}) that, as an iteration operator
$\alpha_{\lambda,(\sim)}=O(1)e^{-\kappa {\mid
  \lambda \mid}^b}$, uniformly on compact sets, as $\lambda \rightarrow -\infty$
and this holds for every iteration, that is
{\small $$ N^{\alpha,{\alpha}}(I_{{\alpha}_{\lambda,(\sim)}}^N)
\leq C e^{-\kappa {\mid \lambda
  \mid}^b} N^{\alpha,{\alpha}}(I_{
{\alpha}_{\lambda,(\sim)}}^{N-1}) \leq
C' e^{-\kappa (N-1) {\mid \lambda
  \mid}^b} N^{\alpha,{\alpha}}(I_{{\alpha}_{\lambda,(\sim)}})
  \leq C'' e^{-\kappa N {\mid \lambda \mid}^b}$$}
and we see that $\mid I^N_{{\alpha}_{\lambda,(\sim)}} \mid
\rightarrow 0$ as $N \rightarrow \infty$, for large negative
$\lambda$. The behavior in the ''bad'' variable, has been
estimated to
$e^C$, for some positive constant $C$. Thus the remainder
$I_{R_{\lambda,M,(\sim)}}=O(1)e^{C-\kappa M \mid \lambda
  \mid^b} \rightarrow 0$, uniformly on compact sets, as $M \rightarrow
\infty$ in $\B_{\alpha}$, for $\lambda$ sufficiently large and negative. Further
$ \mid \mid I_{R_{\lambda,M,(\sim)}} \mid - \mid I_{R_{\lambda,M}^{\sim}} \mid \mid \leq C(M,K)$
for some positive constant $C(M,K)$, such that $C(M,K) \rightarrow 0$ as $M \rightarrow \infty$.
That is,$$ \mid \mid I_{R_{\lambda,M,(\sim)}} \mid - \mid I_{R_{\lambda,M}^{\sim}} \mid \mid \leq C_1
 \sum_{N \geq M} C^N \mid \big[ , \big]_{N+1}' \mid + B_N \mid \big[, \big]_{N+1}' \mid$$
for some positive constants $C_1,C$ and $B_N$. The sum can be estimated with,
 $e^{-M \kappa \mid \lambda \mid^b } \sum_{N \geq M}( C^N+B_N )$
or with $e^{C- M \kappa \mid \lambda \mid^b} \rightarrow 0$ as $M \rightarrow \infty$.
 Thus $\mid I_{R_{\lambda,M}^{\sim}} \mid \rightarrow 0$,
as $M \rightarrow \infty$.
We claim that the fundamental solution $g_{\lambda}$ is on the form
$g_{\lambda}=K^+_{\lambda} + \big[ u_{\lambda},K^+_{\lambda} \big]$
and to prove that this representation actually holds in
$H'$, we prove that the remainder in (\ref{g-lambda}),
$R^{\sharp}_{\lambda,N} \rightarrow 0$ in this space.

\vsp

For a (real) or complex vector space $E$, we can define a Banach space, $ Exp_{\alpha,\rho}(E)$, over E with
$$ \parallel G \parallel_{\alpha,\rho}=\sup_{z \in E} \mid G(z)e^{-\alpha \rho(z)} \mid$$
where $\alpha \in \R$ ${}^+$ and $\rho$ is a (real) or complex norm on $E$. We define $ Exp(E)$ as the
inductive limit of the spaces $ Exp_{\alpha,\rho}(E)$. We note the following \cite{Mart} Ch. 2, Proposition 1.5

\emph{
The Fourier-Borel transform establishes a topological vector space isomorphism between the convolution algebra $H'(E)$
and the algebra of entire functions of exponential type on $E^*$, $Exp(E^*)$.}

Assume $R_{\lambda,M}^{\sharp}$ is represented in $H'$, as a
measure $\mu_{M,\epsilon}$, with support in a ball with respect to
a complex norm $\rho$ and radius $(1+ \epsilon)\alpha$. It follows
from the estimates above, that $ \sup_z \mid e^{-C}
\mathcal{F}$ $I_{R_{\lambda,M,(\sim)}} \mid < \infty$, for a
positive constant $C$, for every $M$. Here
$\mathcal{F}$ denotes the Fourier-Borel transform. The
Paley-Wiener-Schwartz theorem, gives $\sup_z \mid e^{-C \mid z
\mid} \mathcal{F}$ $I_{R_{\lambda,M}} \mid < \infty$, for all $M$ and
we see that $\mathcal{F}$ $R_{\lambda,M} \in \mbox{ Exp }_{C,\mid z
\mid}$. The constant can be chosen as $C-\kappa M \mid \lambda
\mid^b$ and by proposition the cited result in section (\ref{sec:constr_g_lambda}), we see that
$I_{R_{\lambda,M}} \rightarrow 0$ in $H'$ as $M \rightarrow
\infty$. Obviously, $\F$ $R_{\lambda,M} \rightarrow 0$ in $\mbox{ Exp }_{C,\mid \xi \mid}$,
that is with a real norm, as $M \rightarrow \infty$. We can now assume $\parallel \mu_{M,\epsilon} \parallel
\rightarrow 0$, as $M \rightarrow \infty$. Thus
$ \mathcal{F}$ $R_{\lambda,M}^{\sharp}(\zeta)=\int
e^{<\zeta,\hat{z}>}d \mu_{M,\epsilon}$ and $\mid \mathcal{F}$ $R_{\lambda,M}^{\sharp}(\zeta) \mid \leq \parallel
\mu_{M,\epsilon} \parallel e^{(1+\epsilon)\alpha \rho^*(\zeta)}$, where
$\rho^*$ denotes the dual norm, that is $\rho^*(\zeta)=\sup_{\mid z
\mid \leq 1} \mid <z,\zeta> \mid$. We conclude that $\parallel
\mathcal{F}$ $R_{\lambda,M}^{\sharp} \parallel_{(1 +
\epsilon) \alpha, \rho^*} \rightarrow 0$, as $M \rightarrow
\infty$ and $g_{\lambda,M} \rightarrow g_{\lambda}$ in $H'$.

\subsection{ Error in the estimation }
\label{sec:Error}
  Finally, we need to estimate the difference
  $u_{\lambda,(\sim)}-\widetilde{u}_{\lambda}$, where
  $\widetilde{u}_{\lambda}$ is $u_{\lambda}*''_x\varphi*''_y\psi$ and
  $u_{\lambda}$ according to (\ref{u-lambda-in-B-alpha}). Using the
  representation (\ref{u-lambda-in-B-alpha}), we see that every term is on the
  form $\Big[, \Big]'\otimes \text{operator}\delta_{x''}$ so we could
  put $u_{\lambda}(x,y)=p_{\lambda,(y'')} \otimes q\delta_{x''}(x,y)
  + R_{\lambda}(x,y)$. Thus $\widetilde{u}_{\lambda}- \widetilde{R}_{\lambda} = (p_{\lambda,(y'')} \otimes q \delta_{x''})*''_x \varphi *''_y \psi =
  (p_{\lambda,(y'')} \otimes q(D_{y''}) \delta_{x''})*_{x,y}''(\varphi\psi)(x,y)=
  {\varphi}(x'')\Big( p_{\lambda,(y'')}q(-D_{z''})
  {\psi}(y''-z'') \Big) \vert_{z''=0}$, for $x \in \R^{\nu}$. The conditions on the coefficients gives that
  $\mid u_{\lambda,(\sim),M} - \widetilde{u}_{\lambda,M} \mid \leq
  C(M,K)$ with $K$ as in the beginning of the section and this
  relation holds for every $M$.
Further, in $H'$,
\begin{equation}
  g_{\lambda}(x,y)= h_{\lambda}\otimes \delta_{x''}(x,y) +
\Big[ p_{\lambda,(y'')},h_{\lambda} \Big]'\otimes q\delta_{x''}(x,y) +
R^{\sharp}_{\lambda}(x,y) \label{g-lambda}
\end{equation}

We have
$u_{\lambda,M+1}=\alpha_{\lambda} + \Big[
u_{\lambda,M},\alpha_{\lambda} \Big]$. If $g_{\lambda,M}$ is derived from
(\ref{fsol}) with $u_{\lambda,M}$ instead of $u_{\lambda}$, then
$g_{\lambda}-g_{\lambda,M}= \Big[
u_{\lambda}-u_{\lambda,M},K^+_{\lambda} \Big]$ and if we compare
(\ref{u-lambda-in-B-alpha}) with
(\ref{g-lambda}), $R^{\sharp}_{\lambda,M}= \Big[
R_{\lambda,M},K^+_{\lambda} \Big]$.
For the remainder terms, we get $\mid \Big[ R_{\lambda,(\sim),M}, K^+_{\lambda,(\sim)} \Big] -
\Big[ R_{\lambda,M}, K^+_{\lambda} \Big]^{\sim} \mid $ $\leq C_1(M,K)$, with $C_1(M,K) \rightarrow 0$ as $M \rightarrow
\infty$. A calculation similar to the one for $u_{\lambda,(\sim)}$ above, gives $\mid g_{\lambda,(\sim),M} -
\widetilde{g}_{\lambda,M} \mid \leq C_2(M,K)$. We conclude that the
difference $g_{\lambda,(\sim)} - \widetilde{g}_{\lambda}=O(1)$ as
$\lambda \rightarrow - \infty$.

\subsection{ Construction of fundamental solution to ${P(y,D}_{\textit{y}'})
- \lambda$ }
\label{sec:K-lambda-est}
The same method can be used to construct a fundamental solution
$K_{\lambda}(x,y)=k_{\lambda}\otimes \delta_{x''}(x,y)$ to the
variable coefficients operator $P(y,D_{y'}) - {\lambda}$. The same
arguments gives $$ K_{\lambda}(x,y)=K^+_{\lambda}(x,y) + \Big[
u^0_{\lambda},K^+_{\lambda} \Big](x,y)$$ with $u^0_{\lambda}$ such that
$$\beta_{\lambda}=u^0_{\lambda} - \Big[\beta_{\lambda},u^0_{\lambda} \Big]$$
The Neumann series for $u^0_{\lambda}$ is particularly simple. Let
$\beta_{\lambda}(x,y)=$ $\sum_j b_{\lambda,j}(x',y') \otimes
\delta_{x''}(y'')$ then
$$ u^0_{\lambda,t}= \sum^{\infty}_{\mid I \mid = 1} \sum_I
\Big[b_{\lambda,i_0},b_{\lambda,i_{\mid I \mid +1}} \Big]'_{\mid I \mid
    + 1} \otimes \delta_{x''}*_x\U_{-2t}*_y\U_{-2t}$$
  with convergence in $\mathcal{L}(\B_{\alpha})$. If we let $B$ denote
  the sum, since $(1- \Delta_x'')^t(1-\Delta_y'')^tu_{\lambda,t}^0(x,y)=u_{\lambda}^0(x,y)$, we get $$K_{\lambda}(x,y)= \Big( h_{\lambda} + \Big[
  B,h_{\lambda} \Big]' \Big)(x',y') \otimes \delta_{x''}(y'')$$

\section{ Construction of a parametrix to the operator ${L_{\lambda}}$}
\label{sec:parametrix}
We have constructed a fundamental solution to the partially formally
hypoelliptic operator $L_{\lambda}$, using a modified ''parametrix''. We finally
prove that use of  parametrices corresponding to our
operator in the parametrix method would not give a
fundamental solution with better properties.

\vsp

Assume as before, that the singularity is in 0. It is in the
parametrix method sufficient to consider contact operators, frozen in
0. As before $L_{\lambda}=P_{\lambda}+R$, and let $R^{\star}(z,D)=\sum
\Big[ c_j(z) - c_j(z',0) \Big]R_j(D)$ for
$R(z,D)=\sum_jc_j(z)R_j(D)$. We now have $L^{\star}_{\lambda}=$ $\Big(
L^{\star}_{\lambda}-P_{\lambda}^{\Sigma} \Big) + P_{\lambda}^{\Sigma}=$ $A_{\lambda} +
P_{\lambda}^{\Sigma}$,
where $P_{\lambda}^{\Sigma}=P_{\lambda}(z',0,D')$.
Then $A_{\lambda}(z)=L_{\lambda}-L^{\Sigma}_{\lambda}$, where
$L^{\Sigma}_{\lambda}$ has constant coefficients on the ''bad'' side and
variable coefficients only for operators strictly weaker or equivalent to
$P_{\lambda}^{\Sigma}$. Further,
$A_{\lambda}(z',0)=0$.
Assume $K_{\lambda}^{\Sigma}$ a fundamental solution to the operator
$P_{\lambda}^{\Sigma}$. If $L_{\lambda}$ is assumed partially hypoelliptic, then $-\alpha_{\lambda}=A_{\lambda}K_{\lambda}^{\Sigma} \in
C^{\infty}(\R^{\nu})$, since if  $\Sigma=\{z ; z''=0 \}$
and $\mbox{ sing supp }K_{\lambda}^{\Sigma} \subseteq \Sigma$, then $\mbox{ sing supp }(R-R^{\Sigma})K_{\lambda}^{\Sigma}
\subseteq \Sigma$ but $R-R^{\Sigma}=0$ in $\Sigma$, so
$(R-R^{\Sigma})K_{\lambda}^{\Sigma} \in C^{\infty}$. Thus $K_{\lambda}^{\Sigma}$ is a parametrix to $L^{\star}_{\lambda}$,
that is $L^{\star}_{\lambda}K_{\lambda}^{\Sigma}=\delta_0-\alpha_{\lambda}$.

\vsp

For the more general operator, we use the null-space to the remainder operator, ($R^{\star}$ above), in the
construction of a parametrix with singular support on $\Sigma$.
We can construct a fundamental solution to $L^0_{\lambda}=P^0_{\lambda}+\sum_{j=1}^{r}P_j(0,D')Q_j(D'')$
with singular support (and in fact support) on $\Sigma$, as follows. Let $E_{0,\lambda}$ be
the fundamental solution to $P_{\lambda}^0$ and $E_j, j \neq 0$,
solutions to the homogeneous equations on $L^{2}($ $\R^{\textit{n}})$ (or $\E)$,
$P_jE_j=0$. Let $E_{\lambda}(y)=\sum_{j=0}^m E_j(y')\otimes \delta_0(y'')$, then we can prove that
$L^0_{\lambda}E_{\lambda}=\delta_0$.
We prove the argument for $L^{\Sigma}_{\lambda}$, using (\cite{Ho_LPDO} Th. 13.3.3 (1983)). Assuming $P_j(z',D')E_j=0, j
> 0$, $P_{0,\lambda}(z',D')E_{0,\lambda}=\delta_0$,
we can find a linear mapping $\mathcal{L}_{\textit{k}} : \E \rightarrow \E$ ($L^{2} \rightarrow
L^{2}$),
such that for $j > 0
$ or $k > 0$,\\
$P_j(z',D')$ $\mathcal{L}_{\textit{k}}$ $P_k(z',D')E_j=$ $P_j(z',D')E_j=0$, for $E_j$
adjusted to a compact set.

\vsp

Assume $L_{\lambda}^{\Sigma}$ formally self-adjoint, $L_{\lambda}^{\Sigma}=P_{\lambda}^{\Sigma}+R^{\Sigma}$,
where $R^{\Sigma}$ is a
tensor product, of a variable coefficients operator
and a constant coefficients operator. Let
$E_{\lambda}=K_{\lambda}^{\Sigma}+H$, such that $P_{\lambda}^{\Sigma}K_{\lambda}^{\Sigma}=K^{\Sigma}_{\lambda}
P_{\lambda}^{\Sigma}=I$ and
$R^{\Sigma}H=HR^{\Sigma}=0$,
 on $L^{2}$ ($\E$). We then have $L_{\lambda}^{\Sigma} E_{\lambda}=I + P_{\lambda}^{\Sigma}H + R^{\Sigma}
 K_{\lambda}^{\Sigma}$.
We now claim that $R^{\Sigma}K_{\lambda}^{\Sigma}=P_{\lambda}^{\Sigma}H=0$ in $L^{2}$ ($\E$). Proof:
$R^{\Sigma}$ is surjective
on $L^{2}$ ($\E)$, that is for $f \in L^{2} (\E)$, there is a $\varphi \in L^{2}$ $(\E)$, such that
$R^{\Sigma} \varphi = f$. Further $L_{\lambda}^{\Sigma} \varphi=P_{\lambda}^{\Sigma} \varphi + f$ and $\varphi=\varphi +
HP^{\Sigma}_{\lambda} \varphi + K_{\lambda}^{\Sigma}f$, that is $K_{\lambda}^{\Sigma}f=-HP_{\lambda}^{\Sigma}\varphi$
and $R^{\Sigma}K^{\Sigma}_{\lambda}f=0$. The second identity is trivial
and we have $L^{\Sigma}_{\lambda}E_{\lambda}=I$.
Further, $L_{\lambda}
E_{\lambda}=L^{\Sigma}_{\lambda}E_{\lambda} +
A_{\lambda}E_{\lambda}=I +
\beta_{\lambda}$, where $\beta_{\lambda} \in C^{\infty}$ and we have
a parametrix to $L_{\lambda}$.

\vsp

We need to estimate the solutions $E_j$, to the homogeneous equations.
Let's consider the bounded mappings $E,E_0,E_1$ on $L^{2}(\R^{\textit{n}})$ associated to the operator
$L_{\lambda}$, according to
$ P_{\lambda}^0(D')E_0f = f$ for $f \in L^{2}$,
$ L_{\lambda}^{\Sigma}Ef=L_{\lambda}(z',0,D')Ef = f$ for $f \in L^{2}$ and
$ (P_{\lambda}^0(D')+R^{\Sigma}(z',D'))E_1f=f$ for $f \in L^{2}$.
Assume $g \in L^{2}$, such that $R^{\Sigma}(z',D')g=0$. For
$P^{\Sigma}_{\lambda}g=f$, we get $L_{\lambda}^{\Sigma}(z',D')(g - Ef) = 0$
and since $L_{\lambda}^{\Sigma}$ is hypoelliptic, $g=Ef+\eta$ in
$L^{2}$ with $\eta \in
C^{\infty}$, such that $L_{\lambda}^{\Sigma} \eta=0$. Here
$R^{\Sigma}$ denotes one of the strictly weaker operators in the
representations of $L^{\Sigma}_{\lambda}$ and $g$ corresponds to a
solution to the homogeneous equation, for this operator.

\vsp

Assume $f$ a $L^{2}-$ function such that $R^{\Sigma}Ef=0$, is $Ef \in C^{\infty}$?
Proof:
Immediately $P_{\lambda}^{\Sigma}Ef-f=0$, extending with $Ef$, we get
$(P_{\lambda}^{\Sigma}-I)Ef+(E-I)f=0$ and conclude that $Ef \in C^{\infty}$. Using that
$g=Ef + \eta$, $\eta \in C^{\infty}$, we see that $g \in C^{\infty}$.
Further, we can write $P^0_{\lambda}(E_1-E_0)f=-R^{\Sigma}E_1f$. We have earlier noted that if the left side is assumed in $H^{s}$
then the right side will be in $H^{s+\sigma}$, for some positive $\sigma$.
Repeating this procedure gives $P^0_{\lambda}(E_1-E_0)f \in C^{\infty}$. We can write $E_1=E_0B$ where $B=(I+A)^{-1}$ and $A=R^{\Sigma}E_0$.
Since $(I+A):C^{\infty} \rightarrow C^{\infty}$, we get
$-Af \in C^{\infty}$. Further $Bf \in C^{\infty}$ and finally $f \in
C^{\infty}$.

\vsp

We have earlier proved that
\begin{equation}
  \sup \mid D^{\beta'}_{\xi'}
\big[ \frac{R^0(\xi')}{M(\xi')-\lambda} \big] \mid = O(1) \mid
\lambda \mid^{-c} \qquad \lambda \rightarrow - \infty \label{par-est}
\end{equation}
for some positive number $c$. We wish to study the relevant
commutators over $\R^{\textit{n}}$. First, $E_0(y')$ is most
easily defined as convolution with a fundamental solution to
$P_{\lambda}^0(D')$ (compare \cite{Ho_LPDO}(1983) Th.
13.2.1). Let $C_{E_0}$ be the commutator $E_0 \psi - \psi E_0$,
$\psi \in C^{\infty}_0(\R^{\textit{m}})$, such that $\psi=1$
on some suitable set. A Taylor expansion of order $m$, for $\psi$
in $x$, gives $$M_{\alpha}( C_{E_0}f )\leq C \sup \mid
D^{\beta'}_{\xi'} \big[ \frac{1}{M(\xi')-\lambda} \big] \mid
M_{\alpha}(f)= O(1) \mid \lambda \mid^{-c} M_{\alpha}( f) $$ where
$\mid \beta' \mid \leq m$ and the norm is taken over
$\R^{\textit{n}}$. For the commutator $C_A$, the same type of
argument gives $N^{{\alpha},{\alpha}}(C_A)=O(1) \sup \mid
D^{\beta'}_{\xi'} \big[ \frac{R^{\Sigma}(\xi')}{M(\xi')-\lambda}
\big] \mid =$ $O(1) \mid \lambda \mid^{-c}$, as $\lambda  \rightarrow -\infty$, where $R^{\Sigma}$ is
assumed frozen in some point on $\Sigma$. For the iterated
operators, we have $C_{A^n}=A^{n-1}C_A +C_AA^{n-1}$. Using the
Banach algebra property of $\B_{\alpha}$, we can at least say
that $N^{{\alpha},{\alpha}}(C_B)=O(1)\mid \lambda \mid^{-c}$ as
$\lambda \rightarrow -\infty$.

\vsp

Assume $\psi$ a test function, such that $\psi=1$ on an open set
containing the singularity, then
$$\psi g =  E \psi f - C_E f + \psi \eta$$
The commutator $C_{E_1}$, can be rewritten $C_{E_0}B+E_0C_B$ and is
easily estimated in operator norm to $O(1)\mid \lambda \mid^{-c}$. For
the mapping $E$, we have the same form $E=F_0B$, where now $F_0$ is a
fundamental solution to the contact operator $L_{\lambda}^0$, but
since this operator is equivalent to $P_{\lambda}^0$, we get the same
estimate for $C_E$ and $E$. All that remains is to estimate $\eta$, a
solution to the homogeneous equation corresponding to
$P_{\lambda}^{\Sigma}$. But this follows immediately from the estimates we have
produced for the fundamental solutions to $P_{\lambda}^{\Sigma}$.
Expressions on the form (\ref{par-est}), can be treated as in the
proof of Lemma \ref{Prop4} and we see that
$$N^{{\alpha},{\alpha}}(\mbox{ exp }(\kappa \mid \lambda
\mid^b(z_j'-x_j'))C_A)=O(1)\mid \lambda \mid^{-c} \qquad \lambda
\rightarrow - \infty$$
with $j,\kappa$ and $\alpha$ as in proposition 5.1.
The commutators $C_{E_0},C_{E_1},C_E$, can be estimated in the same
way. Finally, we must have that $\eta=0$ in a neighborhood of the singularity,
so $\eta=O(1)e^{- \kappa \mid \lambda \mid^b}$, as $\lambda \rightarrow
-\infty$, uniformly on compact sets in $\R^{\textit{n}}$.
The fundamental solution $K_{\lambda}^{\Sigma}$ was estimated in
section \ref{sec:K-lambda-est}.

\newtheorem{He_par}{ Lemma }[section]
\begin{He_par} \label{He-par}
 If $E$ is parametrix to a differential operator $P$, such that for every $V$( =neighborhood $x$), $PE - \delta_{x} \equiv 0$ in $V$,
then $P$ is not hypoelliptic.
\end{He_par}

Proof: Assume $P$ hypoelliptic, with a parametrix $E$.  We then
have that $(I_{E}-I)$ is locally regularizing. If locally $I_{PE}=I$, then also locally
$u - Pu \in C^{\infty}$, for all $u \in \D$. But, since $P$
is hypoelliptic, the same must hold for $P-I$ and we have a
contradiction.$\Box$

\vsp

$\bf{Remark:}$ This problem is however easily handled. In the
constant coefficients case, we can assume $E$ a fundamental solution to the hypoelliptic operator $P$ and choose a
test function $\zeta \in \mathcal{D}$ such that $\mbox{ supp } \zeta \cap \mbox{ supp }E \neq \emptyset$.
 Assume further that $\zeta$=1 in $U_{\epsilon}$ an $\epsilon-$(neighborhood 0), such that the commutator $C_PE \neq 0$ in $\mbox{ supp }\zeta \backslash
U_{\epsilon}$, then $\zeta E$ is a parametrix to $P$ and $P\zeta E - \delta_{x} \neq 0$ in $\R^{\nu}$ $\backslash U_{\epsilon}$.
In the variable coefficients case, if $P$ is hypoelliptic, then $E$ is very regular. We can assume $\mbox{ deg }P >
0$,
so that $E$ has support in some (neighborhood of $0) \backslash \{ 0 \}$. Any such neighborhood will do. We use that $C_PE \in C^{\infty}$
and that $I_{PC_E}=-I_{C_PE} \in C^{\infty}$ and since $P$ is hypoelliptic, $C_E \in C^{\infty}$.

\vspace{.5cm}

However, we make a small modification of the parametrix constructed in this section, so that the
remainder is regularizing. Assume as before,
that $L_{\lambda}^{\Sigma}$ is the variable
coefficients operator with support on ${\Sigma}$ and
$A_{\lambda}=L_{\lambda}-L_{\lambda}^{\Sigma}$.
Consider, instead of $E_{\lambda}$, $K_{\lambda}^{\Sigma}+K_{\lambda}^{\delta}$, where
$K_{\lambda}^{\delta}=\sum_j E_j \otimes U_j^{\delta}$, for $U_j^{\delta}$ very regular distributions,
mapping $L^{2} \rightarrow L^{2}$ with support on
$\Sigma' \times U_{\delta}$, a neighborhood of $x''$ and with $E_j$, such that
$L_{\lambda}^{\Sigma}\sum_jE_j=0$.
We use a commutator, to modify the parametrix as before, $C_{K_{\lambda}^{\delta}}=
K_{\lambda}^{\delta} \psi - \psi K_{\lambda}^{\delta}$, for some suitable test function $\psi \in C^{\infty}_0$ and $E_{\lambda}^{\delta}=K_{\lambda}^{\Sigma}
+ C_{K^{\delta}_{\lambda}}$. We have that $L_{\lambda}E_{\lambda}^{\delta}=
L_{\lambda}^{\Sigma} K_{\lambda}^{\Sigma} + A_{\lambda}C_{K^{\delta}_{\lambda}}$. If in the
commutator, the test function is chosen with support on $\Sigma' \times U_{\delta}$ and $=1$ on
$V_{\epsilon} \times U_{\delta}$, for some neighborhood of $x'$,$V_{\epsilon} \subset \Sigma'$, it follows that
$A_{\lambda}C_{K^{\delta}_{\lambda}} \in C^{\infty}$. This is a parametrix in the usual sense,
with regularizing remainder.

\section{ Hypoellipticity in $L^{2}$ and $\mathcal{D'}$ }
\label{sec:L-2-HE-to}
We have discussed the parametrix construction in $L^{2}$ and $\mathcal{D'}$. Obviously, an operator homogeneously hypoelliptic in $\mathcal{D'}$,
must be homogeneously hypoelliptic in $L^{2}$. We give the following result in the opposite direction,
\newtheorem{L2_HE}[He_par]{ Proposition }
\begin{L2_HE}
Assume $P$ a variable coefficients, constant strength, differential operator, with a representation $P(x,D)=P_0(x,D) + \sum_j R_j(x,D)$,
such that $R_j \prec P_0$ and $\sum_j R_j$ with $\sigma > 0$ for all frozen operators. Further, that the operator is defined as a
constant coefficients, hypoelliptic type operator outside a compact set and that $\mbox{ Re P} \sim
P_0$, then $P$ is hypoelliptic in $\mathcal{D'}$.
\end{L2_HE}
Assume however first, that $P=P_0 +  \lambda R$, with $R \prec
\prec P_0$, a constant coefficients differential operator and that
$P_0$ has $\sigma >0$. It is trivial that, if
$\mathcal{N}$ $(P_0+\lambda I) \neq \{ 0 \}$, for some $\lambda \in
\mathbf{C}$ of finite modulus, then $P_0$ is hypoelliptic in $L^{2}$.
According to Fredholm's alternative, the condition is satisfied for all $P_0$.
Since $R$ is strictly weaker than $P_0$, also $P$ has $\sigma >0$.
 If we can let $P_0=\mbox{ Re P }$, for the conditions above, the operator $P$ is
hypoelliptic in $L^{2}$. The parametrix construction now gives,
for these conditions, a parametrix in $\mathcal{D'}$, with singular
support confined to $\Delta$, the diagonal in $\R^{\nu}
\times \R^{\nu}$.

\vsp

For a more formal construction, the following lemma is useful.
\newtheorem{mkt-reg}[He_par]{ Lemma }
\begin{mkt-reg} \label{mkt-reg}
Given an operator $P_0$, as above with $\sigma > 0$ and $E$ any parametrix, $PE - \delta_{x} \in C^{\infty}$ over
$\D$, then $E$ is very regular.
\end{mkt-reg}
Proof: We will show that for the commutator $C_E=E \psi_x - \psi_x E$, for a test function in $C^{\infty}_0(\R^{\nu}$
$\times \R^{\nu})$, $\psi_x =1$ close to $x$, we have $C_E \in C^{\infty}$. Assume $PE-\delta_{x} \in C^{\infty}$ also over $L^{2}$. Since $L^{2}$ with topology
induced by $C^{\infty}$, is a nuclear space, $C_E$ must be defined on $L^{2}$, by a kernel in $C^{\infty}$ and this
kernel will be regularizing also in $\D$. We see that, that $Eu - u \in C^{\infty}$ over $L^{2}$.
Thus, $E\psi_x u - \psi_x u- C_Eu \in C^{\infty}$ over $L^{2}$ and the result follows for $u \in \D$.$\Box$

\vsp

$\bf{Remark:}$ Note that it is a consequence of the results in section
\ref{sec:Lin} and proposition \ref{so}, that a constant coefficients,
self-adjoint operator $P$ with $\sigma > 0$, is hypoelliptic in $\D$, if
and only if we have an inequality
$$
\parallel C_P u \parallel_{H^{0,0}_K} \leq C \parallel u
\parallel_{H^{0,0}_K(P)} \ \text{ for } \  u \in H^{0,0}_K(P)
$$
where $K$ is such that $C_P \neq 0$ but otherwise arbitrary.

 Proof:(of the proposition) Assuming $R=\sum_jR_j$ has
$\sigma
>0$ and that the operator is extended with a constant coefficients
operator, outside a compact set, means that $R$ is hypoelliptic in $L^{2}$. The condition that
$R \prec P_0$, means that $\mbox{ Re }P$ has $\sigma > 0$ and that
$P$ is hypoelliptic in $L^{2}$. Finally, the parametrix construction
gives a parametrix in $\mathcal{D'}$, that is very regular, so
$P$ is hypoelliptic in $\mathcal{D'}$.$\Box$

\vsp

$\bf{Remark:}$ The importance of the requirement that $\mbox{ Re }L \sim
L$, is illustrated by the following example (cf.\cite{Rodino}). The
operator in $\R^{\textsl{2}}$ of order $2 m + 1$, $m \geq 1$,
$P=(D_x + i x^{2k} D_y)^{2m+1} - i x D_y^{2m}$, is for $k \geq 1 $
not hypoelliptic in $\mathcal{D'}$, but it is homogeneously
hypoelliptic in $\mathcal{D'}$ and as a consequence, hypoelliptic in
$L^{2}$.

\vsp

We note, that an operator that does not depend on all variables in space, can be
completed to a hypoelliptic operator, $P_{0,t}(D)=P_0(D')(1 - \Delta'')^t$. Let, for a partially hypoelliptic
operator, $L=P_0 + R$, $\widetilde{L}=P_{0,t}+R$, where $t$ is chosen such that $2 \mbox{ deg }_{x''} R < t$, that
is $R \prec \prec P_{0,t}$. For $\widetilde{L}$, with constant real coefficients, according to what was said
above, $\widetilde{L}$ is hypoelliptic in $\mathcal{D'}$. For $\widetilde{L}$, with variable coefficients, such that
$\mbox{ Re } \widetilde{L} \sim \widetilde{L}$, we also get a hypoelliptic operator in $\mathcal{D'}$. This means that given a
variable coefficients, self-adjoint, formally partially hypoelliptic operator, we can always complete this
operator to a formally hypoelliptic operator.

\vsp

\newtheorem{completing}[He_par]{ Lemma }
\begin{completing}
The constant, real coefficients operator $\widetilde{L}(D)=P_{0,t}(D')+R(D)$ is hypoelliptic in $\D$
\end{completing}
Proof:\\
It is sufficient ( and necessary ) to prove, for the commutator $C_{\widetilde{L}}$, defined as $C_{\widetilde{L}}=\widetilde{L} \phi
 -\phi \widetilde{L} \neq 0$, for a suitable test function $\phi
\in C^{\infty}_0(\R^{\nu})$, such that $\phi=1$ on a small open set in $K$ and $K$ arbitrary, that
$$\parallel C_{\widetilde{L}}u \parallel_{H^{0,0}_K} \leq C \parallel u \parallel_{H^{0,0}_K(\widetilde{L})} \ \text{ for } u \in {H^{0,0}_K}$$
But since this relation obviously holds for $P_{0,t}$ and since $R \prec \prec P_{0,t}$, implies an even stronger
relation, $\parallel R u \parallel_{H^{\sigma',\sigma''}_K} \leq \parallel u \parallel_{H^{0,0}_K(P_{0,t})}$, for
some positive numbers $\sigma',\sigma''$ (Proposition \ref{Mizo}), the result follows.$\Box$

\vsp

Note that for a partially hypoelliptic operator $L=P_0 + R$, we have $WF(u) \subset WF(Lu) \cup \Gamma_x$, $u \in
\mathcal{D'}$, where $\Gamma_x$ is the characteristic set for the operator $P_0$, that is the set of real zero's to
this polynomial. For complex zero's to a constant coefficients, hypoelliptic polynomial, we have with $\xi= \xi' + i \xi''$ and $P(\xi)=0$,
if $\xi''$ is bounded, then $\xi'$ will be bounded. Accordingly, for a constant coefficients, partially hypoelliptic operator,
$P(\xi,\eta)=0$ with $\xi'',\eta$ bounded, implies $\xi'$ bounded (cf.\cite{Gard}, Theorem 1). Let, for frozen $(x,y)$,
 $\Phi_x(L)=\{ (x,y,\xi,\eta); L(x,y,\xi,\eta)=0$
 $\quad 0 \neq \eta \text{ bounded }, \xi'' \text{ bounded } \}$. Then, for a partially hypoelliptic operator, the
 projection of $\Phi_x(L)$ on the real directions, $\xi',\eta'$, is in a bounded set in $\R^{\textit{n}} \times
 \R^{\textit{n}}$. For the completed operator, we get $\Phi_x(\widetilde{L}) { \subset \subsetneq \neq} \Phi_x(L)$, but for the
 corresponding real sets, $\Gamma_x(\widetilde{L})=\Gamma_x(L)$. Note that, since the completing polynomial has no
 real zero's, $\Gamma_x(\widetilde{L})$, is not necessarily bounded.

\vsp

If the operator is of the form of $L_{\lambda}^{\star}$, that is $L_{\lambda}^{\star}=L_{\lambda}^{\Sigma}+
 (P_{\lambda}^{\star} + R^{\star})$, where the two terms have support in complementary sets, then any iterate of the operator is of the
 same form. That is $(L_{\lambda}^{\star})^2 = (L_{\lambda}^{\Sigma})^2 + B_{\lambda}$, where
 $B_{\lambda}=P_{\lambda}^{\star}R^{\star}+ R^{\star}P_{\lambda}^{\star}$ with support outside $\Sigma$ and $(L_{\lambda}^{\Sigma})^N=
 (L_{\lambda}^{\Sigma})^N + D_{\lambda}$, and $D_{\lambda}$ with support outside $\Sigma$. Finally, given an operator $L$,
with variable coefficients and constant strength, (but not necessarily self-adjoint), we have that there is an iteration index
$N_0$, such that $L^N$ is hypoelliptic for every $N \geq N_0$. Because, $\mbox{ Re }L^N \sim (\mbox{ Re }L_{x^0})^N$, where
$L_{x^0}$ is the operator $L$ with frozen coefficients and where the right side is hypoelliptic,
for $N \geq N_0$, according to Lemma \ref{Ruc}. According to Lemma \ref{mkt-reg}, we thus have that $L^N$ is hypoelliptic
for every $N \geq N_0$.

\section{ Propagation of singularities in $\mathcal{D'}$ for the Neumann series}
\label{micro-local-I}
 The conditions on $L_{\lambda}^{\Sigma}$,
give that $E_{\lambda}$ must be regular in $x'$, that is $\varphi E_{\lambda} \in$
$H^{s,-N}_K$, for some $\varphi \in \mathcal{D}$ and $N$ finite. This would mean that the
parametrix has the same regularity properties, as a finite development
of the fundamental
solution constructed in section \ref{sec:constr_g_lambda}. The same method can be used to construct a fundamental solution, using
the parametrix we have now constructed. This time ${\alpha}_{\lambda}$ is replaced by
$A_{\lambda}E_{\lambda} \in
C^{\infty}$. $E_{\lambda}=K^{\Sigma}_{\lambda}+\sum_j E_j \otimes
\delta_{x''}$ with singular support on ${\Sigma}$. We assume
${\Sigma}$ adjusted to a singularity in $x$. In section \ref{sec:constr_g_lambda}, we
estimated $K^{\Sigma}_{\lambda}=O(1)e^{-\kappa \mid \lambda \mid^b}
  \otimes \delta_{x''}$, as ${\lambda} \rightarrow -\infty$,
on compact sets and $K^{\Sigma}_{\lambda}=O(1)\mid \lambda
\mid^{-c}$, as ${\lambda} \rightarrow -\infty$, uniformly on $\R^{\textit{n}} \times
\R^{\textit{n}}$. Using the estimates we have produced for the
solutions to the homogeneous equations, we have formally
$E_{\lambda}=O(1)e^{-\kappa \mid \lambda \mid^b} \otimes
\delta_{x''}$, on compact sets and $E_{\lambda}=O(1)\mid \lambda
\mid^{-c} \otimes  \delta_{x''}$, uniformly on $\R^{\textit{n}} \times
\R^{\textit{n}}$. We can write as before ${\alpha}_{\lambda}(x,z)=\sum_j
a_{\lambda,j,(z'')}(x',z') \otimes q_j(x'',z'')$, where
$a_{\lambda,j,(z'')}=0$ on ${\Sigma}$. Further
$v_{\lambda}={\alpha}_{\lambda}+ \big[
{\alpha}_{\lambda},{\alpha}_{\lambda} \big]+ \ldots=v_{\lambda,N}+ R_{v,\lambda,N}$
and $v_{\lambda} \in C^{\infty}_0$ has no support on ${\Sigma}$.

\vsp

Assuming as before that the variable $z''$ is frozen in $y'' \neq
x''$, then $v_{\lambda}$ will have support on ${\Sigma}$ and we will
get a representation of the fundamental solution, on a form much like
the one constructed in section \ref{sec:K-lambda-est} and with the same regularity properties.
Let's define $\D_{\Sigma}$-convergence as convergence in
$\D$, such that the singular support, if any, is maintained in
${\Sigma}$. Then, we obviously have $v_{\lambda,N} \rightarrow
v_{\lambda}$ in $\D_{\Sigma}$-meaning. Even for the fundamental
solution $g_{\lambda}=E_{\lambda}+ \big[ v_{\lambda},E_{\lambda}
\big]$, we get that $\big[ R_{v,\lambda,N},E_{\lambda} \big] \rightarrow
0$ in $\D_{\Sigma}$-meaning, as $N \rightarrow \infty$.

\vsp

\newtheorem{Dop_to_Cop}{ Lemma }[subsection]
\begin{Dop_to_Cop} \label{Dop_to_Cop}
  For a $u \in \E$ $(\R^{\nu})$, ${\alpha} \in C^{\infty}(\R^{\nu})$ and
  a constant coefficients operator $P(D)$, there is a ${\beta} \in
  C^{\infty}_0(\R^{\nu} \backslash \mbox{ sing supp } {\textit{u}} )$, such
  that in $\E(\R^{\nu}$ $\backslash \mbox{ sing supp } u )$
  $$ {\alpha}(P(D)u)=(P(D)\delta_0)*({\beta}u)$$
\end{Dop_to_Cop}

Using that the coefficients to the operator $L$ are in
$C^{\infty}$, we can write $L(z,D)u=(L(D) \delta_0)*( \beta u
)$, for $u \in \E$ and $\beta$ an operator corresponding to
multiplication with a function in $C^{\infty}(\R^{\nu} \backslash
\mbox{ sing supp }{\textit{u}} )$. That is, assume $E$ the convolution inverse in $\E$, corresponding to
$L(D)\delta_0$. If $L(z,D)=\sum_{j=1}^r \alpha_j(z) L_j(D)$, let $L^T(z,D)=\sum_{j=1}^r L_j(D) \beta_j(z)$, be an
associated operator. We then can find $\beta_j \in C^{\infty}$, such that in $\E$,
$ E*L(z,D)u=$ $E*L^T(z,D)u=\sum^r_{j=1} \beta_j u=\beta u $.
Let's denote the set of multiplication operators, corresponding to a development of
an operator $L(z,D)$, $M(L)$. Thus $\beta \in M(L)$, is the set $\{ \beta_j \}^r_{j=1}$, occurring in the representation of
$L^T$. Assume as before
${\alpha}_{\lambda}(x,z)=A_{\lambda}(z,D)E_{\lambda}(x,z)$, further
that $F_{\lambda}$ is the convolution inverse for
$L_{\lambda}(D)$, in
$L_{\lambda}^{\#}I_{E_{\lambda}}=\delta_{x}$. Assuming $E_{\lambda}$ a
two-sided fundamental solution to $L_{\lambda}(y,D_{y})$, this means
that $L^{\#}={}^tL_{\lambda}$.
According to Lemma \ref{Dop_to_Cop}, we can find a
$\gamma \in M({}^tL_{\lambda})$, with support outside $\Sigma$, such that $\gamma
I_{E_{\lambda}}(u)= F_{\lambda}*u,$ $u \in C^{\infty}_0 (\R^{\nu})$.
Using the Lemma \ref{Dop_to_Cop} one more time, gives $\gamma
I_{{\alpha}_{\lambda}}(u)=$ $F_{\lambda}*A_{\lambda}(D)
\delta_0*{\beta}u$, $u \in \E$, now assuming $\gamma$ with
support outside the singular support for $I_{E_{\lambda}}(u)$ and
$\beta \in M({}^tA_{\lambda})$ with support outside $\mbox{ sing supp } u$.
Note that the operator parts in the operators $A_{\lambda}$ and
${}^tL_{\lambda}$ are the same, so outside $\Sigma$, the operator
$I_{{\alpha}_{\lambda}}$ acting on $C^{\infty}_0$, corresponds to
multiplication with test functions or we could say that it has only
the localizing property.

\vsp

Assume $E_{\lambda,{\delta}}=K^{\Sigma}_{\lambda}+\sum_j E_j \otimes
T_{\delta}$, where $T_{\delta}$ is a very regular measure in
${\E}^{(\textsl{0})}$
with support contained in an open set $U_{\delta}$, such that $U_{\delta} \rightarrow
\{ x'' \}$, as $\delta \rightarrow 0$. We then have, for a fixed $x$,
$ \mbox{ sing supp }E_{\lambda,{\delta}}(x,\cdot) \subset $ $\Sigma$ $\cup
\{ x' \} \times U_{\delta}$. Since the coefficients in ${\alpha}_{\lambda,{\delta}}$ are
0 on $\Sigma$, we have $\mbox{ sing supp }{\alpha}_{\lambda,{\delta}}(x,\cdot) \subset$ $
\{ x' \}
\times U_{\delta}$. Let's write $Z_{\delta}=\{ x' \} \times U_{\delta}$. The support for
${\alpha}_{\lambda,{\delta}}(x,\cdot) \subset$ $W' \times U_{\delta}$, where $W'=\{ y'; y \in W \}$
and where $W$ is a compact set, in which the operator is
formally partially hypoelliptic.
Let's assume, that for the situation where $y''$ is fixed outside
$\Sigma$, ${\alpha}_{\lambda,(y''),{\delta}}$ is a measure and study how
iteration of $I_{{\alpha}_{\lambda},(y''),{\delta}}$, corresponds to
convolution. Using the Lemma \ref{Dop_to_Cop}, on the tensorized
integral operators, and writing formally
$A_{\lambda}(z',y'',D)=$ $P_{\lambda}(z',D') \otimes Q_{\lambda}(D'')$, we have
$\gamma_1
I_{{\alpha}_{\lambda,\delta}}'(v_1)=$ $F_{\lambda}*P_{\lambda}(D')\delta_0*\beta_1
v_1$, $v_1 \in C^{\infty}_0(\R^{\textit{n}})$
and $ \gamma_2
I_{{\alpha}_{\lambda,\delta}}''(v_2)=$ $G*Q(D'')\delta_0*\beta_2 v_2$, $ v_2
\in C^{\infty}_0$ $(\R^{\textit{m}})$
or equivalently $(F_{\lambda} \otimes
G)*A_{\lambda}(D)\delta_0*'\beta_1 v_1 *'' \beta_2 v_2$. We assume
that $\gamma_1 \in M({}^tP_{\lambda})$, $\gamma_2 \in M({}^tQ)$,
$\beta_1 \in M({}^tP_{\lambda})$ and $\beta_2 \in M({}^tQ)$.
We can again identify the operator parts, and we have $\gamma_1
\otimes \gamma_2 I_{{\alpha}_{\lambda,\delta}}(v_1 \otimes
  v_2)=\delta_0*'\beta_1 v_1*''\beta_2 v_2$. Iteration of this
  procedure, gives for test-functions $\varphi \otimes \psi \in C^{\infty}_0(\R^{\nu})$,
  $$ \gamma^{(N)} I_{{\alpha}_{\lambda,\delta}}^N(\varphi \otimes \psi)=
     \delta_0*\gamma^{(N-1)}I_{{\alpha}_{\lambda,\delta}}^{N-1}(\varphi
       \otimes \psi)= \ldots =\gamma^{(0)} \varphi \otimes \psi$$
for $\gamma^{(j)} \in M$ on tensor form, $j=0, \ldots ,N$
and the iteration can be repeated infinitely many times.

\section{ Propagation of singularities in measure topology }
\label{micro-local-II}
Assume $\mu \in {\E}^{(\textit{0})}(\R^{\nu})$, is such that
$\widehat{\mu}$ slowly
decreasing. Using the conditions on $\alpha_{\lambda,\delta}$ and the
Paley-Wiener theorem, we have that
$\widehat{I_{{\alpha}_{\lambda,{\delta}}}(\mu)} / \widehat{\mu}$ is an
entire analytic function. Using a theorem in \cite{Ho_Conv},(Theorem
3.6), we have a $F_{\Sigma} \in \E$, such that
$F_{\Sigma}*\mu=I_{\alpha_{\lambda,\delta}}(\mu)$. We can assume
  $F_{\Sigma}$ on tensor form. Further
  $I^N_{{\alpha}_{\lambda,{\delta}}}(\mu)=F_{\Sigma}*\ldots*F_{\Sigma}*\mu$, where the convolution
 is repeated $N$ times. The singular support for this representation, is contained in
 $\{ x_1 + x_2 + \ldots + x_N + y ; x_j \in \mbox{ sing supp }F_{\Sigma}, y \in \mbox{ sing supp }\mu \}$.
  For a fixed $x$ away from $0$,
  the singular support for $F_{\Sigma}$ is included in $W' \times$ $U_{\delta}$.
In ${\E}^{\textsl{(0)}}$, we have the equalities for the measures
under hand
$$F_{\Sigma}*(\mu_1 \otimes \mu_2)(x)=F_{\Sigma}*(\mu_1 \otimes
\delta_0)*(\delta_0 \otimes \mu_2)=F_{\Sigma}*'\mu_1*''\mu_2$$
although their singular supports may differ. We can use the tensor
form of $F_{\Sigma}$ and these equalities, to get
$I^N_{{\alpha}_{\lambda,{\delta}}}$ as an iteration of
partial convolutions. This
gives a particularly simple displacement of the singular support
for $\mu=\mu_1 \otimes \mu_2$. Note that $I_{L_{\lambda}^{\Sigma}(y,D)E_{\lambda}}(\mu)=\mu$, so
assuming ${}^tL_{\lambda}=L_{\lambda}$ and using that the differential
operator part  in $L_{\lambda}^{\Sigma}$ is $L_{\lambda}$, this means that
$I_{c_{\Sigma}E_{\lambda}}(\mu)=F_{\lambda}*\mu$, where $c_{\Sigma}$ denotes
multiplication with $C^{\infty}-$ functions, derived from the coefficients in
$L_{\lambda}^{\Sigma}$ as before. Further, if $cE_{\lambda}$ denotes
multiplication with the $C^{\infty}$-functions corresponding to the coefficients in $L_{\lambda}(y,D)$, we have
$I_{cE_{\lambda}}(\mu)=F_{\lambda}*\mu-F_{\lambda}*I_{K}(\mu)$, where
$K$ is regularizing
and we shall see that the last term does not effect the wave front
set, that is $WF(I_{cE_{\lambda}}(\mu))=WF(F_{\lambda}*\mu)$.

\newtheorem{Ho_L}[Dop_to_Cop]{ Lemma }
\begin{Ho_L} \label{asymp}
  Assume $X,Y$ open sets in $\R^{\nu}$, and $X'\times X''
  \subset X$, $Y' \times Y'' \subset Y$. Given a measure $\mu \in {\E}^{(\textsl{0})}$ $(X \times Y)$,
  (but not the Dirac-measure), with $\widehat{\mu}$ slowly decreasing, such that
  $\mu=\mu_1 \otimes \mu_2$, for $\mu_1,\mu_2$ very regular in $X'
  \times Y'$ and
  $X'' \times Y''$ respectively, iteration of partial convolution with
  the convolution kernel corresponding to $\mu_1$,
  followed by the kernel corresponding to $\mu_2$, will for $x$ fixed sufficiently far away from
  $0$, outside the diagonal after a finite number
  of steps, give a $C^{\infty}$-function.
\end{Ho_L}

$\bf{Remark:}$ The conditions on $\mu$ are sufficient to conclude that $\mu$ is a parametrix to a
partially hypoelliptic differential operator. Assume $\mu \in {\E}^{(\textsl{0})}$ invertible and otherwise according to the conditions. Then we have
existence of the formal inverse in ${\E}^{(\textsl{0})}$, why we can solve the equation $\big[\mu-1\big]*f
=w \in C^{\infty}$, according to $\sum\mu^j*w=\sum_0^N+\sum_{N+1}^{\infty}$, where the last term is in
$C^{\infty}$ and by localization we can assume $f \in \E$. Now choose the polynomial $P$ as a partially
hypoelliptic polynomial with $\Delta_{\mathbf{C}}(P) \subset Z_{\widehat{f}}$ (lineality). We then have
$$ \mbox{ sing supp }(\mu*f)=\mbox{ sing supp }(f)$$
$$ \mbox{ sing supp }(P(D)f)=\mbox{ sing supp }(f)$$
Chose $g \in \mathcal{D'}_{\textsl{L}^{\textsl{1}}}$ such that $g$ is hypoelliptic in $\mathcal{D'}_{\textsl{L}^{\infty}}$
and such that $P(D)g*f - \delta \in C^{\infty}$ ($\Rightarrow g*f-\delta \in C^{\infty}$), then
$$ P(D)\mu-\delta \sim P(D)g*\mu*f-\delta \sim$$ $$ \sim P(D)g*\big[\mu*f-f\big] + g*\big[ P(D)f-f \big]
+ \big[g*f-\delta \big] \in C^{\infty}$$
where $\sim$ indicates that the singular supports coincide. Note that if we assume $P^Ng*f-\delta \in
C^{\infty}$ we also have $g*\mu^N*\sum_0^N \in C^{\infty}$ and the proposition is that there exists a $g$, hypoelliptic
for convolution, such that
$\mbox{ch sing supp }(\mu^N*\sum_0^N ) \subset - \mbox{ch sing supp }g$

\vsp

Proof: Let's assume $\mu=E_{\lambda,\delta}$, an invertible
measure corresponding to $\alpha_{\lambda,\delta}$, according to
section \ref{micro-local-I}. Using the translation invariance for
constant coefficients operators, we first assume that $L_{\lambda}
\delta_0*F_{\lambda}=\delta_0$. By tensorizing the operator
$L_{\lambda}$, we have that $F_{\lambda}$ has support only in $\{
y_1 \geq 0, y_2 \geq 0, \ldots , y_\nu \geq 0 \}$. The displacement of the
support and singular support, during iterated convolution, will be
only in the direction of non-negative coordinates. Assume $Z$ an open
set containing the singularities for $F_{\lambda}$. The sets
$Z$,$Z+Z$,$Z+Z+Z$,$\ldots$, then constitute a countable covering
of the support for the iterated convolution. Since this support is
compact, we must have a finite sub covering.

\vsp

The singular support for $F_{\lambda}$ is included in $\cup \{y;  y_j=0$ for some $j \}$.
 Using that
$E_{\lambda,\delta}$ is invertible, this means that $ \mbox{ ch sing supp
  }F_{\lambda}^N*E_{\lambda,\delta}=$ $\mbox{ ch sing supp
  }F_{\lambda}^{N+1}*E_{\lambda,\delta}$, for some $N$,
where the left side is $\mbox{ ch sing
  supp }F_{\lambda}*\ldots*F_{\lambda}*E_{\lambda,\delta}$, with $N$ repeated
$F_{\lambda}'s$.

\vsp

As we translate the
singularity to some $x$ far away from $0$, we would prefer to write
$$L_{\lambda}(D) \delta_0*F_{\lambda}(x-y)
=\delta_0(x-y)$$
Replace $x$ with a point, close to the diagonal, but $\neq x$. This
corresponds to a displacement of the support for $F_{\lambda}$, away
from the coordinate axes. Denote the
corresponding kernel $E_{\lambda,\delta}^{\textit{x}}$. The singular
support for $E_{\lambda,\delta}^{x}(z,\cdot)$, will during the iterated
convolution with $F_{\lambda}$, be moved along the diagonal in the
direction of negative coordinates. We still have $\mbox{ ch sing supp
  }F_{\lambda,x}^N*E_{\lambda,\delta}^{x}= \mbox{ ch sing supp }E^{x}_{\lambda,\delta}$, which
after a finite number of iterations is a contradiction. That is,
outside the diagonal, we have $F_{\lambda,x}^N*E_{\lambda,\delta} \in C^{\infty}$.
$\Box$

\vsp

$\bf{Remark:}$
If $y \in ( \mbox{ neighborhood } x )$, for instance $\mid x-y \mid < \epsilon$, for some $\epsilon >0$. Assume $z$ a point on the
distance $\epsilon$ from the origin and with only positive coordinates. Then $y - 2 \epsilon z \notin (\mbox{ neighborhood }
x)$. This is the type of translation suggested in the proof above.

\section{ Asymptotic convergence for the Neumann series $v_{\lambda}$ }
Assume $N$, the smallest positive integer, such that the singular
support is stable, that is
$\mbox{ ch sing supp }v_{\lambda,N,\delta}=\mbox{ ch sing supp
  }v_{\lambda,N+1,\delta}$,
then Lemma \ref{asymp} applied to the operator $A_{\lambda,\delta}$ gives
  that, since the convolution operator corresponding to
  $A_{\lambda,\delta}$, has no support on $\Sigma$ and since this set contains
  the diagonal, we have $v_{\lambda,N} \in
  C^{\infty}(\R^{\nu} \times \R^{\nu})$.
Thus,
{\small $$L_{\lambda}(E_{\lambda,{\delta}}+\big[
v_{\lambda,N+1,{\delta}},E_{\lambda,{\delta}} \big])=\delta_{x}+ {\alpha}_{\lambda,{\delta}} +
v_{\lambda,N+1,{\delta}} + \big[ v_{\lambda,N+1,{\delta}},{\alpha}_{\lambda,{\delta}} \big]+
\text{ a term in } C^{\infty}$$}
Using that $v_{\lambda,N+2,{\delta}}= {\alpha}_{\lambda,{\delta}}+\big[
v_{\lambda,N+1,{\delta}},{\alpha}_{\lambda,{\delta}} \big]$, we see that we have in fact a
parametrix to the operator $L_{\lambda}$.

\vsp

Let's write $v_{\lambda,k}^{\otimes}$, for the
tensorized iteration. We then have
$$WF(L_{\lambda}(E_{\lambda,{\delta}}+\big[
v_{\lambda,k,{\delta}}^{\otimes},E_{\lambda,{\delta}}
\big])-\delta_{x}) \subset WF( L_{\lambda,(y'')}(E_{\lambda,{\delta}}+\big[
v_{\lambda,k,{\delta}}^{\otimes},E_{\lambda,{\delta}}
\big])-\delta_{x})$$
and this would work just as well as a parametrix.

\vsp

Let's write $v_{\lambda,k}^{\textit{p}}$, for the iterated partial
convolution described above. We then have
$$ WF( {\alpha}_{\lambda,{\delta}} ) \subset WF(
{\alpha}_{\lambda,{\delta}}^{\otimes} ) \subset WF(
{\alpha}_{\lambda,{\delta}}^{\textit{p}} ) $$
On $\Sigma$, we have $E_{\lambda,{\delta}}+ \big[
v_{\lambda,{\delta}},E_{\lambda,\delta} \big]$
$\rightarrow E_{\lambda}$, as ${\delta}
\rightarrow 0$ with convergence in $H'$.
Using the representation $g_{\lambda,{\delta}}=E_{\lambda,{\delta}}
+ \big[ v_{\lambda,N,{\delta}}^{\textit{p}}, E_{\lambda,{\delta}} \big]
+ R_{\lambda,N,{\delta}}^{\textit{p}}$, with
$R_{\lambda,N,{\delta}}^{\textit{p}} \in C^{\infty}$ in a sufficiently
small neighborhood of the diagonal, we have
convergence in $C^{\infty}$, for all but a finite number of terms, as
$\delta \rightarrow 0$.

\vsp

Finally, we have earlier established convergence in
$\E$, for the equivalent to $v_{\lambda}^{\otimes}$. In this case
we have that in the ''bad''
variable, any direction may be singular for $v_{\lambda,N}^{\otimes}$, while the
in ''good'' variable, we have a hypoelliptic situation. We say, for
$v_{\lambda,N}^{\otimes} \in \D_{\Gamma}$, that
$v_{\lambda,N}^{\otimes} \rightarrow v_{\lambda}^{\otimes}$ in $\D_{\Gamma}$-meaning,
if the convergence holds in $\D$, while the wave front set is
contained in $\Gamma$. Thus we must have
$$ \sup_V \mid \xi \mid^N \mid (\varphi v_{\lambda,N}^{\otimes}-\varphi
v_{\lambda}^{\otimes}) \widehat{} \mid \rightarrow 0$$
for $N=1,2, \ldots$ and $\varphi \in C^{\infty}_0($ neighborhood $W)$, such that
$\Gamma \cap (\mbox{ supp } \varphi \times V)= \emptyset$ for any closed cone
$V$. But since $V$ only can contain ''good'' directions, the
convergence follows immediately, from what has already been proved.
The same result can be given for the fundamental solution.

\section{ Asymptotic hypoellipticity for the operator $L_{\lambda}$ }
We can now give a definition of asymptotic hypoellipticity. Assume
$g_{\lambda}$ a fundamental solution to the partially formally hypoelliptic, formally self adjoint,
variable coefficients operator $L_{\lambda}$. Assume further that we
have $g_{\lambda} = E_{\lambda,\delta} + \Big[v_{\lambda,N,\delta},
E_{\lambda,\delta} \Big] + R_{\lambda,N,\delta} =
g_{\lambda,N}+r_{\lambda,N}$. According to Lemma \ref{asymp}, we have
for $N$ sufficiently large, $r_{\lambda,N} \in C^{\infty}$.
We first give the following lemma.

\newtheorem{Asymptot-HE}[Dop_to_Cop]{ Lemma }
\begin{Asymptot-HE}
If $L_{\lambda}$ is hypoelliptic, then $g_{\lambda,N}$ is very
regular, for every $N \geq 0$. Conversely, if $L_{\lambda}$ is
partially formally hypoelliptic and $g_{\lambda,N}$ is hypoelliptic
with $r_{\lambda,N} \in C^{\infty}$,
for every $N \geq 0$, then $L_{\lambda}$ is hypoelliptic.
\end{Asymptot-HE}

Proof: \\
Assume $L_{\lambda}$ hypoelliptic. We have
$L_{\lambda}(g_{\lambda}-r_{\lambda,N})=\delta_{x}-L_{\lambda}r_{\lambda,N}$,
for every $N \geq 0$. It is a trivial consequence of Lemma
\ref{asymp}, that  $L_{\lambda}r_{\lambda,N} \in C^{\infty}$, for every
$N \geq 0$. Thus $g_{\lambda} - r_{\lambda,N}$ is a $\mathcal{D'}-$
parametrix to the operator $L_{\lambda}$, which must be very regular.
Assume now $g_{\lambda}-r_{\lambda,N}$ hypoelliptic
 with $r_{\lambda,N} \in C^{\infty}$, then $g_{\lambda,N}$ is a
 parametrix to $L_{\lambda}$ and $$ \mbox{ sing supp
   }I_{(g_{\lambda}-r_{\lambda,N})}L_{\lambda}\mu=\mbox{ sing supp
   }\mu \quad \text{ for every } \mu \in \D$$ and since $g_{\lambda,N}$ is
 hypoelliptic, $ \mbox{ sing supp }L_{\lambda}\mu=\mbox{ sing supp
   }\mu$. $\Box$

\newtheorem{Asymp-HE}[Dop_to_Cop]{ Definition }
\begin{Asymp-HE}[ Asymptotically hypoelliptic operator ]
  Assume $g_{\lambda}=g_{\lambda,N}+r_{\lambda,N}$ a fundamental
  solution to a partially formally hypoelliptic operator as in the
  previous lemma. If, for a finite positive integer $N_0$, we have
  that $r_{\lambda,N} \in C^{\infty}$ and $g_{\lambda,N}$
  hypoelliptic, for all $N \geq N_0$, then the operator $L_{\lambda}$ is said to be
  asymptotically hypoelliptic.
\end{Asymp-HE}

We could also say that, for an asymptotically hypoelliptic operator,
 we have, for large $N$, $L_{\lambda}r_{\lambda,N} \in
 C^{\infty}$. Thus $L_{\lambda}(g_{\lambda}-g_{\lambda,N})
 \rightarrow 0$ as $N \rightarrow \infty$ in $C^{\infty}$ and using an inequality for
 constant strength operators, $g_{\lambda}-g_{\lambda,N} \rightarrow
 0$ in $C^{\infty}$, as $N \rightarrow \infty$. So, $g_{\lambda}$ is
 approximated by parametrices, asymptotically in $C^{\infty}$.

 \newtheorem{Parametrix}[Dop_to_Cop]{Proposition }
\begin{Parametrix}
Given a partially formally hypoelliptic operator
$L_{\lambda}(y,D_{y})=P_{\lambda}(y,D_{y'})+\sum_{j=1}^r P_j(y,D_{y'})Q_j(D_{y''})$, we
have a parametrix, $E_{\lambda}$, on the form
$K_{\lambda}^{\Sigma}+\sum_jE_j \otimes \delta_{x''}$ with
singular support on $\Sigma=\{ z \in \R^{\nu};$ $z''=x'' \}$. Here
$K_{\lambda}^{\Sigma}$ is the fundamental solution to the operator
$P_{\lambda}(y',x'',D_{y'})$ and $E_j$ solutions to the homogeneous
equations $P_j(y',x'',D_{y'})E_j=0$ for $j=1, \ldots ,r$. We have the
estimates $E_{\lambda}=O(1)e^{-\kappa \mid \lambda \mid^b} \otimes
\delta_{x''}$, as $\lambda \rightarrow - \infty$, on compact
sets in $\R^{\nu} \times \R^{\nu}$. Further $E_{\lambda}=O(1)\mid
\lambda \mid^{-c} \otimes  \delta_{x''}$ uniformly on $\R^{\textit{n}}
\times \R^{\textit{n}}$, as $\lambda \rightarrow -\infty$.
\end{Parametrix}

\section{ Conclusions concerning $\text{g}_{\lambda}$}
\label{sec:concl}
The fundamental solution to the constant coefficients operator, $K^+_{\lambda}$, is obviously of exponential
$\rho^*$-type 0, in the bad variable, over $\mathbf{C}^{\nu}$. For $K_{\lambda}^+=T_1 \otimes T_2$, we have
$$ \mid \widehat{T}_1(i \xi') \mid \leq C \sup \mid \big[ \frac{1}{M(- \xi')-\lambda} \big] \mid$$
Thus $K^+_{\lambda}$ is of exponential $\rho^*$-type 0, also in the "good" variable. According to \cite{Mart} Ch. 2,
Corollarium 2, this means that it allows real support. For $u_{\lambda}$, we established in section \ref{sec:Paley-Wiener}
that it is in $\E$, which means that it can be represented  as an analytic functional, by a measure with compact support in $E=\mathbf{C}^{\nu}$.
We will use the same notation $u_{\lambda}$ for these elements. Let's assume
that it is portable by a ball with respect to a complex norm $\rho$ and of radius $\alpha$. According to \cite{Mart} Ch. 2,
Lemma 1, this means that it is of exponential $\rho^*$-type $ \leq \alpha$. The problem is to establish whether
also $u_{\lambda}$ allows real support, in which case the same holds for the representation of $g_{\lambda}$ as analytic functional.
Immediately, if $E^*$ is instead a compact set in $\mathbf{C}^{\nu}$, we have that $u_{\lambda}$ is of exponential $\rho^*$-type 0 and allows
real support. Also, it is of exponential $\rho^*$-type $\alpha$ and allows real support in $H'(\R^{\nu})$.

\vsp

 In the general case, we can at least say, using \cite{Mart} Ch.2, Proposition 1.2, if $E_{\R}$ a real
vector space with complexification $E_{\mathbf{C}}$, $\mathcal{D}$ $(E_{\R})$ is a dense sub algebra of $H'(E_{\mathbf{C}})$
and
$$ u_{\lambda}(\varphi)= \int_{\R^{\nu}} \psi(x) \varphi (x) d x \qquad \varphi \in H(E_{\mathbf{C}}), \psi \in \mathcal{D}(\textit{E}_{\R})$$
The conclusion so far, is that $g_{\lambda}$ exists in $H'(E)$, $E=\mathbf{C}^{\nu}$ and is portable by a ball of radius $\alpha$, with respect to a complex
norm $\rho$. Further, as $B_{\rho,\alpha}$ can be regarded as an analytic variety, using \cite{Mart} (Ch.I Cor. to Th. 2.5),
there exists a unique $H(E)$-convex set $W$ in $B_{\rho,\alpha}$, containing the $H(E)-$convex hull to the support for $g_{\lambda}$.

\vsp

For the next proposition we need: $$
G_{x,\lambda}^{(\alpha',\alpha')}(x',y') = \frac{1}{(2 \pi)^n}
\int \frac {\xi^{2\alpha'} \exp(i(x'-y')\cdot \xi')}{{\mbox Re }
P^x(\xi')-\lambda}d \xi'$$ where we assume the operator $P^x$
formally self-adjoint in $L^{2}(\R^{\nu})$. We write
$g_{\lambda}^{(\alpha,\beta)}(x,y)$ for the derivative
$(iD_x)^{\alpha}(iD_y)^{\beta}g_{\lambda}(x,y)$. We note that the
coefficients corresponding to $L_{\lambda}$, can be extended to
entire analytic functions, using $ \mid \widehat{c}_{\alpha}(z)
\mid \leq C e^{c \mid \mbox{ Im }z \mid}$ and the condition on
self-adjointness means that the symbol $L_{\lambda}(z,\zeta)$ can
be treated as real analytical. Using \cite{Ni_72} Lemma 10 we have
\newtheorem{Prop6}{Proposition }[subsection]
\begin{Prop6} \label{Prop6}
For any positive integer $M$, provided $\varrho$ (as in section \ref{sec:Sob_emb}) $> n+M$,if $\varphi,\psi \in C^{\infty}_0(\R^{\textit{m}})$
with support in a neighborhood of the origin, then there is
for all sufficiently large negative values of $\lambda$, an analytic functional $g_{\lambda}$
on $\mathbf{C}^{\nu} \times \mathbf{C}^{\nu}$ and on $\R^{\nu} \times \mathbb{R}^{\nu}$, with the following properties:
\begin{enumerate}
\item for every  $x \in \R^{\nu}$ the analytic functional on $\R^{\nu}$, $g_{\lambda}(x,\cdot)$ is
a fundamental solution with singularity $x$ to the operator $L(y,D_y)-\lambda$
\item $g_{\lambda}$ is in $H'$ on the form $K^+_{\lambda}(x,y) + \big[ u_{\lambda},K^+_{\lambda} \big](x,y)$, where
$u_{\lambda}(x,y)$ can be represented in $\E(\R^{\nu} \times \R^{\nu})$, as an
  infinite sum of tensor products
  {\small $$ u_{\lambda}=\sum_{\mid I
    \mid = 1}^{\infty} \sum_I
  \Big[\Big[p_{\lambda,i_0,(y'')},p_{\lambda,i_{\mid I \mid + 1},(y'')}
   \Big]_{\mid I \mid + 1}',h_{\lambda}\Big]'\otimes
    Q_{i_0} \ldots Q_{i_{\mid I \mid + 1}}\delta_{x''} $$} \par
\item $g_{\lambda}*_x''\varphi*_y''\psi$ belongs to $C^M(\R^{\nu} \times \R^{\nu}$) and
for every multi-order $\alpha$ with $2\mid \alpha' \mid \leq M$, we
have with some positive constant $c$
{\small $$ g_{\lambda}^{(\alpha,\alpha)}*_x''\varphi*_y''\psi(x,x) =
(1 + O(1)\mid \lambda \mid^{-c})G_{x,\lambda}^{(\alpha',\alpha')}
\Big[(iD_{x''})^{\alpha''} \varphi \Big] \Big[ (iD_{y''})^{\alpha''}
\psi \Big](x,x)$$ \par}
   $\lambda \rightarrow - \infty$, for every x $\in \R^{\nu}$.
\end{enumerate}
\end{Prop6}
If we, in the argument following (\ref{g-lambda}) use the estimates
for $K^+_{\lambda}$, that are proven in \cite{Ni_72} Lemma 10
(compare with the third item in the following proposition), we get
in the first two, (still according to \cite{Ni_72} Lemma 10)
\newtheorem{Prop7}[Prop6]{Proposition }
\begin{Prop7} \label{Prop7}
With conditions as in Proposition \ref{Prop6},
\begin{enumerate}
\item $g_{\lambda}*_x''\varphi*_y''\psi(x,y)=O(1)\mid \lambda \mid^{-c}, \quad \lambda \rightarrow - \infty$
     uniformly on $\R^{\nu} \times \R^{\nu}$, for some positive constant $c$.
\item for $\mid \alpha' \mid \leq M$, then
      $D^{\alpha}_xg_{\lambda}*''_x\varphi*_y''\psi(x,\cdot) \in C^{\infty}(\{ y
      \in \R^{\nu};$
      $y-x \notin$ $0 \times \mbox{supp } \psi \})$, for every
      $x \in \R^{\nu}$. Further, for all multi-index $\alpha, \beta$ , $g_{\lambda}^{(\alpha,\beta)}*_x''\varphi*_y''\psi(x,y)=$
     $O(1)\exp(-\kappa \mid \lambda \mid^b)$, $\lambda
    \rightarrow - \infty$, uniformly on compact subsets in $\R^{\nu} \times
      \R^{\nu}$, where $\kappa$ is a positive constant that may depend on
      the compact subset, and where $b$ is the positive number
      corresponding to $M$
      as in (\ref{b})
\item For the fundamental solution corresponding to the operator $P(y,D_{y'})-\lambda$, we have  the estimates,
$K_{\lambda}(x,y)=O(1) \mid \lambda \mid^{-c} \otimes
  \delta_{x''}$ as $\lambda \rightarrow - \infty$, uniformly on
  $\R^{\nu} \times \R^{\nu}$, for some positive constant
  $c$. Further for $\mid \alpha' \mid \leq M$, $D^{\alpha'}_x \Big(
  h_{\lambda} + \Big[ B, h_{\lambda} \Big] \Big)(x', \cdot) \in
  C^{\infty}(\R^{\textit{n}} \backslash$ ${x'})$. For all multi-index
  $\alpha',\beta'$, $K_{\lambda}^{(\alpha',\beta')}(x,y)=O(1)\exp(-\kappa
  \mid \lambda \mid^b) \otimes \delta_{x''}$, as $\lambda \rightarrow - \infty$,
  uniformly on compact sets in $\R^{\nu} \times \R^{\nu}$
  with $\kappa$ and $b$ as in 10.2.2. Finally,
  $K^{(\alpha',\alpha')}_{\lambda}(x,x)=(1 + O(1)\mid \lambda
  \mid^{-c})G^{(\alpha',\alpha')}_{x,\lambda} \otimes \delta_{x''}(x,x)$ as
  $\lambda \rightarrow -\infty$, for every $x \in \R^{\nu}$
    \end{enumerate}
\end{Prop7}

\section{ Homogeneously hypoelliptic operators }
\label{sec:Hom_HE}

We will as usual use the same notation for the Schwartz kernel and its
corresponding integral operator.
This means for $L_{\lambda}F_{\lambda}=\delta_{x}-\gamma$ and $\gamma
\in C^{\infty}$, that $L_{\lambda}F_{\lambda} \in \Phi$, that is it is
Fredholm. Using standard arguments from the Fredholm theory, we can assume that both $L_{\lambda}$ and $F_{\lambda}$ are Fredholm
operators.

\vsp

We consider the standard projections $P: L^{2} \rightarrow R(L_{\lambda})$,
$Q:H^{s,t} \rightarrow N(L_{\lambda})$. $L_{\lambda} \in \Phi(H^{s,t},L^{2})$ gives a
decomposition $H^{s,t}=X_0 \bigoplus N(L_{\lambda})$ and $L^{2}=Y_0 \bigoplus
R(L_{\lambda})$, where $N(L_{\lambda})$ denotes the solutions to the homogeneous equation
and $R(L_{\lambda})$ the range of $L_{\lambda}$. We can construct an inverse, $E_{\lambda}$, to
$L_{\lambda}$, considered as an operator on $X_0$, which is extended
on $Y_0$ to an operator in $B(L^{2},H^{s,t})$.
Using a fundamental theorem in the Fredholm theory, there is a
$E_{\lambda} \in \Phi(L^{2},H^{s,t})$, such that
$E_{\lambda}L_{\lambda}=I$ on $X_0$ and $L_{\lambda}E_{\lambda}=I$ on
$R(L_{\lambda})$. Let $P^{\bot}=(I-P)$, then
$P^{\bot}(L^{2})=N(E_{\lambda})$. Further,
$L_{\lambda}E_{\lambda}(I-P^{\bot})=(I-P^{\bot})$ or
$L_{\lambda}E_{\lambda}=I-P^{\bot}$, and $P^{\bot}$ is a finite rank
operator (a compact operator). In the same way
$E_{\lambda}L_{\lambda}(I-Q)=(I-Q)$, or
$E_{\lambda}L_{\lambda}=I-Q$. We note that, for
$\sigma > 0$, $L_{\lambda}$ is homogeneously $L^{2}$-hypoelliptic, which means
that $N(L_{\lambda})=Q(H^{s,t})$ $\subset C^{\infty}$, that is $Q$ is
regularizing on $H^{s,t}$ and $E_{\lambda}$ is a left parametrix to
$L_{\lambda}$. If $L_{\lambda}$ is homogeneously $L^{2}$-hypoelliptic, then also its Hilbert-space adjoint is
homogeneously $L^{2}$-hypoelliptic, that is $P^{\bot}(L^{2})=N(L_{\lambda}^{\mbox{ adj }}) \subset C^{\infty}$,
and $P^{\bot}$ is regularizing on $L^{2}$. We conclude that $E_{\lambda}$ is a left and right
parametrix to the operator $L_{\lambda}$.

\vsp

Noting that $L_{\lambda}E_{\lambda}=I-P^{\bot}$ in $L^{2}$ with $P^{\bot}$ regularizing, (we are assuming the projections non-trivial)
we have that $\mbox{ sing supp}_{L^{2}}(L_{\lambda}E_{\lambda}\varphi)=\mbox{ sing supp}_{L^{2}}(\varphi)$. Further, $P=I-P^{\bot}$, means that
$P$ is hypoelliptic on $L^{2}$. The same observations, can be
made for $E_{\lambda}L_{\lambda}=I-Q$, so $Q^{\bot}$ is hypoelliptic. Finally,
$\mbox{ sing supp}_{L^{2}}(\varphi)=\mbox{ sing supp}_{L^{2}}(L_{\lambda}E_{\lambda} \varphi)$ $ \subset \mbox{ sing supp}_{L^{2}}( E_{\lambda} \varphi )$, so $E_{\lambda}$ is
hypoelliptic. Also,
$\mbox{ sing supp}_{L^{2}}(\varphi)=\mbox{ sing supp}_{L^{2}}(E_{\lambda}L_{\lambda} \varphi) $ $\subset \mbox{ sing supp}_{L^{2}}(L_{\lambda} \varphi)$ and we conclude,
\newtheorem{HE-L2}{ Proposition }[subsection]
\begin{HE-L2} \label{HE-L2}
  On $L^{2}$,any homogeneously hypoelliptic operator $L_{\lambda}$ is hypoelliptic and
  conversely.
\end{HE-L2}
The extension of $E_{\lambda}$ to $R(L)$ can be made in different ways. If
$L_{\lambda}$ is assumed homogeneously $L^{2}$-hypoelliptic, then $E_{\lambda}$
can be defined as regularizing on $R(L)$. Considered as an operator on
$L^{2}$, $E_{\lambda}$ is then hypoelliptic.

\vsp

Assume $E_{\lambda}$ a parametrix to a homogeneously hypoelliptic, constant coefficients differential operator $L_{\lambda}$,
that is $L_{\lambda}I_{E_{\lambda}}=I-I_{\gamma}$ on $L^{2}$. Assume $Y_0=N(L_{\lambda}^{-1})$, where
$L_{\lambda}^{-1}$ is the Fredholm-inverse operator and $R(I_{\gamma})=Y_0, N(I_{\gamma})=X_0$.
Further, $I_{E_{\lambda}}L_{\lambda}=I-I_{\eta}$ on $L^{2}$, such that
$N(I_{\eta})=R(L_{\lambda})$, $R(I_{\eta})=N(L_{\lambda})$, then $E_{\lambda}$ works as a
Fredholm-inverse operator to $L_{\lambda}$. Given $X_0$, $R(L_{\lambda})$, by adding a regularizing operator if necessary,
we can find $\gamma,\eta$ with these null-spaces. Since the operator $I_{\gamma}=I_{\gamma}Q$, with $Q$
regularizing, also $I_{\gamma}$ will be regularizing. An analogous argument, gives that also
$I_{\eta}$ is regularizing. Thus, any parametrix to a homogeneously hypoelliptic operator, can be
adjusted to a Fredholm-inverse operator.

\section{ Some remarks on the distribution parametrices }
\label{sec:finite-order}
For a constant strength operator, extended with constant coefficients outside
a compact set, we have seen that Levi's parametrix method, gives parametrices in
$\mathcal{D'}$ ${}_{L^{2}}$. Assume $P_{\lambda}(D)$ a hypoelliptic operator with constant coefficients and $E_{\lambda}$ a
parametrix to this operator. Thus
$$ \parallel I_{P_{\lambda}E_{\lambda}} (\varphi) \parallel_{L^{2}} \leq \parallel \varphi \parallel_{L^{2}} +
\parallel I_{\gamma} (\varphi) \parallel_{L^{2}}$$ and we can assume $\gamma \in C^{\infty}_0$. We have
that $P_{\lambda}E_{\lambda} \varphi \in L^{2}$ and as we shall see, a "converse to H\"older's inequality" gives
that $E_{\lambda} \in L^{2}$. Particularly, $E_{\lambda}:L^{2} \rightarrow L^{2}$. If an
operator has parametrix on the form of a tensor product with the Dirac measure, this will be in
$\mathcal{D'}$ ${}_{L^{2}}^{-m+1,n/2+1}$ and it maps $L^{2} \rightarrow L^{2}$.
More generally, any constant coefficients operator
parametrix in $\mathcal{D'}$ ${}_{L^{2}}^{l}$, for $k + \mid \alpha \mid \geq n/2+1+m$,
$\mid \alpha \mid \leq l$, and $k/2$ the order of the polynomial of growth for the
"multiplier", corresponds to a bounded  integral operator  $L^{2} \rightarrow L^{2}$. For a more general
parametrix in $\mathcal{D'}$ ${}_{L^{2}}^{s}$, corresponding to a variable coefficients operator
$L_{\lambda}$, we have that $I_{E_{\lambda}}(\varphi)=Q(D)F$,
where $Q(D)$ is a constant coefficients polynomial of order $s$ and $F \in L^{2}$. We
then have that, if $L_{\lambda}(x,D)Q(D)F \in L^{2}$ and assuming the order of $Q$ larger than $n/2$,
that is $s \geq n/2 + 1$, then
$$\parallel I_{E_{\lambda}}(\varphi) \parallel_{L^{2}} \leq C \parallel L_{\lambda}(x,D)F \parallel_{L^{2}} < \infty $$

\vsp

Assume $E(\overline{x},y)$ corresponds to the conjugate with respect to $x$ in the kernel, that is
if for instance $E=E_1 \otimes E_2$, we have $E(\overline{x},y)=\overline{E_1} \otimes E_2$.
Let $E(\overline{x},\overline{y}) - E(x,y)=$ $E(\overline{x},\overline{y})-
E(\overline{x},y) + E(\overline{x},y) - E(x,\overline{y}) + E(x,\overline{y}) - E(x,y)=$ $\Sigma_1 +
\Sigma_0 + \Sigma_2$. For a symmetric operator on $L^{2}$, we have $\Sigma_1=\Sigma_2=0$, meaning that the
operator is symmetric separately in the respective variables. For a symmetric operator acting on $\mathcal{D'}^F$,
we have $I_{\Sigma_0} \sim \mbox{ Im }I_E$. If the operator is symmetric on $\mathcal{D'}$ (that is ${}^tI_E=I_E^*$),
we have $\mbox{ Im }I_E=0$, but this is usually not the case for homogeneously hypoelliptic
operators. However, we always have that $I_E: \mathcal{D'} \rightarrow \mathcal{D'}^F+\text{i}\mathcal{D'}^F$.

\vsp

Assume $I_E^*=\overline{I}_E$ on $\mathcal{D'}^F$. Then, for $\varphi \in C^{\infty}_0$ real
$C_{\overline{I}_{E}}=\overline{I}_E \varphi - \varphi \overline{I}_E= \big[ \varphi I_E \big]^*-\big[
\varphi \overline{I}_E \big]$. This is regularizing on $\mathcal{D'}$. Thus, $ -2 i \varphi \mbox{ Im }I_E=
\varphi ( I_E^* - I_E ) -C_{\overline{I}_{E}}$. By extending the right side with $\pm \varphi I$,
we get a regularizing effect on $\E$ as before. The difference between a homogeneously
hypoelliptic and a hypoelliptic operator is that for the latter it is sufficient to consider the
real part of the operator $P$, which means that we can assume $I_E^*=I_E$ on $\mathcal{D'}$ and this
gives a regularizing effect for $\mbox{ Im }I_E$ on $\mathcal{D'}$.

\newtheorem{he-2}[HE-L2]{ Proposition }
\begin{he-2}
Any constant coefficients, homogeneously hypoelliptic differential operator on $\mathcal{D'}$ is a hypoelliptic
operator on $\mathcal{D'}^F$.
\end{he-2}

Proof: Assume $P$ homogeneously hypoelliptic on $\mathcal{D'}$ with parametrix $E$.
We have considered the operators
$$ I_E: C^{\infty} \cap \mathcal{D'}^F \rightarrow \textsl{C}^{\infty}$$
meaning $P(D)E-I \in C^{\infty}$ over $\mathcal{D'}^{F}$
$$ C_{I_E}: C^{\infty} \cap \mathcal{D'} \rightarrow \textsl{C}^{\infty} $$
If $P(D)u \in L^{2}_{loc}$ and $P(D)u \in C^{\infty}(\Omega)$ for all open sets $\Omega \subset
\R^{\textit{n}}$, then for any test function $\varphi$
$$ \mbox{ sing supp } \varphi u=\mbox{ sing supp } \Big( I_E \varphi P(D) u -C_{I_E} P(D) u \Big)$$
and we have $\varphi u \in C^{\infty}$. $\Box$

\section{ Partially hypoelliptic operators considered as Fredholm operators }
Also, if
$\sigma$ arbitrary, that is $L_{\lambda}$ not necessarily
homogeneously $L^{2}$-
hypoelliptic, we can extend the definition of $E_{\lambda}$ with a
regularizing term and we get a hypoelliptic action on, at least part of
$L^{2}$. We now assume the operator $L_{\lambda}$ self-adjoint and partially self-adjoint.

\vsp

If $E_{\lambda} \in \Phi$, is on the form $I-K$, there is a positive
integer $N_0$, such that $N(E_{\lambda}^N)=N(E_{\lambda}^{N_0})$ for all $N \geq N_0$. We can rewrite
$T^{-1}_{{\lambda}}=(I-\lambda E_{\lambda})L_{\lambda}=L_{2 \lambda}$ in $X_0 \bigoplus N(L_{\lambda})$ and
$T^{-1}_{\lambda}=L_{\lambda}(I-\lambda E_{\lambda})=L_{2 \lambda}$ on $R(L_{\lambda}) \bigoplus Y_0$. We get a new decomposition
$Y=N(E_{\lambda}^{N_0}) \bigoplus N(E_{\lambda}^{N_0})^{\bot}$ and a corresponding decomposition of the operator
$T$, $T= T_1 \bigoplus T_2$, with $T_1= \sum_1^{N_0-1} \lambda^j E_{\lambda}^j$ and $T_2$ $L^{2}$-hypoelliptic. That is
$E_{\lambda}:L^{2} \rightarrow H^{s,t}$, where $s,t$ can be chosen arbitrary large, so on
$N(E_{\lambda}^{N_0})^{\bot} \cap R(L_{\lambda})^{\bot}$, $E_{\lambda}$ can be defined as regularizing. We can rewrite the
first expression for $T^{-1}_2$, as
${\lambda}E_{\lambda}=I_{\lambda}-L_{2 \lambda}E_{\lambda}$, where
$I_{\lambda}$ denotes the identity operator in $R(L_{\lambda})$. But
$R(L_{2 \lambda}) \subset R(L_{\lambda})^{\bot}$ and we have that $E_{\lambda}$ is in fact hypoelliptic on
$N(E_{\lambda}^{N_0})^{\bot}$. Thus, $\mbox{ sing supp}_{L^{2}}(T_2 \varphi)=\mbox{ sing
supp}_{L^{2}}(\varphi)$,
for $\varphi \in L^{2}$ and for $u$ sufficiently regular, the integral
$I_{T_2}(u)$, can be estimated in supremum norm without having to
regularize the kernel.

\vsp

We have seen that $E_{2 \lambda}$ is $L^{2}$-hypoelliptic in the bands of
ranges surrounding and including $R(L_{2 \lambda})$, but since $N(E_{N
\lambda})=\{ 0 \} \Rightarrow L_{N \lambda}$ $L^{2}$-hypoelliptic, it
is has only $L^{2}$-action in the
outer bands. Since $N(E_{\lambda}^{N_0})$ is a finite-dimensional space, normed with $L^{2}$-norm, nuclearity gives that this
$L^{2}$-action can be defined by a kernel in $L^{2}$. For $f \in R(L_{2 \lambda})$, we have the estimates,
$\parallel E_{2 \lambda} f \parallel_{L^{2}} \leq C \mid {\lambda}
\mid^{-1} \parallel f \parallel_{L^{2}}$.

\vsp

We can define a "(s,t)"-regularizing operator $C_{s,t}: L^{2} \rightarrow
H^{s,t}$, as $*'\U_{-\it{s}}*''\U_{-\it{t}}$, for
positive real numbers $s,t$. Let $E_{\lambda,t}=C_{0,t}E_{\lambda}$,
where $E_{\lambda}$ is the parametrix constructed in $H^{0,-N}$,
mentioned above. Then, for $t$ sufficiently large ($>N$),
$E_{\lambda,t}$ is a parametrix in $L^{2}$,
corresponding to a
$L^{2}$-hypoelliptic operator. That is the condition that an operator $P_{\lambda}: H^{0,t} \rightarrow L^{2}$ has
no derivatives in the $x''$-variables,
means that the operator $P_{\lambda,-t}$ is hypoelliptic as an
operator $L^{2} \rightarrow L^{2}$.
 Thus $\mbox{ sing supp}_{L^{2}}(P_{\lambda,-t}E_{\lambda,t}
 \varphi)$
 $=\mbox{ sing supp}_{L^{2}}(E_{\lambda,t} \varphi)=$ $\mbox{ sing supp}_{L^{2}}(\varphi)$. Following the argument in
 the beginning of the section, this means that $E_{\lambda,t}$, can be
 extended with ($C^{\infty}$)-regularizing terms to $L^{2}$. Finally,
 $E_{\lambda}=(1-\Delta_{y''})^{t/2}E_{\lambda,t}$, which means that the
 regularizing terms for $E_{\lambda,t}$, will be regularizing for
 $E_{\lambda}$ as well. Note that this does not necessarily mean that $P_{\lambda}$ is hypoelliptic,
 since its parametrix is not hypoelliptic on $R(P_{\lambda})$. Thus, a
 parametrix on the form of a tensor-product with the
 Dirac-measure, can be used to extend the definition of
 $E_{\lambda}$, outside the range of the operator. On the range we
 can use the Fredholm-inverse,$P_{\lambda}^{-1}$ which gives a hypoelliptic action,
 for any differential operator $L_{\lambda}$.

\vsp

 On the bands, where $E_{\lambda}$ is
 regularizing, we have analytic dependence on ${\lambda}$, for
 $v_{\lambda}$ and this gives estimates like
 $$ \parallel E_{\lambda}^N \parallel_{c} \leq C_N \mid {\lambda} \mid^{-N} \qquad
 {\lambda} \  \rm{ finite }$$
 and particularly, $\parallel E_{\lambda} \parallel_{c} \leq C \mid \lambda \mid^{-N} \qquad \forall N$
 and for $\lambda$ large, on these bands.

 \vsp

 Note that the parametrix to $L_{\lambda}$ in $L^{2}$, can be written
 $E_{\lambda}=\bigoplus^{\infty}_{j=-\infty} E_{\lambda}P_{j
 \lambda}$, where $P_{j \lambda}$ is the projection on
  $R(L_{j \lambda})$. Since adding a compact operator to the Fredholm-inverse, does not change the index or the form
  $I-K$, we let $E_{\lambda}$ be
  defined as $L_{\lambda}^{-1} + X$ on $R(L_{\lambda})$ with $X \in
  C^{\infty}$, as $X' \in C^{\infty}$ on $R(L_{\lambda})^{\bot} \cap
  N(E_{\lambda}^M)^{\bot}$, for some suitable $M$ and as a function in $L^{2}$ otherwise.
   Thus, this operator $L_{\lambda}$ is hypoelliptic on
   $N(E_{\lambda}^M)^{\bot}$. That is, for $\varphi \in L^{2}$,
   $\mbox{ sing supp}_{L^{2}}(\varphi)=\mbox{ sing supp}_{L^{2}}(L_{\lambda}E_{\lambda} \varphi) \subset
   \mbox{ sing supp}_{L^{2}}(E_{\lambda} \varphi)$, further $\mbox{ sing supp}_{L^{2}}(L_{\lambda}E_{\lambda} \varphi)=
   \mbox{ sing supp}_{L^{2}}(E_{\lambda}L_{\lambda} \varphi)$
   and the result follows.

\section{ Hypoellipticity in the $\lambda$-infinity}
\label{sec:Hyp-infty}

Assume $E_{\lambda}$ the parametrix to an operator with constant coefficients $L_{\lambda}$.
The parametrix setting gives that, if $L_{\lambda}$ is self-adjoint, $E_{\lambda}$ has Fredholm-index 0,
and that $\lim_{n \rightarrow \infty} \mbox{ dim } N(E_{\lambda}^n)=$ $\lim_{n \rightarrow \infty} \mbox{ codim }
 N( E_{\lambda}^n ) < \infty$. If there is an index $N_0$, such that $E_{N \lambda}=I-K$, for $N \geq N_0$ and for $K$ regularizing, then
 $E_{ M \lambda}=I-K'$ with $K'$ regularizing, for all $\mid M \mid \leq N_0$. According
 the previous paragraph, $E_{N \lambda}$ will be on
 this form as $\mid N \mid \rightarrow \infty$.
 More precisely, assume for instance, $(\delta_x-E_{2 \lambda}) \in \mathcal{R}$, where $\mathcal{R}$ is the
 set of regularizing operators in $L^{2}$. Further that $(\delta_x-E_{2
 \lambda + 1}) \in \mathcal{R}$, which particularly means that $E_{2
 \lambda + 1}$ is pseudo local. The first condition means that
 $L_{2 \lambda + 1} E_{ 2 \lambda} \in \mathcal{R}$. Finally, $E_{2 \lambda + 1}L_{2 \lambda + 1}E_{ 2
 \lambda}(\delta_x-E_{\lambda}) \in \mathcal{R}$, which through the second condition
 means that $(\delta_x-E_{\lambda}) \in \mathcal{R}$. We could say, that
operators on the form $I-K$ in $L^{2}$, constitute a radical subsystem
among $L^{2}$-hypoelliptic operators.
Note that if the kernel $N(E_{\lambda})=Y_0 \neq \{ 0 \}$, the decomposition $L^{2}=R(L_{\lambda})
\oplus Y_0$ indicates that $L_{\lambda}$ still can not be considered as hypoelliptic on $L^{2}$.

\section{ Asymptotically $L^{2}$-hypoelliptic operators }
\label{sec:As-HE}

We now claim that $E_{\lambda}^2$ is a parametrix to the operator
$L_{\lambda}^2$. That is $E_{\lambda}^2=(E_{\lambda}P_{\lambda}+E_{\lambda}(P^{\bot}_{\lambda}))^2$
=$(E_{\lambda}P_{\lambda})^2+(E_{\lambda}P^{\bot}_{\lambda})^2$. If $E_{\lambda}$
is $L^{2}$-hypoelliptic, then the same must hold for $E_{\lambda}^2$. We
make the following definition,

\newtheorem{ahe}{ Definition }[subsection]
\begin{ahe}
Assume $E_{\lambda}$ a parametrix in $L^{2}$, to the operator $L_{\lambda}$
and that $E_{\lambda}$ is $L^{2}$-hypoelliptic on $N(E_{\lambda}^{N_0})^{\bot}$
with $N_0$ chosen as the smallest positive integer, such that the null space remains stable.
We then say that $L_{\lambda}$ is hypoelliptic on $L^{2}$, if $N_0 = 1$
and that it is asymptotically ($N_0$-) hypoelliptic, if $N_0 > 1$.
\end{ahe}

We claim that if $L_{\lambda}$ is asymptotically $N_0$-hypoelliptic, then $L_{\lambda}^{N_0}$
is hypoelliptic in $L^{2}$. Assume $E_{\lambda,N_0}$ the usual $L^{2}$-parametrix to $L_{\lambda}^{N_0}
$. We first have to prove, that
$N(E_{\lambda,N_0}^{N_1})=N(E_{\lambda}^{N_0}) \Rightarrow N_1 =
1$. We have that
$E_{\lambda,N_0}L_{\lambda}^{N_0}=L_{\lambda}^{N_0}E_{\lambda,N_0}=\delta_x-\gamma$,
for some $\gamma \in C^{\infty}$ and $\gamma =0$ on $N(E_{\lambda}^{N_0})$.
Thus $N(E_{\lambda}^{N_0}) = N(E_{\lambda,N_0})$, that
is $N_1 = 1$ and $L_{\lambda}^{N_0}$ is hypoelliptic on $L^{2}$.

\vsp

Note that given a parametrix $E_{\lambda,N_0}$ to a $L^{2}$-hypoelliptic operator $L_{\lambda}^{N_0}$, if $N(E_{\lambda,N_0}) \neq
\{ 0 \}$, we can always add to $E_{\lambda,N_0}$, a solution to the homogeneous equation , $H$, non-zero on
$N(E_{\lambda,N_0})$, so that $N(E_{\lambda,N_0}+ H)$ is trivial.
We can now assume $N(E_{\lambda,N_0})=\{ 0 \}$ and get $E_{\lambda,N_0}-E_{\lambda}^{N_0} \in C^{\infty}$. That is
the parametrix to $L_{\lambda}^{N_0}$, has a regularizing action on $N(E_{\lambda}^{N_0})$.

\vsp

The fact that $L_{\lambda}^{N_0}$, is $L^{2}$-hypoelliptic, renders a spectral kernel in $C^{\infty}$. The relation between
the spectral kernels corresponding the operator and the iterated operator, gives the spectral
kernel corresponding to $L_{\lambda}$ in $L^{2}$. This however, does not imply that
it is in $C^{\infty}$, on $N(E_{\lambda}^{N_0})$.

\vsp

Assume $U_{\lambda}=\{ \xi \in \R^{\nu}, \mid \xi'' \mid \leq \mid \lambda \mid \}$
, then $L^{2}(\R^{\nu})=$ $L^{2}(U_{\lambda}) \bigoplus L^{2}(\R^{\nu}$ $ \backslash U_{\lambda})$, where $L^{2} \ni f = f' +
f''$, and $\mbox{ supp } \hat{f'} \subset U_{\lambda}$. Assuming the frozen operator, $L_{\lambda}$, hypoelliptic
in $x'$, we can assume the corresponding parametrix, $E_{\lambda}=$ $E_{\lambda}' \otimes E_{\lambda}''$, adjusted to $N(E_{\lambda}')=\{ 0
\}$. The operator $E_{\lambda}''$, is well defined after adjusting the $L^{2}$-element, $\hat{f''}$ to a compact set.
Thus $N(E_{\lambda}) \subset L^{2}( \R^{\nu} \backslash U_{\lambda} )$. It is for the variable
coefficient case, sufficient to study operators on tensor form $L_{\lambda}^{\Sigma}$. This only involves an
operation on $E_{\lambda}'$, which does not effect $N(E_{\lambda}')$.

\section{ Some remarks on Weyl's criterion}
\label{sec:Weyl}
\begin{enumerate}
  \item[$\circ$] According to L. Schwartz, a condition equivalent with
    hypoellipticity for a differential operator $P$, is that $P$ and
    ${}^tP$ (the transposed operator) have parametrices, that
    are very regular. For the fundamental solution, we obviously have
    that Weyl's lemma implies that the fundamental solution (kernel)
    is very regular. For the
    opposite implication, there are counter examples in the variable
    coefficients case, for example the following differential operator
    (by Mizohata cf.\cite{Miz_op},\cite{Rodino})
    $$ P=\frac{\delta}{\delta x_1} + ix^h_1\frac{\delta}{\delta x_2}
    \qquad h \mbox{ integer}$$
    with fundamental solution $e(x,y)= \frac{1}{2 \pi} \big(
    x^{h+1}_1/(h+1) + ix_2 - y_1^{h+1}/(h+1) -iy_2 \big)^{-1}$.
    When $h$ is odd, $P$ is not hypoelliptic, since one solution to
    the equation $Pu=0$ is $u(x)= \big( x_1^{h+1}/(h+1) + ix_2 \big)^{-1}$.
    Finally, we wish to remark that the
    condition that the differential operator $P(D)$ is dependent on some
    variables in space, is essential for the opposite implication to
    hold. A trivial counterexample is the identity operator. This operator has
    the property of microlocal hypoellipticity, it is a hypoelliptic
    pseudo differential operator, but as we shall see,
    it is not a hypoelliptic differential operator.

    We have obviously, in the constant coefficients case, that if Weyl's criterion is to be written
    $ \mbox{ sing supp }P(D)u = \mbox{ sing supp }u$ for all $u \in \mathcal{D'}(\Omega)$
    and every $\Omega \subset \R^{\textit{n}}$, that we must consider the situation outside
    $\mathcal{D'}^F$ and the Sobolev-spaces. Assume $P$ with $\sigma > 0$ but not hypoelliptic
    ( that is not self-adjoint ) and consider
    \begin{equation} \label{Nordino}
    \parallel (P^*-P)u \parallel \leq C \parallel Pu \parallel \qquad u \in H^{0,0}_K
    \end{equation}
    for some constant $C$. This criterion can always be satisfied. Using $P^*-\overline{P} + \overline{P} - P$=
    $P^*-\overline{P}-2i \mbox{ Im }P$, if $P$ is localized with $\varphi \in C^{\infty}_0$, real
    and such that $C_{\overline{P},\varphi} \neq 0$ on a compact set, we have $C_{\overline{P},\varphi} - 2i \varphi
    \mbox{ Im }P=\big{[}\varphi P \big{]}^*-\big{[} \varphi P \big{]}$. Since, according to (\ref{Nordino}), $P^*-P \prec \prec P$, if
    $P$ has $\sigma > 0$ and if we assume $\overline{P}=P^*$ hypoelliptic, we must have $\mbox{ Im }P \prec
    \prec P$. Consider, for example $P+iP$, with $P$ hypoelliptic and real, then this is a
    homogeneously hypoelliptic operator, but not necessarily a hypoelliptic operator in $\mathcal{D'}$.
     With the additional condition that the operator is self-adjoint, we must however have that it is hypoelliptic.

    \vsp

     Consider now parametrices $E$ to
    partially hypoelliptic operators. We have seen that $E^N: \mathcal{D'} \rightarrow \mathcal{D'}^F$, for some $N$,
    but this does not mean that $E: \mathcal{D'} \rightarrow \mathcal{D'}^F$. For example, $T^2=(\gamma -i
    \delta_x)^2=\gamma^2-\delta_x-2i\gamma: \mathcal{D'} \rightarrow \mathcal{D'}^F$, but $T: \mathcal{D'}
    \rightarrow \mathcal{D'}$, particularly $\mbox{ Im }T: \mathcal{D'} \rightarrow \D$.
    Further, $\mid Tu \mid^2 \in \mathcal{D'}^F$,  $u \in \mathcal{D'}$ does not imply that $Tu \in
    \mathcal{D'}^F$, for instance $u$ real, if $w=(\gamma + i)u(\gamma - i)u$ and $(\gamma + 1)u(\gamma
    -1)u=-v$. Using that $v \in \mathcal{D'}^F$ and $w-v \in \mathcal{D'}^F$, we see that $w \in
    \mathcal{D'}^F$. We have already seen that geometric ideals have the property that the ideal is its complexification if
    and only if it is radical and since the class of parametrices to partially hypoelliptic
    operators are associated to geometric ideals, we get the same type of behavior for this
    class. Note that a necessary condition for hypoellipticity in $\mathcal{D'}$ is thus that $\mbox{
    Im }P \prec \prec P$ and this condition will be satisfied by iteration, that is $C_P^N \prec
    \prec P^N$ from some $N$, implies $\mbox{ Im }P^N \prec \prec P^N$. More precisely, assume $(,)$
     a complex scalar product over a Hilbert space $H$. Thus, $\mbox{ Re }(x,y)=\mbox{ Re}(ix,iy) \quad x,y \in H$
      and $(x,y)=\mbox{ Re }(x,y)-i\mbox{ Re }(ix,y)$. Then, for a constant coefficients differential operator
      $P=P_1+iP_2$, $$ \mbox{ Im }\big[ (P^2 \varphi,\varphi)-(P^*P \varphi,\varphi) \big]=2 i \mbox{ Re }(P_1 \varphi,P_2 \varphi)$$
      for $\varphi \in C^{\infty}_0$ arbitrary. If $P^*=P$, we must have
      \begin{equation} \label{pnrf} i\mbox{ Re }(P_1 \varphi,P_2 \varphi)=0 \end{equation}
      Conversely, if (\ref{pnrf}) then $(P-P^*) \bot P^*$ over $C^{\infty}_0$. Thus, if $E$ is the projection
      $R(P^*)^{\bot} \rightarrow N(P^*)$, then $P^*E(P-P^*)\varphi=0$ for $\varphi \in N(P^*)^{\bot}$. Particularly,
      if $\varphi$ is chosen as the mollifier (cf. the proof of Prop. \ref{Mizo}), we have $\parallel P^2 \varphi_n \parallel \rightarrow \mid P(\xi) \mid^2$ as $n
      \rightarrow \infty$.

      \newtheorem{ortho-strength}{ Lemma }[section]
      \begin{ortho-strength} \label{ortho-strength}
      Assume $P$ and $Q$ constant coefficients differential operators such that $P \prec Q$, $Q$
      hypoelliptic and $(P\varphi,Q\varphi)=0$ for all $\varphi \in N(Q)^{\bot}$. Then $P \prec \prec Q$.
      \end{ortho-strength}

      Proof: For a Fredholm inverse $E$ to $Q$ we have that $PE: L^{2} \rightarrow H^{\sigma}$ for a
      $\sigma \geq 0$. For an appropriate $\varphi$, H\"older's inequality gives that $\mid \xi
      \mid^{\sigma}P(\xi)\widehat{E} \rightarrow 0$ as $\mid \xi \mid \rightarrow \infty$, which is
      interpreted as $P \prec \prec Q$. $\Box$

      \vsp

      Assume now $P_2 \prec \prec P_1$, we can then find an entire $f$ with $P_2=fP_1$ and
      $$ (P_1 \varphi, P_2 \varphi)= \int f \mid P_1 \varphi \mid^2 d \xi $$
      with $\mid f \mid \rightarrow 0$ as $\mid \xi \mid \rightarrow \infty$. If we assume $f$ adjusted
      to $f'$ with support outside a ball containing the origin, we could say $(P_2 \varphi,P_1 \varphi)
      \sim 0$. For all $\varphi \in C^{\infty}_0$ we have $P_2 \varphi \in C^{\infty}_0$ and this means
      that there exists a $\gamma_{\varphi}$ regularizing, such that $P_2(I-\gamma_{\varphi})\varphi=0$.
      Further, $(P_2P_1E_1^{\varphi}\varphi,E_1^{\varphi}\varphi)=0$ and
      $P_1E_1^{\varphi}=I-\gamma_{\varphi}$. Assume $P_1$ hypoelliptic, then for all $\psi \in
      C^{\infty}_0$, there is a $\varphi \in C^{\infty}_0$, such that $E_1^{\psi}\varphi=\psi$
      ($E_1^{\psi}$ is chosen according to the support of $\psi$). If we choose $\psi$ so that $\parallel
      \psi \parallel=1$ in the mollifier, we have that $\parallel P^2 \varphi_n \parallel \rightarrow
      \mid P(\xi) \mid^2$ as $n \rightarrow \infty$. If $P_1$ is only partially hypoelliptic, we can
      consider $Q_N=P_1^N+iP_2^N$, such that $\mbox{ Im }\big[(Q_N^2 \varphi,\varphi)-(Q_N^*Q_N \varphi,\varphi)\big]=2 i
      \mbox{ Re }(P_1^N \varphi,P_2^N \varphi)$ and the above argument gives that $\parallel Q_N^2 \varphi_n
      \parallel \rightarrow \mid Q_N(\xi) \mid^2$ as $n \rightarrow \infty$.

      \vsp

     Assume $E_1 \in C^{\infty}(\R^{\nu}
   \backslash$ $0)$ a parametrix, that is with
   $\gamma \in \mathcal{D}$, such that $\gamma = 1$ in a neighborhood 0,
    $$ P(D)*(\gamma E_1* \varphi) = \zeta*\varphi + \varphi \qquad
    \zeta \in \mathcal{D}$$
   which means
   $$ \varphi= \gamma E_1*P(D)\varphi-\zeta*\varphi$$
   This is sufficient to conclude that $\mbox{sing supp }\varphi \subset
   \mbox{sing supp }P(D)\varphi$.
   If we use Leibniz' formula to construct a parametrix $\gamma E_1$,
   with $\gamma$ as before, that is,
   $$ P(D)*(\gamma E_1)=\gamma P(D)*E_1 + \zeta =\delta_0 + \zeta \qquad \zeta \in \mathcal{D}$$
   where $\zeta=0$ in a neighborhood 0, the same conclusion holds, unless
   $E_1$ has no support in the complement of the origin ($\zeta \equiv
   0$). In this case the corresponding operator is not necessarily a
   hypoelliptic differential operator.

   Note that the fundamental solution to a constant coefficients, hypoelliptic differential operator, is very
   regular in $\mathcal{D'}^{F}$, but it is not hypoelliptic in $\mathcal{D'}$ according to the argument
   above, since this would mean that the Dirac measure is hypoelliptic.
 \item[$\circ$]
The regularity behavior for PHE operators, can be expressed using
Sobolev spaces, in the following way:

\vsp

We have seen that a fundamental solution $g_{\lambda}$ with
singularity in 0, to the operator $P(D_{x'})-\lambda$, can be
given on the
form
$ g_{\lambda}=\mathcal{F'}^{-\textrm{1}} \Big(\frac{\textrm{1}}{\textrm{P}(\xi')-\lambda}
\Big) \otimes \delta_{\textrm{0}}$. We then have, for $\psi \in C^{\infty}_0$ so that $\psi=1$ in a
neighborhood $U$ in $\R^{\nu}$ of $0$,
$$ \psi u = g_{\lambda}* \Big( \psi L_{\lambda}u + ( L_{\lambda} \psi - \psi L_{\lambda})u -
QR(\psi u) \Big)=g_{\lambda}*\Big( \psi L_{\lambda}u + B_{\lambda}u -R'(\psi u) \Big)$$

We have, $B_{\lambda}=0$ in $U$ and $R'(\psi u)=\psi R' u +
R''u$ with $R''u=0$ in $U$. According to Proposition \ref{PHE}, we have
$\psi L_{\lambda} u= \psi P_{\lambda} u + \psi R' u$, so $\psi u= g_{\lambda}* \psi P_{\lambda} u$.
Is $P_{\lambda}$ hypoelliptic? Since
$C^{\infty}(U)=$ $\cap_{s,t}H_{loc}^{s,t}(U)$,
it would be sufficient to prove $\parallel \psi u \parallel_{s,t} \leq
C \parallel \psi P_{\lambda}(D_{x'}) u \parallel_{s,t}$, given arbitrary real numbers $s,t$.
Let $f_{\sigma}(D_{x'},D_{x''})=$ $\Big(P(D_{x'})-\lambda \Big) \Big( 1 -
\Delta_{x''} \Big)^{\sigma}$, for a non-negative number $\sigma$. We then
see that $F_{\sigma} = \F$ $(f_{\sigma} \delta_{\textrm{0}})$
is a weight function, defining a Banach space $H_{F_{\sigma}}$ through
the norm $\parallel u \parallel_{F_{\sigma}}^2$ \\
$=\int \mid \F$ $(u)(\xi) \mid^2 F_{\sigma}(\xi) d \xi$
and $F_{\sigma}^{-1}$ gives the antidual space to $H_{F_{\sigma}}$.
According to [\cite{Mal}, Lemma I.2.2]:

\newtheorem{Malg}{Proposition}[section]
\begin{Malg} \label{Malg}
Given two weight functions $h,f$,
if
$$\lim_{\xi \rightarrow \infty} \frac{h(\xi)}{f(\xi)}=0$$
then for a positive constant C, $\parallel \cdot \parallel_{h} \leq C \parallel \cdot
\parallel_{f}$.
\end{Malg}

\newtheorem{Mizo}[Malg]{Proposition}
\begin{Mizo} \label{Mizo}
Assume $M(D_{x'}),N(D_{x'})$  constant coefficients operators, where $M$
is assumed hypoelliptic over $\R^{\textit{n}}$, then $N$ is
  strictly weaker than $M$
if and only if there is a positive number $\sigma$, such that
\begin{equation} \parallel N(D_{x'})f \parallel_{s + \sigma,t} \leq C_K \parallel
M(D_{x'})f \parallel_{s,t} \qquad \text{ for all } f \text{ in }
\E(\textsl{K}) \label{sw}
\end{equation}
Here $s$ and $t$ are arbitrary real numbers.
\end{Mizo}
Proof:
The implication $N \prec \prec M$ $\Rightarrow$ (\ref{sw}), is proved in \cite{Miz_repr}.
For the opposite implication, it is sufficient to consider the case
$s=t=0$. We use the mollifier (cf.\cite{Schech}) $\varphi_k(x)=\mbox{e}^{i \xi
  \cdot x} \psi(x/k)/k^{n/2}$, $k=1,2,\ldots$, where $\psi \in C^{\infty}_0$,
such that $\parallel \psi \parallel_{0,0}=1$. Using that the
right side norm in (\ref{sw}) is equivalent to $\parallel \cdot
\parallel_{H^{0,0}_K(M)}$, we see that
$$\parallel N(D_{x'}) \varphi_k \parallel_{\sigma,0} \leq C_K \parallel
\varphi_k \parallel_{H^{0,0}_K(M)}$$ Taking the limit as $k
\rightarrow \infty$, we have that
$$\frac{(1 + \mid \xi' \mid^2)^{\sigma}\mid N(\xi') \mid}{1 + \mid
  M(\xi') \mid} \qquad \text{ bounded for } \xi' \in \R^{\textit{n}}$$
This is sufficient to conclude $N \prec \prec M$.
 $\Box$

 \vsp

Proposition \ref{Mizo} gives particularly, that a necessary condition for
hypoellipticity in $\D$, for an operator $P(D)$, is that $Id \prec \prec
P(D)$, this means in our case, $Id \prec \prec f_{\sigma}(D)$, for all $\sigma
\geq 0$. Proposition \ref{Malg} gives that
$ \parallel f^{-1}_{\sigma}*u
\parallel_{s,t} \leq C \parallel u \parallel_{s,t}$, for $u \in
H^{s,t}$. So,
$$\parallel \psi u \parallel_{s,t}=\parallel g_{\lambda}* \psi P_{\lambda}(D_{x'})u
\parallel_{s,t} \leq C \parallel \psi P_{\lambda}(D_{x'}) u \parallel_{s,t+\sigma}$$
but in the case of a partially hypoelliptic operator, we have to assume $\sigma > 0$ and if $P_{\lambda}(D_{x'})u \in C^{\infty}(U)$, we do not necessarily have that $u
\in C^{\infty}(U)$.
Assume $U$ a neighborhood of $x$, $\psi$ a test function with support
away from $0$, further that $u$
is a fundamental solution to $L_{\lambda}(D)$, we then have
$\parallel \psi u \parallel_{s,t} =$ $\parallel g_{\lambda}* \Big( \psi
L_{\lambda}(D)u + B_{\lambda}u + R'(\psi u) \Big)
\parallel_{s,t}$.
So on $U$, a fundamental solution to a
partially hypoelliptic operator cannot have better regularity properties
than $g_{\lambda}$. Note that for the identity operator,
the necessary condition for hypoellipticity in $\mathcal{D'}$, is not satisfied.
  \end{enumerate}

\section{Some results on the spectral kernel for partially hypoelliptic operators}

\subsection{Partial regularity for the spectral kernel}
\label{sec:Reg}
In this section, we follow the arguments of Nilsson \cite{Ni_61}. Assume $H=L^{2}(\R^{\nu})$ is our Hilbert space.
For the spectral family $\{ E({\lambda}) \}$, associated to
a realization, $\A_L$, corresponding to the formally self-adjoint operator $L(x,D_x)$ and for
$(\lambda_1,\lambda_2 \text{]} \subset \R$, we define the operators
$E(\lambda_1,\lambda_2)=E(\lambda_2)-E(\lambda_1)$,
as projections on the subspace of $H$,
$H(\lambda_1,\lambda_2)=H(\lambda_2) \ominus H(\lambda_1)$ (the
minus sign denotes the orthogonal
complement, of $H(\lambda_1)$ in $H(\lambda_2)$). For a given closed interval $I$
and a corresponding partition of finitely many subintervals $\{ I_j
\}^N_{j=1}$, each of length $\leq \epsilon$, we can
write $H_I= \bigoplus^N_{j=1}H_{I_j}$, where each of the subspaces is invariant
for $\A_L$. For $\lambda_j \in I_j$ and for every $x
\in H_{I_j}$, we have  $$ \parallel (\A_L - \lambda_j)x \parallel
= \parallel ((\lambda - \lambda_j)\chi_{I_j})(\A_L)x \parallel
\leq \sup_I \mid (\lambda - \lambda_j)\chi_{I_j}(\lambda) \mid \parallel x \parallel
\leq \epsilon \parallel x \parallel $$ that is $x$ is an $\epsilon-$approximative eigenvector to
$\A_L$ (cf.\cite{Ni_92}).

\vsp

We note that the spectral family or orthogonal spectral resolution,
uniquely determined by the operator $\A_L$, is a regular countably
additive spectral measure, non-decreasing and such that $E({\lambda}) \rightarrow 0$,
as ${\lambda} \rightarrow - \infty$ and $E({\lambda}) \rightarrow I$,
as $\lambda \rightarrow + \infty$. According to the spectral theorem,
we have
$\A_Lu = \int_{\R} \lambda d E(\lambda)u$  with strong convergence,
for $u$ in the domain of $\A_L$.

\subsection*{ The spectral resolution} \label{sec:spec-res}
We can show (analogously to [15] Theorem 3), using Proposition
\ref{Prop2}, for $kr > n/2 + l$, $\mid \alpha \mid \leq l$ and for a
compact set $K$,
there is an integer $N$ such that
\begin{equation} \sup_K \mid D^{\alpha}f*''\U_{-\textsl{2N}} \mid \leq \textsl{C}_{\textsl{K}} \parallel \textit{f}
\parallel_{\textsl{H}^{\textsl{s,-N}}_{\textsl{K}}(\textsl{P}^{\textit{r}})}
\label{o0} \end{equation}
Also, if $\tilde{f}=f*''\varphi$, $\varphi$ a test function with
support in neighborhood 0
\begin{equation}
   \sup_K \mid D^{\alpha}\tilde{f}(y) \mid \leq C(K)\big{(}
   \parallel \A_{P^r}f \parallel_H + \parallel f \parallel_H \big{)} \quad f \in H(\lambda_1,\lambda_2) \label{o1}
\end{equation}

\vsp

 The following result follows immediately from [8,II Theorem 10.4.8], where our additional
 condition, is due to the fact that we let $\alpha \rightarrow \infty$
 \newtheorem{Ho}{ Lemma }[subsection]
 \begin{Ho}
\label{Ho}
 Given a hypoelliptic, constant coefficients  operator $P$, adding $\alpha Q$, where $\alpha$
 is a complex constant, and $Q$ a strictly weaker, constant
 coefficients operator, such that $\mid P(\xi)+\alpha Q(\xi) \mid \neq
 0$ for $\mid \xi \mid$ sufficiently large, gives an operator
 equivalent with $P$
 \end{Ho}

\subsection*{ The constant coefficients case }
In section \ref{sec:parametrix} , we constructed a fundamental
solution with singularity in $x$, $h_{\lambda}$, to the constant
coefficients operator
$L_{\lambda}(D_{y})=L(D_{y})-\lambda=P_{\lambda}+R$, with singular
support on $\Sigma_x=\{ z; z''=x'' \}$. For a suitable test
function $\zeta \in C^{\infty}_0(W)$, $\zeta=1$ on $\Sigma_x \cup
F_2$, we can construct the parametrix $G_{\lambda}=\zeta
h_{\lambda}$ to $L_{\lambda}$, as $L_{\lambda}(\zeta h_{\lambda})=\delta_{x} -
\eta_{\lambda}$, where $\eta_{\lambda}$ is in $C^{\infty}$. This
gives a representation formula (similar to [15]), for $u$
sufficiently regular and for concentric spheres $F_4 \subset F_3
\subset F_2 \subset F_1
 \subset W$, where $W$ is a neighborhood of the singularity $x$.
 For a $\psi \in C^{\infty}_0$,
 such that $\psi = 1$ on $F_2$, we have for x $\in F_2$
\begin{eqnarray}
    u(x)= \int_{F_1\backslash F_2} \big( \overline{B_{\lambda}(x,y)} +
    \overline{\eta_{\lambda}(x,y)} \big) u(y)dy + \nonumber\\
    \int_{F_1} \overline{\psi(y)G_{\lambda}(x,y)}L_{\lambda}(D_{y})u(y)dy
\end{eqnarray}
where
$B_{\lambda}(x,y)=L_{\lambda}(D_{y})(1-\psi(y))G_{\lambda}(x,y)$.
According to the proof of Proposition \ref{Prop7}.3,
$\mid B_{\lambda}(x',y') \mid \leq C \exp(-\kappa \mid \lambda
\mid^b)$,
for $x \in F_3$, as $\lambda \rightarrow - \infty$ and the same
estimate holds for $B_{-\lambda}$,
as $\lambda \rightarrow \infty$,
on compact sets in $\R^{\textit{n}}$.
\newtheorem{ordo}[Ho]{ Notation }
\begin{ordo} \label{ordo}
 We write $B_{\lambda}(x,y)=O(1)\exp(-\kappa \mid \lambda
\mid^b) \otimes \delta_{x''}$, meaning $\mid B_{\lambda}*''_{x}\varphi*_{y}\psi \mid
\leq$ $ C \exp(-\kappa \mid \lambda \mid^b) \mid \varphi(x'') \psi(y'') \mid$, for
$\varphi \otimes \psi \in C^{\infty}_0(\R^{\textit{m}} \times \R^{\textit{m}})$.
\end{ordo}
The same estimate holds for $\eta_{\lambda}$, since it only involves
derivatives of $g_{\lambda}$ and of a test function. Further
$G_{\lambda}^{(0,\beta')}=O(1)\mid \lambda \mid^{-c} \otimes \delta_{x''}$, as $\lambda
\rightarrow - \infty$, uniformly on $\R^{\nu} \times \R^{\nu}$.
Let's localize $u$ with a $\phi \in C^{\infty}_0$ $(\R^{\nu})$,
$\phi = 1$ in $F_3$ and regularize with $\varphi \in
C^{\infty}_0(\R^{\textit{m}})$ with support in a neighborhood of 0. We then have,
for $\widetilde{u}_{\delta}=u*''\varphi_{\delta}$, $\parallel
\varphi_{\delta} \parallel_{L^1}=1$,
\begin{equation}
  \parallel \widetilde{u}_{\delta} \parallel_{H(F_3)} \leq C e^{-\kappa \mid
\lambda \mid^b}
\parallel u \parallel_{H(F_1 \backslash F_3)} + C'\mid \lambda
\mid^{-c} \parallel
L_{-\lambda}u \parallel_{H(F_3)} \qquad \text{ as } \lambda
\rightarrow \infty \label{Ni} \end{equation}
For finite ${\lambda}$, Lemma \ref{Ho} gives that $L-\lambda \sim_{x'}
L+\lambda$, so this estimate holds also for $L_{\lambda}$ with
$\lambda$ large and positive, if we can prove that the equivalence constant is independent
of $\lambda$. Note that we can assume $u$ with support contained in a bounded
domain, which means that the equivalence implies $\parallel
L_{-\lambda}u \parallel_{H} \leq C \parallel L_{\lambda}u
\parallel_{H}$,
for a constant independent of $u$ and $\lambda$. For the
independence, we need the following result, that can be found in
\cite{Boi}.
\newtheorem{Boj}[Ho]{ Lemma }
\begin{Boj} \label{Boj}
  For constant coefficients operators $Q,M$, such that $M$
  hypoelliptic and $Q \prec \prec M$, there are positive constants $C$
  and $k$, such that
  $$ \mid Q(\xi') \mid \leq C \tau^{-k}(1 + \mid \xi' \mid )^{-k}(\tau +
  \mid M(\xi') \mid )$$
  for every $\xi'$ in $\R^{\textit{n}}$ and $\tau \geq 1$ real.
\end{Boj}

First, let's denote
$ L^{-2}=\mid L \mid^2 + \mid \lambda \mid^2$
where $L$ is regarded as a polynomial over $\R^{\textrm{n}}$.
Lemma \ref{Boj} gives that $1/L^- \leq C/\mid \lambda
\mid^{\sigma}$, for some positive number $\sigma$ and with the
constant $C$ independent of $\lambda$. If $ L^- \in \mathcal{H}_{\sigma}$
and $\sigma \geq 1$, $\mid \lambda \mid/\mid L \mid \leq
C$, for every large $\lambda$. Otherwise, this result holds for the iterated
operator $L^r$ and we have $L^{-,(r)} \sim_{x'} L^{r}$, for large $\lambda$. The same
result gives also that $(L+{\lambda})^r \sim_{x'} L^{-,(r)}$, for any
large $\lambda$. We have the following result,
\newtheorem{param}[Ho]{ Proposition }
\begin{param} \label{param}
  Given a hypoelliptic constant coefficients operator
  $P$, assume
  $P \in \mathcal{H}_{\sigma}$, for $\sigma \geq 1$. Then, for $\lambda$ complex
  $$ P \pm \lambda \in \mathcal{H}_{\sigma} \qquad \mid \lambda \mid
\rightarrow \infty$$
\end{param}
We have a condition that the type operator $M \in \mathcal{H}_{\sigma}$,
for $\sigma > n$, which means that the operator need not be iterated
any further, in order to be partially hypoelliptic independently of $\lambda$.

\vsp

For $f \in H(\lambda - \epsilon, \lambda) $,  $\epsilon
> 0$, we have that f is an ${\epsilon}$-approximative eigenvector to L and
${\lambda}$, so
$ \parallel L_{\lambda}(D)f \parallel_H \leq \epsilon \parallel f
\parallel_H $  and using (\ref{Ni})
$$ \parallel \widetilde{f}_{\delta} \parallel_H \leq C \big( \exp(-\kappa \mid \lambda \mid^b)
 + \mid \lambda \mid^{-c}\epsilon \big) \parallel f \parallel_H$$
where ${\lambda}$ is assumed large and positive. The
constant $C$ is not dependent on $u$, $\lambda$ or on the test
functions used in the regularization.

\vsp

For $f \in H(\lambda - k\epsilon,\lambda)$, k a positive integer, we can write $f=\sum^k_1f_j$, where
$f_j \in H(\lambda - j\epsilon,\lambda - (j-1)\epsilon)$, so by selecting $\epsilon$ and $k$ in a suitable way, we get
$$ \parallel \widetilde{f}_{\delta} \parallel_H \leq C \exp(-\kappa \mid \lambda \mid^b)
\parallel f \parallel_H $$
where the constant $C$ is independent of $f$,$\lambda$ and the
regularization and $\lambda$
is assumed large and positive.

\vsp

For $f\in H(\lambda - 1,\lambda)$, we also have $\A_{P^r}f \in H(\lambda - 1,\lambda)$
so $$ \parallel \widetilde{\A_{P^r}f}_{\delta} \parallel_H \leq C_r \exp(-\kappa \mid
\lambda \mid^{b}) \parallel f \parallel_H $$
Thus for $\mid \alpha \mid \leq l$ and kr $>$ n/2 + l
$$ \sup_{\overline{F_4}} \mid D^{\alpha} f*''\varphi(y) \mid \leq C(\alpha)
\exp(- \kappa \mid \lambda \mid^{b})
\parallel f \parallel_H $$ where $F_4$ is a sphere
with center $x_0$, $\varphi$ a test function with support in a
sufficiently small nbhd of 0 and
$C$ is independent of $f$,$\lambda$ and the regularization, but not of $F_4$. Here $\lambda$
is assumed large and positive.  Finally, this result follows for a general
compact set $K$, by the Heine-Borel theorem.

\subsection*{ The variable coefficients case }

Let's now consider the variable coefficients case. In section
\ref{sec:parametrix}, we constructed a parametrix also to
the variable coefficients operator $L_{\lambda}(y,D_{y})$ $F_{\lambda}=
\delta_x-\gamma_{\lambda}$, with ${\gamma}_{\lambda} \in C^{\infty}$,  where we have
adjusted the singularity to $x$, such that
$P_{\lambda}^{\Sigma}K^{\Sigma}_{\lambda}=\delta_x$.
We get a representation formula similar to (\ref{Ni}),
where we assume $W$ a nbhd of $x_0$, in which the operator is partially formally hypoelliptic,
$L_{\lambda}F_{\lambda}=\delta_{x}+\gamma_{\lambda}$ and
$B_{\lambda}(x,y) = L_{\lambda}(y,D_{y})(1 -
\psi(y))F_{\lambda}(x,y)$ with support on $F_1 \backslash$ $F_2$.
To simplify the calculations we use a representation
$$u*''\varphi_{\delta}(x)=I_{\overline{\psi
  B_{\lambda}^{\delta}}}(u)(x)+I_{\overline{\psi
  F_{\lambda}^{\delta}}}(L_{\lambda}(\psi u))(x)- I_{\overline{\psi \gamma_{\lambda}^{\delta}}}(u)(x)$$
where $B_{\lambda}^{\delta}=B_{\lambda}*''\varphi_{\delta}$ and
analogously for $\gamma_{\lambda}^{\delta}$.
In order to produce an estimate like (\ref{Ni}), we need a fine
estimate of $\gamma_{\lambda}$.
Using \cite{Miz_repr} (Cor. 2 to Prop 2.1) and Proposition \ref{Mizo}, we can give the following result

\newtheorem{rest}[Ho]{ Lemma }
\begin{rest} \label{rest}
Given a variable coefficients operator $P$, with
coefficients in $C^{\infty}(\R^{\nu})$ and $=0$ on $\Sigma_{x}=\{ (y',y'');
y''=x'' \}$, we have
$$ \parallel P(y,D_{y})T \parallel_{s+\sigma,-N} \leq \epsilon \parallel M(D_{y'})T
\parallel_{s,-N'} \qquad M(D_{x'})T \in H^{s,-N'}_K$$
where $s$ is a real number, $N,N'$ positive integers and $\sigma$ a real
number that can be chosen as positive if $P \prec \prec_{x'} M$ and as
zero if $P \sim_{x'} M$. Finally $\epsilon$ can be chosen arbitrarily
small as the support for $T \rightarrow \Sigma_{x}$
\end{rest}

\vsp

$\bf{Remark:}$ The set $\Sigma_{x}$ can be chosen in different ways, but assuming
$E_{\lambda}= K^{\Sigma}_{\lambda} + \sum_j E_{\lambda,j} \otimes \delta_{x''}$, where
$E_{\lambda,j}$ is an arbitrary solution to the homogeneous equation,
$\Sigma_{x}$ according to the Lemma, seems to be the natural choice.

\vsp

In section \ref{sec:parametrix}, we saw that the remainder corresponding to the parametrix,
is on the form of operators $C_{\lambda},D_{\lambda}$ with
coefficients vanishing on $\Sigma_{x}$, acting on
$E_{\lambda}$, that is
$ \gamma_{\lambda}=C_{\lambda}E_{\lambda}+D_{\lambda} E_{\lambda}$,
where $C_{\lambda} \sim_{x'} M$ and $D_{\lambda} \prec
\prec_{x'} M$.
If we mollify
$E_{\lambda}$ appropriately, we have
$E_{\lambda}^{\delta}=E_{\lambda}*''\varphi_{\delta}$ with $\mbox{
  supp }E_{\lambda}^{\delta} \rightarrow \Sigma_{x}$,  as $\delta
\rightarrow 0$.
We then have the following estimate of the remainder term, for
${\lambda}$ sufficiently large and for some positive constant $c$,
\begin{equation}
  \parallel \int \overline{\phi \psi \gamma_{\lambda}^{\delta}(x,y)}u(y)d y \parallel_H\leq {\epsilon_1}
\mid \lambda \mid^{-c}
\parallel \phi \parallel_H \parallel \psi \parallel_H \parallel u
\parallel_H
\label{rest-eps}
\end{equation}
where we have used Lemma \ref{rest}, Leibniz' formula, Cauchy-Schwarz' inequality and the estimate
in Prop. \ref{Prop7}. Here, ${\epsilon}_1$ is dependent on the value of
the coefficients corresponding to $L_{\lambda}$ in a nbhd of
$\Sigma_{x}$ and ${\epsilon}_1 \rightarrow 0$ as $\delta \rightarrow
0$. Further, $\epsilon_1$ is dependent on the mollifier,
but it is not dependent on the support of $u$.

\vsp

An argument similar to the constant coefficient case, gives for  $\mid \alpha \mid \leq l$ and kr $> n/2 + l$,
$$ \sup_{K} \mid D^{\alpha} f*''\varphi_{\delta}(y) \mid \leq C(\alpha,\delta)
\exp(- \kappa \mid \lambda \mid^{b})
\parallel f \parallel_H $$
and still, the constant $C$, is dependent on the mollifier and on the
compact set $K$, but not on $f$ or $\lambda$.

\vsp

{ $\bf{Remark}$ (1):
The right hand side in (\ref{o1}) can be used as a definition of a
norm. We prefer
in this case to work with the Hilbert space $H^{0,-N}_K$, $K$ a compact set in $\Omega$ and we let
$\mid f \mid_{r,N} =$  $ ( \parallel f \parallel_{H^{0,-N}_K} + \parallel P^rf \parallel_{H^{0,-N}_K} )$.
Let $H^{0,-N}_K(P^r)$ denote the Hilbert space of elements in $D(P^r)$, r $\geq$ 0, such that
$\mid \cdot \mid_{r,N} < \infty$. The argument above applied to the spaces $H^{0,-N}_K(P^r)$ gives, for $\U_{-N}$ as in section \ref{sec:Sob_emb}}
\begin{equation}
  \sup_{K} \mid D^{\alpha} f*''{\U}_{-N}(y) \mid \leq C
 \mid f \mid_{r,N} \leq C_{\delta}' \exp(- \kappa \mid \lambda \mid^{b})
 \otimes \delta_{x''} \mid f
\mid_{0,N}  \label{Odh} \end{equation}
The constant in the last expression, is dependent on the choice of mollifier.
{ Note that in section \ref{sec:Sob_sp}, we proved that $\parallel \cdot
\parallel_{H^{s,-N'}_K(L)}$ is norm equivalent to $\parallel
\cdot \parallel_{H^{p_s,-N}_K}$ and we can show that this implies
$\parallel  \cdot \parallel_{H^{s,-N'}_K(\A_L)}$ is
norm equivalent to $\parallel  \cdot
\parallel_{H^{s,-N}_K(\A_P)}$, so the inequalities we can prove for
$\A_P$ also holds for $\A_L$, after adjusting the order of
the Sobolev space in the ''bad'' variable. Since according to section
\ref{sec:Sob_emb},
the iteration of the operator is done to satisfy a condition on the
''good'' variable, we prefer to work with realizations of the operator $P$.}

\vsp

For $f\in H$ and for the resolution corresponding to a realization of
the operator L, we can write $E(\lambda)f=f_1+f_2+\ldots$, for $\lambda
\neq - \infty$. For the spectral family corresponding to the iterated
operator, $E_r(\lambda)$ and r odd,  we have $E_r(\lambda^r)=E(\lambda)$.
Through (\ref{Odh}) and the partial norm equivalence proved
in section \ref{sec:Sob_sp},  $\parallel f \parallel_{H^{0,-N'}_K(L^r)} \leq C \mid f
\mid_{r,N}$, we know that
$$ \sup_{K} \mid D^{\alpha}E(\lambda)f \mid \leq
C(K,\alpha,\delta) \exp(-\kappa \mid \lambda \mid^{b}) \otimes \delta_{x''} \parallel f
\parallel_H $$  for $f\in H$.

\subsection*{ The spectral kernel}
Assuming the operator  $\A_{L^r}$,
$\U_{\textit{0},-\textsl{N}}$-partially formally
self-adjoint, we can construct the spectral
kernel, in the Hilbert space $H^{0,-N}_K$. \\ (Since we also assume that
the operator $L^r(y,D_y)=$ $\sum^t_{j=1}P_{j,(r)}(y,D_{y'})Q_{j,(r)}(D_{y''})$, is
formally self-adjoint in $L^{2}$, this means a requirement that the
operators $P_{j,(r)}(y,D_{y'})$ commute with the
weight operators). The first estimate in
(\ref{Odh}), gives a bound for the ''resolution '', $E_{N}$ on $H^{0,-N}_K$,
\begin{equation} \sup_{\overline{F_4}} \mid D^{\alpha}E_{N}(\lambda)g \mid \leq C \mid g
\mid_{0,N} \label{agmon} \end{equation}
Define $T^{(\alpha)}_{\lambda,N}(x)g=D^{\alpha}E_{N}(\lambda)g(x)$, as in
\cite{Ni_61} (the derivatives are here taken in distribution sense).
Schwartz kernel theorem (and for $T^{(\alpha)}_{\lambda,N}(x)$
interpreted as an evaluation functional, Riesz's representation theorem),
gives existence of a $f_{\lambda}^{(\alpha)}(x,\cdot) \in H^{0,-N}_K$, such that
$T^{(\alpha)}_{\lambda,N}(x)g =$ $(g,f^{(\alpha)}_{\lambda}(x,\cdot))_{0,-N}$, for $x \in K$
and $g \in H^{0,-N}_K$ . We then  have an implicitly defined spectral
function in $L^{2}(\R^{\nu})$,
$e_{\lambda}(x,y)=f_{\lambda}*''_y{\U}_{-N}$
$(x,y)$,
on a compact set. Note that $e_{\lambda}$ does not necessarily have
compact support. The situation outside the compact set is dealt with
in the end of this section.  On the other hand, if $e_{\lambda}$ is the spectral
function constructed on a compact set in $L^{2}$,
our spectral kernel is given by $ f_{\lambda}^{(\alpha)}(x,y) = (1-
\Delta_{y''})^{N/2}e_{\lambda}^{(\alpha)}(x,y)$.

\vsp

For the domain of definition to $T_{\lambda,N}$, we note that the kernel in
$H^{0,-N}_K$, is defined as
\begin{displaymath}
  f_{\lambda}(x,\cdot) = \left\{ \begin{array}{ll}
    f_{\lambda}^L & x \in K \\
    f_{\lambda}^M & x \notin K \end{array} \right.
\end{displaymath}
where $f_{\lambda}^L,f_{\lambda}^M$, are the kernels corresponding
to the operators $L$ and $M$ respectively. Further, in the case where
$x \in K$, we have
\begin{displaymath}
  f_{\lambda}(x,y) = \left\{ \begin{array}{ll}
    f_{\lambda}^L(x,y) & y \in K \\
    f_{\lambda}^M(x,y) & y \notin K \end{array} \right.
\end{displaymath}
If in the scalar product $(f,g)_{0,-N}$, the weight is brought to one
side, which would give an equivalent definition of the space
$H^{0,-N}$, that is for $g \in L^{2}_K$, $(g,f)_{0,-N}^{*}=(g,(1-\Delta_{y''})^{-N}f)_{0,0}$ $=(g,f)_{0,-N}$, then
$T_{\lambda,N}$ can be defined on $L^{2}_K$.
For an element $g \in L^{2}(\R^{\nu})$, we
have that $g=g_1+g_2$, with $g_1 \in L^{2}_K$. So
$$(g,f_{\lambda}(x,\cdot))_{0,-N/2}^{*}=(g_1,f_{\lambda}^L(x,\cdot))_{0,-N/2}^{*}+(g_2,f_{\lambda}^M(x,\cdot))_{0,-N/2}^{*}$$
and we see that $T_{\lambda,N}$ is defined on $L^{2}(\R^{\nu})$.

\vsp

Now define a mapping from $\mathcal{D}$ into $\D$ with kernel
$e_{\lambda}^{(\alpha)}(x,y) \in \D(\R^{\nu} \times
\R^{\nu})$, $T^{(\alpha)}_{\lambda}(x)f =
D^{\alpha}E(\lambda)f(x)$ .
For the partial regularity, we note that for any test function
$\varphi$, if $\psi=(1 - \Delta_{y''})^{-N} \varphi$, then using
(\ref{o0}), we get
$$ \sup_K \mid D^{\alpha}E_{\lambda}*''\psi f \mid \leq C_{N,{\alpha}}
\parallel \varphi \parallel_{L^1} \parallel f \parallel_H$$
the constant may depend on ${\lambda_1},{\lambda_2}$, but it is not
dependent on the choice of $\varphi$. For the partially regularized
resolution, we have the following exponential estimate,
\begin{equation}
 \mid T^{(\alpha)}_{\lambda}(x)f \mid \leq
C(x,\alpha,\delta) \exp (-\kappa \mid \lambda \mid^{b}) \otimes \delta_{x''} \parallel f
\parallel_H \label{Riesz}
\end{equation}
for $x$ in $\R^{\nu}$, for $f \in H$ and $\lambda$
negative. The constant is again dependent on the way we mollify. Further $\parallel e_{\lambda}^{(\alpha)}(x,\cdot) \parallel \leq C(x,\alpha,\delta)
\exp(-\kappa \mid \lambda \mid^{b})\otimes \delta_{x''}$, for $x$ in
$\R^{\nu}$. Just as in connection with the representation
(\ref{Ni}), we note that the estimate implies $ \parallel
e_{\lambda}^{(\alpha)}(x',\varphi_{\delta},y',\psi_{\delta}) \parallel
\leq$ $
C(x,\alpha,\delta)\exp(-\kappa \mid \lambda \mid^{b})C_{\varphi \otimes \psi}$.
The estimate (\ref{agmon}) still holds for the weighted spectral
function $e^{(\alpha)}_{\lambda}(x,\cdot)$.

\vspace{.5cm}

In the estimate following (\ref{Riesz}) above, if $\tilde{x}$ a point
sufficiently close to
$x$, we must have $\varphi(\tilde{x}'')$ arbitrarily close to
$\varphi(\widetilde{x}'')$, which means, for $\lambda$ sufficiently large, that $\{
e^{(\alpha)}_{\lambda}(\cdot,\varphi_{\delta}(\cdot),y',\psi_{\delta}(y'')) \}$ is equicontinuous in
$x$, for every $\alpha,\varphi$,$\psi$ and $y$. This holds for
every $x \in \R^{\nu}$ and through a sequence $y_n \rightarrow
y$, implies continuity in $(x,y)$ for the
family of functions.

\vsp

{ $\bf{Remark}$ (2): According to Proposition \ref{Prop6} $.2$ and (\ref{g-lambda}), we have that
$$ g_{\lambda}^{({\alpha'},{\beta'})}(x,y)=O(1)\exp(-\kappa \mid \lambda
\mid^b) \otimes ( 1+ q) \delta_{x''} +
O((R^{\#}_{\lambda})^{(\alpha',\beta')}(x,y))$$}
{ Had we in (\ref{Ni}) instead used a finite development of $g_{\lambda}$ such that, say
$\mbox{ deg }q=k$, then for appropriate test functions, we would get
the same estimates $B'_{\lambda}(x,y)=O(1)\exp(-\kappa\mid \lambda \mid^b)\otimes
\delta_{x''}$. Since $g_{\lambda}(x,y',\cdot) \in {\D}^F(\R^{\textit{m}})$,
we can let $k$ go to infinity, which would
allow an infinite development of $g_{\lambda}$ and we would still
get the same estimate in (\ref{Riesz}), however we will not develop
this approach any further.}

\vsp

Also
\begin{equation} \mid D^{\beta'}_ye^{(\alpha')}_{\lambda}(x',y') \mid \leq
C_{\beta'}\exp(-\kappa \mid \lambda \mid^{b}) \parallel
e^{(\alpha')}_{\lambda}(x',\cdot) \parallel'= \label{k1} \end{equation}
$$=O(1)\exp(-\kappa_1 \mid \lambda \mid^{b}) $$ on compact sets in
$\R^{\textit{n}}$ and for $x'$ in $\R^{\textit{n}}$.
This is interpreted as regularity in $\R^{\textit{n}}$, that is in
the "good" variable. The notation $\parallel \cdot \parallel'$, indicates
that the norm is taken over $\R^{\textit{n}}$.

\vsp

Outside $K$ ($ \subset {\Omega}$), the operator is $M(D')$
(we can assume $M(\xi') > 1$ for all $ \xi'$). The spectral kernel
then becomes, $f_{\lambda}(x,y)=\widetilde{f_{\lambda}}(x-y)$, where
$\widetilde{f_{\lambda}}(z)=k_{\lambda}(z') \otimes \delta_0(z'')$ and
$$ k_{\lambda}(z') = (2\pi)^{-n} \int_{M(\xi') < \lambda}
e^{i z' \cdot \xi '} d \xi ' {\label M}$$ and we know $k_{\lambda} \in
C^{\infty}(\R^{\textit{n}} \times \R^{\textit{n}})$. The following theorem and the preceding argument is close to \cite{Ni_72}
Theorem 1.

\newtheorem{Th1}[Prop6]{Theorem }
\begin{Th1}
  Assuming the operator partially formally self adjoint, we can
  for every real $\lambda$, implicitly define an element
  $e_{\lambda}(x,\cdot) \in H$, which mollified with test functions with support sufficiently close to the origin,
  is in $C^{\infty}$( $\R^{\nu} \times \R^{\nu}$).
  Further $e_{\lambda}^{(\alpha,0)}(x,\cdot)$ (distribution sense
  derivatives) is in
  $H$, for all $\alpha$,$x \in \R^{\nu}$
  and $$E(\lambda)u(x)= \int e_{\lambda}(x,y)u(y)dy \quad \text{for} \  u \in
  H \ \text{and} \ x \in \R^{\nu}$$
  For appropriate test functions, we have the estimates
  $D^{\beta}_y e_{\lambda}^{(\alpha,0)}(x,y)=O(1)\exp(-\kappa \mid \lambda
  \mid^{b}) \otimes \delta_{x''}$, uniformly on compact sets in $\R^{\textit{n}}  \times
  \R^{\textit{n}}$
  and $\parallel e_{\lambda}^{(\alpha,0)}(x',\varphi_{\delta},\cdot,\psi_{\delta}(\cdot)) \parallel = O(1)\exp(-\kappa \mid
  \lambda \mid^{b})$, uniformly on compact sets
  in $\R^{\nu}$, as $\lambda \rightarrow - \infty$, for x in $\R^{\nu}$.
\end{Th1}
Note that the theorem implies  the representation:
$$ E(\lambda)v(x)=\int f_{\lambda}*_y''{\U}_{-N}(x,y)v(y)dy \quad
\text{for} \quad v \in L^{2}, \quad x \in
\R^{\nu}$$

\subsection*{Asymptotic behavior of the spectral kernel}
\label{sec:asymp}
In this section we study the regularizations $\tilde{e}_{\lambda}(x,y)=e_{\lambda}*_x''\varphi*''_y\psi(x,y)$,
where $\varphi,\psi$ are appropriate test functions, of the spectral kernel corresponding to the partially formally hypoelliptic operator with variable coefficients $L$,
first with the assumption that $L \geq I$. The corresponding integral
operator is denoted $\tilde{E}_{\lambda}$.
We first note some immediate results (derived from \cite{Bedal}):

\vsp

     For an interval $\Delta$ in $\R$ and $f \in L^{2}(\R^{\nu})$, we have
     $$ \text{var}_{\Delta}
     \delta^{\alpha \beta}\tilde{e}_{\lambda}(x,y) \leq  \Big [
     \text{var}_{\Delta} \delta^{\alpha \alpha}\tilde{e}_{\lambda}(x,x)
     \text{var}_{\Delta} \delta^{\beta \beta}\tilde{e}_{\lambda}(y,y) \Big ]^{1/2}$$
     Here $\text{var}_{\Delta}$ denotes the variation over the
     interval $\Delta$ and
     $\delta^{\alpha \beta}$ derivation in x and y.
     Further $\text{var}_{\Delta}\delta^{\alpha}\tilde{e}_{\lambda}(f,y) \leq \mid \tilde{E}_{\Delta}f \mid
     \text{var}_{\Delta} \delta^{\alpha \alpha}\tilde{e}_{\lambda}(y,y)^{1/2} $ for $f \in L^{2}(\R^{\nu})$.

\vsp

    For test functions $\psi=\overline{\varphi}$, we have
     $ \delta^{\alpha \alpha} \tilde{e}_{\lambda}(x,x) \geq 0$ and it is
     a non-decreasing function of $\lambda$.
     The function $\delta^{\alpha \beta} \tilde{e}_{\lambda}(x,y)$ is locally of bounded variation as a function
     of $\lambda$, for $(x,y) \in \R^{\nu} \times \R^{\nu}$

\vsp

     For a realization of the operator, that is a self-adjoint spectral operator $\A_L$,
     we have that the resolvent operator G($\lambda) = (\A_L - \lambda I )^{-1}$
     is a bounded spectral operator on $L^{2}(\R^{\nu})$ and $$
     G(\lambda)=\int^{\infty}_{1} \frac{dE_{\mu}}{\mu - \lambda}
     \qquad \lambda < 1$$
     $\tilde{G}(\lambda)$ corresponding to the regularized spectral
     kernel, can also be represented as an integral operator with
     kernel $\tilde{G}_{\lambda}(x,y)$
     $$ \tilde{G}(\lambda)f(x) = \int \tilde{G}_{\lambda}(x,y)f(y)dy = \int \int (\mu - \lambda)^{-1} d\tilde{e}_{\mu}(x,y)f(y)dy$$
     where according to Stieltjes formula, $\tilde{E}_{\mu}f(x) = \int \tilde{e}_{\mu}(x,y)f(y)dy$ and $f \in L^{2}(\R^{\nu})$

\vsp

     Assuming we can prove an a priori estimate
     $$ \sup_{x,y \in K} \mid \tilde{e}_{\lambda}^{(\alpha,\beta)}(x,y) \mid \leq C_{K,\alpha,\beta}(1 + \lambda)^c \qquad \lambda \geq 1$$
     then with condition the constant $c$ is $<1$, we have that $\tilde{G}_{\lambda}(x,y)$
     is continuous on $\R^{\nu} \times \R^{\nu}$.

\vsp

     Proof of the estimate: Assume with notation as in section \ref{sec:Reg},
     that $\varphi \in H(0,\lambda)$.
     According to the inequality (\ref{o1}) in that section, we have for $kr > n/2 + m$ ($\mid \beta \mid \leq m$)
     \begin{equation} \sup_{y \in K} \mid D_y^{\beta} \widetilde{\varphi}(x,y) \mid \leq C_{K,\beta}( \parallel P^r \varphi \parallel_{H}
     + \parallel \varphi \parallel_{H}) \qquad x \in \R^{\nu} \label{a-priori} \end{equation}
     The inequality for a $\lambda$-approximative eigenvector to
     $P^r$ and $0$, where $r$ is assumed odd, gives
     $$ \parallel P^r \varphi \parallel_H \leq \lambda \parallel \varphi \parallel_H$$
     Using the relation $E(\lambda)=E_r(\lambda^r)$ and (\ref{a-priori}),
     we get for the resolution corresponding to the non-iterated operator,
     $$ \sup_{y \in K} \mid D^{\beta}_y \widetilde{\varphi}(x,y) \mid \leq C_{K,\beta}(1+ \lambda^{1/r})\parallel \varphi \parallel_H
     \leq C'_{K,\beta}(1 + \lambda)^{1/r} \parallel \varphi \parallel_H
     \quad x \in \R^{\nu}$$
     Finally, again using chapter \ref{sec:Reg}
     $$ \sup_{x,y \in K} \mid D_x^{\alpha}\tilde{E}(\lambda) D_y^{\beta} \widetilde{\varphi}
     \mid \leq C_{K,\alpha,\beta}(1+ \lambda)^{1/r} \parallel \varphi \parallel_H$$
     and by choosing $\widetilde{\varphi}$ as the regularized spectral function $\tilde{e}_{\lambda}(x,y)$,
     knowing that this function is in $L^{2}$, the estimate follows.$\Box$

\subsection*{ Estimate for the Green kernel }
     Using partial integration, we can write for $\mid \alpha' + \beta'
     \mid \leq M$ , for all $\alpha'',\beta'' $ and for
     $f,g \in C^{\infty}_0$
     \begin{equation} ( \tilde{G}_{\lambda}^{(\alpha,\beta)}f,g ) =
       \int^{\infty}_1 (\mu - \lambda)^{-1}
       d(\tilde{E}^{(\alpha,\beta)}_{\mu}f,g) \qquad \lambda < 1
       \label{Green} \end{equation}
     For the Green kernel corresponding to the operator with
     constant coefficients $P^x(D_{y''}) - \lambda$ {\small
     \begin{equation} \tilde{G}^{(\alpha,\alpha)}_{x,\lambda}(x,x) = \int^{\infty}_1 (\mu - \lambda)^{-1}
     d \tilde{e}_{x,\mu}^{(\alpha,\alpha)}(x,x)=\Big[ (iD)^{\alpha''}
     \varphi \Big]
     \Big[ (iD)^{\alpha''} \psi \Big] \int \frac{{\xi'}^{2
         \alpha'} d \xi'}{{\mbox Re } P^x(\xi')- \lambda} \label{gx}
   \end{equation} \par}
     where the integral is finite, for all $\alpha''$ and $ s > n + 2
     \mid \alpha' \mid$ and where we have used that $e_{x,\lambda}$ must
     be a tensor product with $\delta_{x''}$.

\vsp

A modification of Nilssons's article (cf.\cite{Ni_80}), gives that we can
use an estimate of the regularized fundamental solution
$\widetilde{g}_{\lambda}$, to produce an estimate of the difference between the resolvent operators
corresponding to the operator ${\mbox Re } P^x - \lambda$ and the variable coefficients operator respectively,
in terms of the first operator. More precisely
\begin{equation} \tilde{G}^{(\alpha,\alpha)}_{x,\lambda}(x,x) - \tilde{G}^{(\alpha,\alpha)}_{\lambda}(x,x) =
O(1)\mid \lambda \mid^{-c}\tilde{G}_{x,\lambda}^{(\alpha,\alpha)}(x,x) \quad \lambda \rightarrow - \infty \label{t0} \end{equation}
The left side in (\ref{t0}) can be written $\int^{\infty}_1 (t - \lambda )^{-1} d\sigma(t)$, where
$ \sigma(t) = \tilde{e}^{(\alpha,\alpha)}_{x,t}(x,x) - \tilde{e}^{(\alpha,\alpha)}_{t}(x,x)$,
which is a monotone non-decreasing function of $t$.

\subsection*{ Estimate for the spectral kernel }
Tauberian theory
applied to $\sigma(\mu)$, leads to an estimate of the regularized spectral function corresponding
to the partially formally hypoelliptic operator, in terms of the regularized spectral function
corresponding to the operator ${\mbox Re } P^x - \lambda$.
Note that
$$ \tilde{e}_{x,\lambda}^{(\alpha,\alpha)}(x,x) = \Big[ (iD)^{\alpha''} \varphi \Big]
\Big[ (iD)^{\alpha''} \psi \Big]\int_{ {\mbox Re } P^x(\xi') < \lambda} {\xi'}^{2 \alpha'} d \xi'$$
where $\lambda \in \R$. Assuming $\psi
= \overline{\varphi}$, we can prove (cf.\cite{Ni_80}),
\newtheorem{Lem5}{Lemma}[section]
\begin{Lem5} \label{Lem5}
There is a complex constant C, a rational number a and an integer t,
$0 \leq t \leq n - 1$, such that
$$ \tilde{e}_{x,\lambda}^{(\alpha,\alpha)}(x,x) =
C(1+o(1)){\lambda}^a(\text{log} \lambda)^t \quad \lambda \rightarrow \infty$$
Further $\tilde{e}_{x,\lambda}^{(\alpha,\alpha)}(x,x)$ is infinitely differentiable for
large $\lambda$ and $$ \frac{d \tilde{e}_{x,\lambda}^{(\alpha,\alpha)}(x,x)}{d \lambda} = o(1) \lambda^{a-1} (\text{log} \lambda )^t$$
as $\lambda \rightarrow \infty$
\end{Lem5}
\newtheorem{Lem6}[Lem5]{Lemma }
\begin{Lem6}
For a positive constant $c$ as in (\ref{t0}), we have
$$\sup_{\lambda \leq \mu \leq \lambda + \lambda/(c log \lambda)}
\int^{\mu}_{\lambda} d\sigma(w) =
o(1){\lambda}^a( log \lambda)^{r-1}$$
\end{Lem6}
Proof: (cf.\cite{Tsu}) Immediately from Lemma \ref{Lem5}
$$\frac{d}{dt} \tilde{e}_{x,t}^{(\alpha,\alpha)}(x,x) = o(1)t^{a-1}(\log{t})^r \quad t \rightarrow \infty$$
The properties of $\sigma(t)$ imply $\int_t^{\mu} d \tilde{e}_t^{(\alpha,\alpha)}(x,x) \geq 0$.
If $t \leq \mu \leq t + t/c \log{t} $, we get $ \int^{\mu}_t d \sigma(t) \leq
\int_t^{\mu} d \tilde{e}_{x,t}^{(\alpha,\alpha)}(x,x) \leq C\int^{\mu}_t t^{a-1} (\log{t} )^r dt
\leq C'(t/\log{t})t^{a-1}(\log{t} )^r = C' t^a (\log{t} )^{r-1} \Box$

\vsp

Immediately from (\ref{gx}), (\ref{t0}) and Lemma \ref{Lem5}, we get
$$ \int^{\infty}_1 \frac{ d \sigma(\mu)}{  \mu + \lambda } =
O(1) \mid \lambda \mid^{-c}\lambda^{a-1}(\textrm{log} \lambda )^r \quad \lambda \rightarrow \infty$$
That is if we use the Stieltjes transform, we see that
$$ \int^{\infty}_1 \frac{d \widetilde{e}_{x,\mu}^{(\alpha,\alpha)}(x,x)}{\mu + \lambda}
= \mathcal{L}(\mathcal{L}(\frac{\textit{d}}{\textit{d} \mu}\widetilde{\textit{e}}_{\textit{x},\mu}^{(\alpha,\alpha)}(\textit{x,x})))$$
and a well known result in Abelian theory of the Laplace transform,
gives that $\overline{\lim}_{\mu \rightarrow \infty}
\mathcal{L}(\frac{\textit{d}}{\textit{d} \mu} \widetilde{\textit{e}}_{\textit{x},\mu}^{(\alpha,\alpha)}) \leq
\overline{\lim}_{\mu \rightarrow
  \infty}(\frac{\textit{d}}{\textit{d} \mu} \widetilde{\textit{e}}_{\textit{x},\mu}^{(\alpha,\alpha)})$
and the result follows since $\mu \leq \lambda$.
This means that the conditions in Ganelius Tauberian theorem (cf.\cite{Gan})
are satisfied
for the operator $L \geq I$. However this restriction can be discarded if we study
instead the operator $ \A_L^r + kI$, for $k,r$ sufficiently
large, both numbers are
assumed $ \geq 0$, $r$ is an even integer and $k$ is real, since this
translation does not effect the asymptotic behavior.
For details we refer to \cite{Ni_72}.

\vsp

\newtheorem{Th2}[Lem5]{Theorem }
\begin{Th2} \label{Th2}
For every multi-index $\alpha$ and every $x \in W$, $W$ a compact set
where the operator $L_{\lambda}$ is assumed partially formally
hypoelliptic, we have
$$ \tilde{e}_{\lambda}^{(\alpha,\alpha)}(x,x) = (1+ o(1)(log \lambda)^{-1})
\tilde{e}^{(\alpha,\alpha)}_{x,\lambda}(x,x) \qquad \lambda
\rightarrow \infty$$
\end{Th2}

\cite{Zuily}

\bibliographystyle{amsplain}
\bibliography{ref}
\end{flushleft}
\end{document}